\documentclass{article}
\usepackage{mathrsfs}
\usepackage{color}
\usepackage{epsfig, graphicx}
\usepackage{latexsym,amsfonts,amsbsy,amssymb}
\usepackage{amsmath,amsthm}
\usepackage{framed}
\usepackage{setspace}
\makeatletter % '@' is now a normal "letter" for TeX

\@addtoreset{equation}{section}
\makeatother % '@' is restored as a "non-letter" character for TeX
\allowdisplaybreaks % For long formula
\newtheorem{theorem}{Theorem}[section]
\newtheorem{lemma}{Lemma}[section]
\newtheorem{corollary}{Corollary}[section]
\newtheorem{remark}{Remark}[section]

\newtheorem{proposition}{Proposition}[section]
\newtheorem{assumption}{Assumption}[section]
\newcommand{\comm}[1]{{\color{black}#1}}   %red
\newcommand{\revise}[1]{{\color{black}#1}}  %blue
\begin{document}
\title{An Efficient Adaptive Finite Element Method for Eigenvalue
Problems\footnote{This work was supported in part
by the National Key Research and Development Program of China (2019YFA0709601), Beijing Natural Science Foundation (Z200003),
National Natural Science Foundations of China (NSFC 11771434), the National Center for Mathematics and Interdisciplinary Science, CAS.}}
\date{}
\author{
Qichen Hong\footnote{LSEC, ICMSEC, Academy of Mathematics and Systems Science,
Chinese Academy of Sciences, Beijing 100190, P.R. China (hongqichen@lsec.cc.ac.cn)}, \ \
Hehu Xie\footnote{LSEC, NCMIS, Institute
of Computational Mathematics, Academy of Mathematics and Systems
Science, Chinese Academy of Sciences, Beijing 100190,
China,  and School of Mathematical Sciences, University
of Chinese Academy of Sciences, Beijing, 100049, China (hhxie@lsec.cc.ac.cn)}\ \ \ and \ \
Fei Xu\footnote{LSEC, ICMSEC, Academy of Mathematics and Systems Science,
Chinese Academy of Sciences, Beijing 100190, P.R. China (xufei@lsec.cc.ac.cn)}}
\maketitle
%--------------------------------------------------------------------------------
\begin{abstract}
The aim of this paper is to propose an efficient adaptive finite element method for eigenvalue
problems based on the multilevel correction scheme and inverse power method.
This method involves solving associated boundary value problems on each adaptive
partitions and very low dimensional eigenvalue problems on some special meshes which are
controlled by the proposed algorithm. Since we
Hence the efficiency of solving eigenvalue problems can be improved to be similar to the
adaptive finite element method for the associated boundary value problems.
The convergence and optimal complexity is theoretically verified and numerically demonstrated.
\end{abstract}

{\bf Keywords.}\ Eigenvalue problem, multilevel correction method, inverse power,
adaptive finite element method, convergence, optimality

{\bf AMS Subject Classification:} 65F15, 65N15, 65N25, 65N30, 65N50
%============================================================================================

\section{Introduction}
The finite element method is one of the widely used discretization schemes for solving
 eigenvalue problems.
The adaptive finite element method (AFEM) is a meaningful approach which can generate a
 sequence of optimal triangulations by refining those elements where the errors, as the
 local error estimators indicate, are relatively large. The
AFEM is really an effective way to make efficient use of given computational resources.
 Since Babu\v{s}ka and Rheinboldt
\cite{BabuskaRheinboldt}, the AFEM has been an active topic, many researchers are attracted
to study the AFEM (see, e.g., \cite{ArnoldMukherjeePouly,BabuskaRheinboldt,BabuskaVogelius,
BrennerScott,MorinNochettoSiebert_2002,Nochetto,WuChen} and the references cited therein)
 in the past 40 years. So far, the convergence and optimality of the AFEM for
  boundary value problems have been obtained and
understood well (see, e.g., \cite{BinevDahmenDeVore,CasconKreuzerNochettoSiebert,Dofler,
DoflerWilderotter,MekchayNochetto,MorinNochettoSiebert_2000,MorinNochettoSiebert_2002,
Stevenson_2007,Stevson_2008} and the references cited therein).

It is well known that the eigenvalue problem is one of the fundamental problems in computational
mathematics and large scale eigenvalue problems always occur in discipline of sciences and engineering
such as materials science, quantum chemistry or physics, structure mechanics, biological system,
data and information fields, etc. However, it is always a very difficult task to solve high-dimensional eigenvalue
problems which come from practical physical and chemistry sciences, and there is a strong demand
by engineers and scientists for efficient eigenvalue solvers.
Besides for the boundary value problems, the AFEM is also a very useful and efficient
way for solving eigenvalue problems (see, e.g.,
\cite{BeckerRannacher, DuranPadraRodrguez,HeuvelineRannacher,Larson,MaoShenZhou,Verfurth}).
The AFEM for eigenvalue problems has been analyzed in some papers
(see, e.g., \cite{DaiXuZhou,GianiGraham,HeuvelineRannacher} and the references cited therein).
Especially, \cite{DaiXuZhou} give an elaborate analysis of the convergence and optimality
for the adaptive finite element eigenvalue computation based on the methods and results
in \cite{CasconKreuzerNochettoSiebert}. In \cite{GianiGraham}, authors also give the
analysis of the convergence for the eigenvalue problems by the AFEM. The optimality in AFEM only
means the scale of the discretization is optimal. But the computing efficiency of AFEM does not arrive the optimality.

\comm{In order to improve the efficiency of the AFEM for eigenvalue problems,
the purpose of this paper is to propose and analyze a type of AFEM
to solve the eigenvalue problems based on the adaptive refinement technique and the recent work on
the multilevel correction method \cite{LinXie,Xie_IMA,Xie_JCP,Xie_Nonconforming,Jia_Xie_Xu_Xie,Xie_Xu_FullMulti}.
Compared with the standard AFEM which includes  solving eigenvalue problems on each refined mesh, we
only needs to solve the associated linear boundary value problem on each refined mesh and
some very low dimensional eigenvalue problems at some special adaptive steps which is controlled by the AFEM proposed in this paper.
Furthermore, the dimension or scale of the included in the new AFEM is fixed all through the adaptive refine process.
%Instead of solving eigenvalue problems on each adaptive mesh, we only need to solve
%boundary value problems on each adaptive mesh and some very low dimensional eigenvalue problems
%on special meshes which is controlled by the AFEM proposed in this paper.
Thus, in this new scheme, the cost of solving eigenvalue problems is almost the same as
solving the associated boundary value problems and the overall efficiency of eigenvalue solving can be improved.
Here, we also prove the convergence and quasi-optimal complexity of the new AFEM for the eigenvalue problems.
}

The rest of the paper is arranged as follows. In Section 2, we shall describe some basic
notation and the adaptive multilevel correction algorithm for the second order elliptic
eigenvalue problem. We then give the analysis of convergence and complexity of
the proposed AFEM in Section 3 and Section 4, respectively.
%In this paper, we will use the induction
%method to prove the convergence of the proposed AFEM.
%The preparation of the induction method is provided in Section 4.
%Then Section 5 is devoted to giving the convergence and optimality results for the
%AFEM proposed in this paper.
In Section 5, some numerical experiments are presented
to test the theoretical analysis. Finally, some concluding remarks are given in the
last section.

%--------------------------------------------------------------------------------------
\section{\revise{Multilevel correction adaptive finite element method}}
Let $\Omega\subset \mathbb{R}^d$ $(d\geq 1)$ denotes a polytopic bounded domain with Lipschitz continuous boundary.
In this paper, standard notation for Sobolev spaces $W^{s,p}(\Omega)$ and their associated norms
and seminorms  (see e.g., \cite{Adams,Ciarlet}) will be used. We denote $H^{s}(\Omega)=W^{s,2}(\Omega)$
and $H_0^{1}(\Omega)=\big\{v\in H^{1}(\Omega):v|_{\partial \Omega}=0\big\}$,
where $v|_{\partial \Omega}$ is understood in the sense of trace,
$\|v\|_{s,\Omega}=\|v\|_{s,2,\Omega}$ and $\|v\|_{0,\Omega}=\|v\|_{0,2,\Omega}$.
Throughout this paper, let $V:=H_0^{1}(\Omega)$. We consider the finite element
discretization on the shape regular family of nested conforming meshes $\{\mathcal{T}_k\}$
over $\Omega$: there exists a constant $\gamma^*$ such that
$$\frac{h_{T}}{\rho_T}\leq \gamma^*,\ \ \ \ \  \forall T\in\bigcup\limits_{k}\mathcal{T}_k,$$
where $h_T$ denotes the diameter of $T$ for each $T\in \mathcal{T}_k$, and $\rho_T$ is the
diameter of the biggest ball contained in $T$, $h_k:=\max\{h_T:T\in \mathcal{T}_k\}$.
In this paper, the notation $\mathcal{E}_k$ is used to denote the set of interior faces (edges or sides) of
$\mathcal{T}_k$.

%-----------------------------------------------------------------------------------------------------
\subsection{Preliminaries}

In this paper, we are concerned with the following second order elliptic eigenvalue problem
\begin{eqnarray*}\label{eigenproblem}
\left\{
\begin{array}{rcl}
Lu:=-\nabla (A\cdot \nabla u)+\varphi u&=&\lambda u\ \ \ {\rm in}\ \Omega,\\
u&=&0\ \ \ \ \ {\rm on}\ \partial\Omega,\\
\int_{\Omega}\big(A\nabla u\cdot\nabla u+\varphi u^2\big)d\Omega&=&1,
\end{array}
\right.
\end{eqnarray*}
where $A=(a_{ij})_{d\times d}$ is a symmetric positive definite matrix with
$a_{ij}\in W^{1,\infty}(\Omega)$\ $(i,j=1,\cdots,d)$, and $0\leq \varphi \in L^{\infty}(\Omega)$.

We first define a bounded bilinear form
\begin{eqnarray*}
a(u,v)=\int_{\Omega}\big(A\nabla u\cdot\nabla v+\varphi u v\big)d\Omega.
\end{eqnarray*}
From the properties of $A$ and $\varphi$, the bilinear form $a(\cdot,\cdot)$ is bounded over $V$
$$|a(w,v)|\leq C_{a}^2\|w\|_{1,\Omega}\|v\|_{1,\Omega},\ \ \ \forall w,v \in V,$$
and satisfies
\begin{eqnarray}\label{eqn_equality_anorm_1norm}
c_{a}\|w\|_{1,\Omega} \leq \|w\|_{a,\Omega} \leq C_{a} \|w\|_{1,\Omega},
\end{eqnarray}
where the energy norm $\|\cdot \|_{a,\Omega}$ is defined by $\| w \|_{a,\Omega}=\sqrt{a(w,w)} $, $c_{a}$ and $C_{a}$
are positive constants.
Then the corresponding variational form can be written as:
Find $(\lambda, u)\in \mathbb{R}\times V$ such that $\|u\|_{a,\Omega}=1$ and
\begin{eqnarray} \label{eigenvalue problem}
a(u,v)=\lambda(u,v),\ \ \ \forall v\in V.
\end{eqnarray}
As we know, the eigenvalue problem (\ref{eigenvalue problem}) has a countable sequence of real eigenvalues
$$0<\hat{\lambda}^1<\hat{\lambda}^2\leq\hat{\lambda}^3\leq\cdots$$
and corresponding orthogonal eigenfunctions
$$\hat{u}^1,\hat{u}^2,\hat{u}^3,\cdots, $$
which satisfy  $a(\hat{u}^i,\hat{u}^j)=\delta_{ij}$, $i,j=1,2,\cdots$.
Here we use $\hat{\lambda}^i$ and $\hat{u}^i$ to denote the $i$-th exact eigenvalue and eigenfunction, respectively.

Let $V_k\subset V$ be the corresponding family of nested finite element spaces  of continuous
piecewise polynomials over $\mathcal{T}_k$ of fixed degree $m\geq 1$, which vanish on the boundary of $\Omega$,
and are equipped with the same norm $\| \cdot \|_{a,\Omega}$ of space $V$.
The standard finite element discretization for (\ref{eigenvalue problem}) is:
Find $(\tilde{\tilde{\lambda}}_k,\tilde{\tilde{u}}_k)\in \mathbb{R}\times V_k$
such that $\|\tilde{\tilde{u}}_k\|_{a,\Omega}=1$ and
\begin{eqnarray}\label{diseig}
 a(\tilde{\tilde{u}}_k,v_k)=\tilde{\tilde{\lambda}}_k(\tilde{\tilde{u}}_k,v_k),\ \ \ \forall v_k\in V_k.
\end{eqnarray}
The eigenvalues of (\ref{diseig}) can also be ordered as an  increasing sequence
$$0<\tilde{\tilde{\lambda}}^1_k < \tilde{\tilde{\lambda}}^2_k
\leq\cdots\leq \tilde{\tilde{\lambda}}^{n_k}_k,\ \ \ n_k={\rm dim}V_k,$$
and the corresponding orthogonal eigenfunctions
$$\tilde{\tilde{u}}^{1}_k,\ \tilde{\tilde{u}}^{2}_k,\ \cdots,\ \tilde{\tilde{u}}^{n_k}_k $$
satisfying $a(\tilde{\tilde{u}}_k^i,\tilde{\tilde{u}}_k^j)=\delta_{ij}$, $i, j=1, 2, \cdots, n_k$.

Based on the finite element space $V_k$, we define the Galerkin projection $R_k:V\rightarrow V_k$ by
\begin{eqnarray*}\label{projection}
a(v-R_kv,v_k)=0,\ \ \ \forall v \in V,\ \forall v_k\in V_k.
\end{eqnarray*}
Then the following boundness holds
\begin{eqnarray*}\label{Proj_Bounded}
\|R_k v\|_{a,\Omega}\leq \|v\|_{a,\Omega},\ \ \ \forall v\in V.
\end{eqnarray*}
%------------------------------------------------------------------------------------------------------------
\begin{lemma}\label{Aubin-Nitsche}(\cite{BrennerScott,Ciarlet})
The duality argument leads to the following inequality
\begin{eqnarray*}
\|(I-R_k)v\|_{0,\Omega} \leq C_{an} \eta_a(V_k)\|(I-R_k)v\|_{a,\Omega},
\end{eqnarray*}
where $C_{an}$ is a positive constant and the quantity $\eta_a(V_k)$ is defined as follows:
\begin{eqnarray}\label{definition_quantity_eta}
\eta_a(V_k)=\sup\limits_{f\in L^2(\Omega),\|f\|_{0,\Omega}=1}\inf\limits_{v_k\in V_k}\|L^{-1}f-v_k\|_{a,\Omega}.
\end{eqnarray}
\end{lemma}
Let $K:L^2(\Omega)\rightarrow V$ be the operator defined by
\begin{eqnarray*}\label{postK}
a(Kw,v)=(w,v),\ \ \ \forall w,  v\in V.
\end{eqnarray*}
Then the eigenvalue problems (\ref{eigenvalue problem}) and (\ref{diseig}) can be written as
\begin{eqnarray}\label{Operator_Eigenvalue_Problem}
 u=\lambda Ku,\ \ \ \ \tilde{\tilde{u}}_k=\tilde{\tilde{\lambda}}_kR_kK\tilde{\tilde{u}}_k.
\end{eqnarray}
For any $v\in V$, since $a(Kv,Kv)=(v,Kv)$ and from (\ref{eqn_equality_anorm_1norm}), we have
\begin{eqnarray}\label{eqn_K_norm}
\|Kv\|_{a,\Omega} \leq \frac{\|v\|_{0,\Omega}}{c_a} \leq \frac{\|v\|_{a,\Omega}}{c_a^2}.
\end{eqnarray}

For the aim of error estimate, we define
\begin{eqnarray*}
M(\hat{\lambda}^i)=\big\{v\in V:v \text{ is an eigenfunction of (\ref{eigenvalue problem}) }
\text{corresponding to the eigenvalue }\hat{\lambda}^i\big\}
\end{eqnarray*}
and the quantity
\begin{eqnarray}\label{definition_quantity_delta}
\delta_{W}(\hat{\lambda}^i)
=\sup\limits_{v\in M(\hat{\lambda}^i),\|v\|_{a,\Omega}=1}\inf\limits_{w\in W}\|v-w\|_{a,\Omega},
\end{eqnarray}
where $W$ is a finite dimensional space.

From \cite{BabuskaOsborn_1989,BabuskaOsborn_1991,Chatelin,Ciarlet}, it is known that
$\eta_a(V_k)\rightarrow 0$, $\delta_{V_k}(\hat{\lambda}^i)\rightarrow 0$ as $h_k\rightarrow 0 $
and the following error estimates of finite element method for eigenvalue problems hold.
\begin{lemma}(\cite{BabuskaOsborn_1989,BabuskaOsborn_1991,Chatelin,Ciarlet})\label{lemma_error_estimate_eigenproblem}
Let $(\tilde{\tilde{\lambda}}_k^i,\tilde{\tilde{u}}_k^i)\in \mathbb{R}\times V_k$ be the solution of
(\ref{diseig}) for $1\leq i\leq n_k$. Then there exist an exact eigenpair $(\hat{\lambda}^i,\hat{u}^i)$
of (\ref{eigenvalue problem}) and constants $\bar{C}_{ea}$, $\bar{C}_{e0}$ and $\bar{C}_{e\lambda}$ such that
\begin{eqnarray*}
\|\tilde{\tilde{u}}_k^i-\hat{u}^i\|_{a,\Omega}
&\leq& \bar{C}_{ea}\delta_{V_k}(\hat{\lambda}^i),\label{Error_u_u_k}\\
\|\tilde{\tilde{u}}_k^i-\hat{u}^i\|_{0,\Omega}&\leq& \bar{C}_{e0}
\eta_a(V_k)\|\tilde{\tilde{u}}_k^i-\hat{u}^i\|_{a,\Omega},\label{Error_u_u_k_-1}\\
|\tilde{\tilde{\lambda}}_k^i-\hat{\lambda}^i|
&\leq& \bar{C}_{e\lambda}\|\tilde{\tilde{u}}_k^i-\hat{u}^i\|_{a,\Omega}^2. \label{Error_lambda_k}
\end{eqnarray*}
\end{lemma}

%------------------------------------------------------------------------------------------------------------
\subsection{Adaptive multilevel correction algorithm}
Now we follow the classic routine to define the a posteriori error estimator for boundary value
problem with $f\in L^2(\Omega)$ as the right hand side term. \revise{For any $f\in L^2(\Omega)$ and $v_k\in V_k$, let us}
define the element residual $\mathcal{R}_T(f,u_k)$ and the jump residual $\mathcal{J}_E(u_k)$ by
\begin{eqnarray*}
\mathcal{R}_T(f,v_k):=f-Lv_k=f
 +\nabla\cdot(A\nabla v_k)-\varphi v_k \ \ \text{in}\ T\in \mathcal{T}_k,\nonumber\\
\hskip-0.4cm\mathcal{J}_E(v_k):=-A\nabla v_k^+\cdot\nu^+-A\nabla v_k^-\cdot\nu^-
:=[[A\nabla v_k]]_E\cdot \nu_E \ \ {\rm on} \  E\in \mathcal{E}_k,\label{Edge_Jump}
\end{eqnarray*}
where $E$ is the common side of elements $T^{+}$ and $T^-$ with outward normals
$\nu^+$ and $\nu^-$, $\nu_E=\nu^-$. Then we can define the local error indicator
$\eta_k(f,v_k;T)$ for the element $T\in\mathcal{T}_k$ by
\begin{eqnarray}\label{definition_errorestimate}
\eta_k^2(f,v_k;T):=h_T^2\|\mathcal{R}_T(f,v_k)\|_{0,T}^2
+\sum\limits_{E\in \mathcal{E}_k,E\subset \partial T}h_E\|\mathcal{J}_E(v_k)\|^2_{0,E},
\end{eqnarray}
and for a submesh $\mathcal{T}'\subset\mathcal{T}_k$ by
\begin{eqnarray*}
\eta_k(f,v_k;\mathcal{T}'):=\left(\sum_{T\in \mathcal{T}'}\eta_k^2(f,v_k;T)\right)^{1/2}.
\end{eqnarray*}
Thus $\eta_k(f,v_k;\mathcal{T}_k)$ denotes the error estimator of finite element approximation $v_k$ with respect to $\mathcal{T}_k$.

For $f\in L^2(\Omega)$, we define the data oscillation as
\begin{eqnarray*}
{\rm osc}(f;\mathcal{T}_k) := \left(\sum_{T\in \mathcal{T}_{k}}\|h_T(f-P_Tf)\|_{0,T}^2\right)^{\frac{1}{2}},
\end{eqnarray*}
where $P_T$ is the $L^2$-projection operator to polynomials of some degree on $T$.
It is obvious that the following inequality holds
\begin{eqnarray*}
{\rm osc}(f-Lv_k;\mathcal{T}_k) \leq \eta_k(f,v_k;\mathcal{T}_k).
\end{eqnarray*}

For convenience, we use the notation $\eta(f;\mathcal{T}_k):=\eta_k(f,R_kKf;\mathcal{T}_k)$.
There exist the following reliability and efficiency for the a posterior error estimator $\eta(f;\mathcal{T}_k)$
(see, e.g., \cite{MekchayNochetto,MorinNochettoSiebert_2002,Verfurth}):
%------------------------------------------------------------------------------------------------------------
\begin{lemma}(\cite{MekchayNochetto,MorinNochettoSiebert_2002,Verfurth})\label{lemma_upper_lower_bound_boundaryvalue}
For $f\in L^2(\Omega)$, there exist a constant $C_{\rm up}$,
solely depending on regularity constant $\gamma^*$ and coercivity
constant $c_a$ in (\ref{eqn_equality_anorm_1norm}), such that
\begin{eqnarray*}
\|(I-R_k)Kf\|_{a,\Omega}\leq C_{\rm up}\eta(f;\mathcal{T}_k),
\end{eqnarray*}
and a constant $C_{\rm low}$, solely depending on regularity constant $\gamma^*$ and continuity
constant $C_a$ in (\ref{eqn_equality_anorm_1norm}), such that
\begin{eqnarray}\label{Lower_Bound_Posteriori_Error}
C_{\rm low}\eta(f;\mathcal{T}_k)\leq \|(I-R_k)Kf\|_{a,\Omega}+{\rm osc}(f-LR_kKf;\mathcal{T}_k).
\end{eqnarray}
\end{lemma}
Before introducing our AFEM, we first introduce some modules for preparation:
\begin{itemize}
\item $(\mu,w)=\textsf{ESOLVE}(W)$:
Solve the eigenvalue problem (\ref{diseig})  in the finite element space $W$ and output the
discrete eigenpair $(\mu,w)\in \mathbb{R}\times W$.
\item $w=\textsf{LSOLVE}(f,W)$:
Solve the linear boundary value problem in the finite element space $W$ with the right hand side term $f$,
namely, the output $w\in W$ satisfies the following boundary value problem
\begin{eqnarray*}
a(w,v)=(f,v),\ \ \ \ \forall v\in W.
\end{eqnarray*}
\item $\{\eta_k(f,v_k; T)\}_{T\in\mathcal{T}_k}=\textsf{ESTIMATE}(f,v_k,\mathcal{T}_k)$:
Compute the error indicator  on each element $T\in\mathcal T_k$ as (\ref{definition_errorestimate}).
\item $\mathcal{M}_k=\textsf{MARK}(\theta,\{\eta_k(f,v_k; T)\}_{T\in\mathcal{T}_k},\mathcal{T}_k)$:
Construct a minimal subset $\mathcal{M}_k$ from $\mathcal{T}_k$ by
selecting some elements in $\mathcal{T}_k$ such that
\begin{eqnarray*}
\eta_k(f,v_k;\mathcal{M}_k)\geq\theta\eta_k(f,v_k;\mathcal{T}_k)
\end{eqnarray*}
and mark all the elements in $\mathcal{M}_k$.
\item $(\mathcal{T}_{k+1},V_{k+1})=\textsf{REFINE}(\mathcal{M}_k,\mathcal{T}_k)$:
Output a conforming refinement $\mathcal{T}_{k+1}$ of $\mathcal{T}_k$ where at least all
elements of $\mathcal{M}_k$ are refined and construct the finite element space $V_{k+1}$
over $\mathcal{T}_{k+1}$.
\end{itemize}
Then we present a type of AFEM to compute the eigenvalue problem in the multilevel
correction framework which is the main contribution of this paper.
%-------------------------------------------------------------------------------------------------------
\begin{framed}
\begin{center}
\textbf{Adaptive Algorithm $C$ }
\end{center}
\begin{spacing}{1.5}
\noindent
\begin{itemize}
\item Given parameters $0<\theta_1<1$ and $0<\theta_2<1$, a coarse mesh $\mathcal{T}_H$ with mesh size $H$ and
construct the finite element space $V_H$.
\item Refine the mesh $\mathcal{T}_H$ to obtain an initial mesh $\mathcal{T}_{1}$ and the finite element space $V_1$ by the
regular way.
\item Solve an eigenvalue problem: $(\lambda_1,u_1)=\textsf{ESOLVE}(V_1)$.
\item Set $\bar{f}_2^{(0)}=\bar{f}_1=\tilde u_1=u_1$, $\bar{u}_1=\frac{u_1}{\lambda_1}$, $n_{1}=1$, $k=2$, $\ell=1$, $j=0$.
\item Compute $\{\eta_1(\bar{f}_1,\bar{u}_{1}; T)\}_{T\in\mathcal{T}_1}=\textsf{ESTIMATE}(\bar{f}_1,\bar{u}_{1},\mathcal{T}_1)$ and set $\eta_1=\eta_1(\bar{f}_1,\bar{u}_{1};\mathcal{T}_1)$.
\end{itemize}
Do the following iteration:
\begin{enumerate}
\item $\mathcal{M}_{k-1}=\textsf{MARK}(\theta_1,\{\eta_{k-1}(\bar{f}_{k-1},\bar{u}_{k-1}; T)\}_{T\in\mathcal{T}_{k-1}},\mathcal{T}_{k-1})$;
\item $(\mathcal{T}_{k},V_{k})=\textsf{REFINE}(\mathcal{M}_{k-1},\mathcal{T}_{k-1})$;
\item Linear solving: $\bar{u}_k^{(j)}=\textsf{LSOLVE}(\bar{f}_k^{(j)},V_{k})$ and compute
\begin{eqnarray*}
u_k^{(j)}=\frac{\bar{u}_k^{(j)}}{\|\bar{u}_k^{(j)}\|_{a,\Omega}},\ \ \ \ \
\lambda_k^{(j)} =\frac{a(u_k^{(j)},u_k^{(j)})}{(u_k^{(j)},u_k^{(j)})};
\end{eqnarray*}
\item $\{\eta_k(\bar{f}_k^{(j)},\bar{u}_{k}^{(j)}; T)\}_{T\in\mathcal{T}_k}=\textsf{ESTIMATE}(\bar{f}_k^{(j)},\bar{u}_{k}^{(j)},\mathcal{T}_k)$;
\item \begin{itemize}
\item If $\eta_k(\bar{f}_k^{(j)},\bar{u}_{k}^{(j)}; \mathcal{T}_k) \leq \theta_2^{j+1} \eta_{\ell}$, then solve $(\tilde{\lambda}_{k}^{(j)},\tilde{u}_{k}^{(j)})=\textsf{ESOLVE}(V_H+{\rm span}\{u_{k}^{(j)}\})$, set $\bar{f}_{k}^{(j+1)}=\tilde{u}_k^{(j)}$, $j=j+1$ and go to step 3;
\item Else, set $\bar{f}_k=\bar{f}_k^{(j)}$, $\bar{u}_k=\bar{u}_k^{(j)}$, $u_k = u_k^{(j)}$ and $\lambda_k = \lambda_k^{(j)}$;
\end{itemize}
\item \begin{itemize}
\item If $j>0$, solve an eigenvalue problem
$(\tilde{\lambda}_{k},\tilde{u}_{k})=\textsf{ESOLVE}(V_H+{\rm span}\{u_{k}\})$ and set $\bar{f}_{k+1}^{(0)}=\tilde{u}_k$; Then set $\eta_{\ell+1} = \eta_k(\bar{f}_k^{(j)},\bar{u}_{k}^{(j)}; \mathcal{T}_k)$, $n_{\ell+1}=k$, $\ell=\ell+1$, $j=0$;
\item Else, set $\bar{f}_{k+1}^{(0)}=u_k$;
\end{itemize}
\item Let $k=k+1$ and go to step 1.
\end{enumerate}
\end{spacing}
\end{framed}
%---------------------------------------------------------------------------------------------------------
\begin{remark}
Here we use the iterative or recursive
bisection (see, e.g., \cite{Maubach,Traxler}) of elements with the minimal refinement condition in the
procedure \textsf{REFINE}. The marking strategy adopted in {\bf Adaptive Algorithm $C$}
was introduced by D\"{o}rfler \cite{Dofler} and Morin et al. \cite{MorinNochettoSiebert_2002}.

Different from the standard AFEM for eigenvalue problems, {\bf Adaptive Algorithm $C$} has no requirement to
solve the eigenvalue problems on the adaptively refined meshes $\mathcal T_k$ which can improve the efficiency since
eigenvalue solving need much more computations than solving the associated linear boundary value problem.
This point is the main contribution of this paper.
\end{remark}
%-----------------------------------------------------------------------------------------------------------
In {\bf Adaptive Algorithm $C$}, we denote by index $k$ the order of nested mesh, index $j$ the order of iteration on a fixed mesh.
Besides, we use $n_{\ell}$ ($\ell>1$) to denote the order of the mesh on which the decision condition
$\eta_{n_{\ell}}(\bar{f}_{n_{\ell}}^{(0)},\bar{u}_{n_{\ell}}^{(0)};\mathcal{T}_{n_{\ell}}) \leq \theta_2 \eta_{\ell-1}$ is satisfied, where reference estimate $\eta_{\ell-1} =\eta_{n_{\ell-1}}(\bar{f}_{n_{\ell-1}},\bar{u}_{n_{\ell-1}};\mathcal{T}_{n_{\ell-1}}) $ is defined recursively, and $\ell$ is introduced to denote the order of subset $\{\mathcal{T}_{n_{\ell}}\}_{\ell\in \mathbb{N}} \subset \{\mathcal{T}_k\}_{k\in \mathbb{N}}$.
For easier to understanding,  we visualize {\bf Adaptive Algorithm $C$} with two flow charts.

The first flow chart translates {\bf Adaptive Algorithm $C$} to the figurative language.

The second flow chart shows the overall behavior of {\bf Adaptive Algorithm $C$}.
The processes between two dotted lines are over meshes from $\mathcal{T}_{n_{\ell}+1}$ to $\mathcal{T}_{n_{\ell+1}}$.

In this paper, $(\tilde{\lambda}_{n_{\ell}},\tilde{u}_{n_{\ell}})$ is the eigenpair approximation
by the eigenvalue solving module $\textsf{ESOLVE}(V_1)$ for $\ell=1$ and
$\textsf{ESOLVE}(V_H+{\rm span}\{u_{n_{\ell}}\})$ for $\ell\geq 2$, $(\tilde{\lambda}_{n_{\ell}}^{(j)},\tilde{u}_{n_{\ell}}^{(j)})$ is the eigenpair approximation by module $\textsf{ESOLVE}(V_H+{\rm span}\{u_{n_{\ell}}^{(j)}\})$ for $\ell\geq 2$ and $j\geq 0$, $(\lambda_k,u_k)$ denotes the eigenpair
approximation by the linear solving.

From {\bf Adaptive Algorithm $C$}, over each mesh $\mathcal{T}_k$ ($k\in \mathbb{N}$) we solve a linear
boundary value problem with $\bar{f}_k$ as the right hand side term, and then refine the mesh based
on the corresponding error indicator. The following theorem
can be derived directly from the convergence of \revise{standard AFEM for boundary value problems}, which
has been proved by Cascon et al \cite{CasconKreuzerNochettoSiebert}.
%----------------------------------------------------------------------------------------------------------
\begin{theorem}(\cite{CasconKreuzerNochettoSiebert})\label{Adaptive_Convergence_Source_Theorem}
 Let $\{\bar{f}_k\}_{k\in \mathbb{N}}$ and $\{\mathcal{T}_{k}\}_{k\in\mathbb{N}}$ be produced by
 {\bf Adaptive Algorithm $C$}.
Then, there exist constants $\gamma>0$ and $\hat{\alpha}\in (0,1)$, depending on the regularity
constant $\gamma^*$, the data $D=(A,\varphi)$ and the parameters $\theta_1$ used in {\bf Adaptive Algorithm $C$},
such that any two consecutive iterates $k$ and $k+1$ have the property
\begin{eqnarray*}\label{Adaptive_Convergence_Source}
\|(I-R_{k+1})K\bar{f}_k\|^2_{a,\Omega}+\gamma\eta^2(\bar{f}_k;\mathcal{T}_{k+1}) \leq \hat{\alpha}^2
 \Big( \|(I-R_k)K\bar{f}_k\|^2_{a,\Omega}+\gamma\eta^2(\bar{f}_k;\mathcal{T}_k)\Big).
\end{eqnarray*}
%where $\mathbb{N}=\{1, 2, 3, \cdots\}$.
\end{theorem}
%----------------------------------------------------------------------------------------------------------
In our analysis, we also need following lemmas.
\begin{lemma}(\cite{CasconKreuzerNochettoSiebert})\label{quasi optimality of the total error}
There exists a constant $\hat{C}_D$ depending only on data $D$
and regularity constant $\gamma^*$ such that
\begin{eqnarray*}
&&\|(I-R_k)Kf\|_{a,\Omega}^2+{\rm osc}^2(f-LR_kKf;\mathcal{T}_k) \nonumber\\
&\leq& \hat{C}_D\inf\limits_{v_k \in V_k}\big(\|Kf-v_k\|_{a,\Omega}^2
+{\rm osc}^2(f-Lv_k;\mathcal{T}_k)\big),\ \ \ \ \forall f\in L^2(\Omega).
\end{eqnarray*}
\end{lemma}
In this paper, we assume that the marking parameter $\theta_1$ satisfies $\theta_1\in(0,\theta_*)$
with $\theta_*$ defined in Assumption 5.8 of \cite{CasconKreuzerNochettoSiebert}.
%--------------------------------------------------------------------------------------------------------
\begin{lemma}(\cite{CasconKreuzerNochettoSiebert})\label{Error_estimate_Lower_Bound_Lemma}
For any function $f\in L^2(\Omega)$, let $\mathcal{T}_{k,*}$ and $V_{k,*}$ be a refinement of $\mathcal{T}_k$
and the finite element space over the mesh $\mathcal{T}_{k,*}$, such that the Ritz-Galerkin approximation of $Kf$
satisfies the energy decreasing property
\begin{eqnarray*}
&&\|(I-R_{{k,*}})Kf\|^2_{a,\Omega}+{\rm osc}^2(f-LR_{k,*}Kf;\mathcal{T}_{k,*}) \leq \tilde{\xi}^2_0\big(\|(I-R_{k})Kf\|^2_{a,\Omega}
+ {\rm osc}^2(f-LR_{k}Kf;\mathcal{T}_{k})\big)
\end{eqnarray*}
with $\tilde{\xi}_0^2\in(0,\frac{1}{2})$.
Then the set $\mathcal{T}_{k}\backslash(\mathcal{T}_{{k,*}}\cap\mathcal{T}_{k})$ of refined elements
satisfies the D\"{o}rfler property
\begin{eqnarray*}
\eta^2_{k}\big(f,R_{k}Kf;\mathcal{T}_{k}\backslash(\mathcal{T}_{{k,*}}\cap\mathcal{T}_{k})\big)
\geq \tilde{\theta}^2\eta^2_{k}(f,R_{k}Kf;\mathcal{T}_k),
\end{eqnarray*}
where $\tilde\theta=\theta_*\sqrt{1-2\tilde{\xi}_0^2}$.
\end{lemma}
%------------------------------------------------------------------------------------------------------------
\begin{remark}\label{remark_fk_vanish_in_osc}
In {\bf Adaptive Algorithm $C$}, since the right hand side term $\bar{f}_k$ of boundary value
problems are always piecewise polynomials, by the definition of data oscillation, we have
\begin{eqnarray}\label{eqn_vanish_f_osc}
{\rm osc}(\bar{f}_k-L\bar{u}_k;\mathcal{T}_k) = {\rm osc}(L\bar{u}_k;\mathcal{T}_k).
\end{eqnarray}
And similar properties of (\ref{eqn_vanish_f_osc}) will be used in the following parts of this paper.
\end{remark}
%------------------------------------------------------------------------------------------------------------
For simplicity, we only give the analysis for the first eigenpair approximation by the adaptive
multilevel correction algorithm ({\bf Adaptive Algorithm $C$}) in this paper.
\revise{Then the symbols $u_k$, $\bar u_k$, $\bar u_k^{(j)}$ and $u_k^{(j)}$ denote approximations for the first exact eigenfunction $\hat u^1$,
$\lambda_k$ and $\lambda_k^{(j)}$ are approximations to the first exact eigenvalue $\hat\lambda^1$}.
 Let us define the spectral projection  $E:V\rightarrow M(\hat{\lambda}^1)$ as follows
\begin{eqnarray}\label{Spectral_Projection}
a(v-Ev,w) = 0,\ \ \ \ \forall w\in M(\hat{\lambda}^1),\ \forall v\in V.
\end{eqnarray}
From the \revise{definitions} (\ref{definition_quantity_eta}) and (\ref{definition_quantity_delta}),
it is easy to show that $\delta_{\tilde{V}_{H}}(\hat{\lambda}^1) \leq \delta_{V_H}(\hat{\lambda}^1)$
and $\eta_a(\tilde{V}_{H})\leq \eta_a(V_H)$, where $V_H \subset \tilde{V}_{H}$.
Hence, the following properties are direct results of Lemma \ref{lemma_error_estimate_eigenproblem}.
\begin{corollary}\label{Error_Estimate_Crude_Lemma}
For each obtained eigenpair approximation
$(\tilde{\lambda}_{n_{\ell}}^{(j)},\tilde{u}_{n_{\ell}}^{(j)})\ (\ell >1,j\geq0)$ and $(\tilde{\lambda}_{n_{\ell}},\tilde{u}_{n_{\ell}})(\ell>1)$ in {\bf Adaptive Algorithm $C$}, the following estimates hold
\begin{eqnarray*}
\|\tilde{u}_{n_{\ell}}^{(j)}-E\tilde{u}_{n_{\ell}}^{(j)}\|_{a,\Omega}&\leq& C_{ea}\min\big\{\delta_{V_H}(\hat{\lambda}^1),
\|\hat u^1-u_{n_{\ell}}^{(j)}\|_{a,\Omega}\big\},\label{Error_u_u_k}\\
\|\tilde{u}_{n_{\ell}}^{(j)}-E\tilde{u}_{n_{\ell}}^{(j)}\|_{0,\Omega}&\leq& C_{e0}\eta_a(V_H)\|\tilde{u}_{n_{\ell}}^{(j)}-E\tilde{u}_{n_{\ell}}^{(j)}\|_{a,\Omega},\label{Error_u_u_k_-1}\\
\|\tilde{u}_{n_{\ell}}-E\tilde{u}_{n_{\ell}}\|_{a,\Omega}&\leq& C_{ea}\min\big\{\delta_{V_H}(\hat{\lambda}^1),
\|\hat u^1-u_{n_{\ell}}\|_{a,\Omega}\big\},\label{Error_u_u_k}\\
\|\tilde{u}_{n_{\ell}}-E\tilde{u}_{n_{\ell}}\|_{0,\Omega}&\leq& C_{e0}\eta_a(V_H)\|\tilde{u}_{n_{\ell}}-E\tilde{u}_{n_{\ell}}\|_{a,\Omega},\label{Error_u_u_k_-1}
\end{eqnarray*}
where $C_{ea}=\bar{C}_{ea}$ and $C_{e0} = 2\bar{C}_{e0}$.
\end{corollary}

%----------------------------------------------------------------------------------------------------------------
\subsection{Error estimate for eigenfunction approximation} % by the adaptive multilevel correction algorithm}

For $k\in \mathbb{N}$, $j\geq0$, we define $f_k^{(j)}$ and $w_k^{(j)}$ by
\begin{eqnarray}\label{definition_w_kj_f_kj}
f_{k}^{(j)}= \frac{ \bar{f}_k^{(j)}}{\|K\bar{f}_{k}^{(j)}\|_{a,\Omega}}\ \ \ \  \text{and}\ \ \ \ w_k^{(j)}=Kf_k^{(j)},
\end{eqnarray}
and $f_k$ and $w_k$ by
\begin{eqnarray}\label{definition_w_k_f_k}
f_{k}= \frac{ \bar{f}_k}{\|K\bar{f}_{k}\|_{a,\Omega}}\ \ \text{and}\ \ w_k=Kf_k.
\end{eqnarray}
It is obvious that $\|w_k^{(j)}\|_{a,\Omega}=\|w_k\|_{a,\Omega}=1$. From step 3 of {\bf Adaptive Algorithm $C$},
the following properties hold
\begin{eqnarray}\label{eqn_relation _tildeu_k_2_w_k}
u_k^{(j)}=\frac{R_kK\bar{f}_k^{(j)}}{\|R_kK\bar{f}_k^{(j)}\|_{a,\Omega}}=\frac{R_{k} w_{k}^{(j)}}{\|R_k w_k^{(j)}\|_{a,\Omega}}
\ \ \ \text{and}\  \ \ u_k=\frac{R_kK\bar{f}_k}{\|R_kK\bar{f}_k\|_{a,\Omega}}=\frac{R_{k} w_{k}}{\|R_k w_k\|_{a,\Omega}}.
\end{eqnarray}

%Let ${\color{red}\eta_a(V_H)}:=\eta_a(V_H)+\delta_{V_H}(\hat\lambda^1)$.
%Obviously from Lemma \ref{Error_Estimate_Crude_Lemma}, we know ${\color{red}\eta_a(V_H)}\ll 1$ when $H$ is small enough.
%------------------------------------------------------------------------------------------------------------------
\begin{lemma}\label{lemma_general_v-Ev_split}
For any functions $v\in V$ and $w\in V$ satisfying $\|v\|_{a,\Omega}=\|w\|_{a,\Omega}=1$, we have the following
inequality
\begin{eqnarray}\label{Inequality_4}
\|v-Ev\|_{a,\Omega} \leq \inf_{\tau \in \mathbb{R}}\|w-\tau v\|_{a,\Omega} + \|w-Ew\|_{a,\Omega}.
\end{eqnarray}
\end{lemma}
%-------------------------------------------------------------------------------------------------------------------
\begin{proof}
By the definition of $E$ in (\ref{Spectral_Projection}), we have $\|v-Ev\|_{a,\Omega} \leq \|v-\tau_0 Ew\|_{a,\Omega}$
for any $\tau_0 \in \mathbb{R}$. The triangle inequality leads to the following estimate
\begin{eqnarray}\label{eqn_general_v-Ev_split}
\|v-\tau_0 Ew\|_{a,\Omega} \leq \|v-\tau_0 w\|_{a,\Omega} + |\tau_0|\|w-Ew\|_{a,\Omega}.
\end{eqnarray}
We choose $\tau_0$ such that $\|v-\tau_0 w\|_{a,\Omega} = \inf_{\tau \in \mathbb{R}}\|v-\tau w\|_{a,\Omega}$.
It is obvious that $|\tau_0|\leq 1$.
Then from (\ref{eqn_general_v-Ev_split}), we obtain
\begin{eqnarray*}
\|v-Ev\|_{a,\Omega} &\leq& \inf_{\tau \in \mathbb{R}}\|v-\tau w\|_{a,\Omega} + \|w-Ew\|_{a,\Omega}
=\inf_{\tau \in \mathbb{R}}\|w-\tau v\|_{a,\Omega}+\|w-Ew\|_{a,\Omega}.
\end{eqnarray*}
This is the desired result (\ref{Inequality_4}) and the proof is complete.
\end{proof}
%-------------------------------------------------------------------------------------------------------------------
From (\ref{eqn_relation _tildeu_k_2_w_k}),
we know that $\inf_{\tau \in \mathbb{R}}\|w_k^{(j)}-\tau u_k^{(j)}\|_{a,\Omega}=\|w_k^{(j)}-R_kw_k^{(j)}\|_{a,\Omega}$.
Then the following theorem is a direct result of Lemma \ref{lemma_general_v-Ev_split}.
%------------------------------------------------------------------------------------------------------------------
\begin{theorem}\label{theorem_split_u_k-Eu_k}
Let $u_k^{(j)}\in V$ and $u_k\in V$ be produced by {\bf Adaptive Algorithm $C$}, $w_k^{(j)}$ and $w_k$ be defined by (\ref{definition_w_kj_f_kj}) and (\ref{definition_w_k_f_k}), respectively.
Then we have
\begin{eqnarray}
\|u_k^{(j)}-Eu_k^{(j)}\|_{a,\Omega} &\leq& \|(I-R_k)w_k^{(j)}\|_{a,\Omega} + \|w_k^{(j)}-Ew_k^{(j)}\|_{a,\Omega},\label{eqn_split_u_kj-Eu_kj}\\
\|u_k-Eu_k\|_{a,\Omega} &\leq& \|(I-R_k)w_k\|_{a,\Omega} + \|w_k-Ew_k\|_{a,\Omega}.\label{eqn_split_u_k-Eu_k}
\end{eqnarray}
\end{theorem}
%--------------------------------------------------------------------------------------------------------------
Theorem \ref{theorem_split_u_k-Eu_k} establishes a basic relation  between the error estimates $\|u_k-Eu_k\|_{a,\Omega}$ of the
finite element approximation produced by {\bf Adaptive Algorithm $C$} and $\|(I-R_k)w_k\|_{a,\Omega}$ of the associated finite
element projection.
%--------------------------------------------------------------------------------------------------------------
\begin{lemma}\label{lemma_equality_of_2kind_error}
For any given function $v\in V$ satisfying $\|v\|_{a,\Omega}=1$ and $a(v,\hat{u}^1)\geq0$,
the following inequality holds
\begin{eqnarray}\label{Inequality_5}
\|v-\hat{u}^1\|_{a,\Omega} \leq 2\|v-Ev\|_{a,\Omega}.
\end{eqnarray}
\end{lemma}
%------------------------------------------------------------------------------------------------------------------
\begin{proof}
The triangle inequality implies the following estimate
\begin{eqnarray}\label{Inequality_6}
\|v-\hat{u}^1\|_{a,\Omega} \leq \|v- Ev\|_{a,\Omega} + \|Ev - \hat{u}^1\|_{a,\Omega}.
\end{eqnarray}
Then we only need to prove that
\begin{eqnarray*}
\|Ev- \hat{u}^1\|_{a,\Omega} \leq \|v- Ev\|_{a,\Omega}.
\end{eqnarray*}
Since $\|Ev\|_{a,\Omega} \leq \|v\|_{a,\Omega}=1$ and $a(v,\hat{u}^1)\geq0$, the following inequality holds
\begin{eqnarray}\label{Inequality_7}
\|Ev-\hat{u}^1\|_{a,\Omega}^2 = (1-\|Ev\|_{a,\Omega})^2 \leq 1-\|Ev\|_{a,\Omega}^2 = \|v- Ev\|_{a,\Omega}^2.
\end{eqnarray}
The desired result (\ref{Inequality_5}) can be deduced by combining (\ref{Inequality_6}) and (\ref{Inequality_7}).
%and the proof is complete.
\end{proof}
%--------------------------------------------------------------------------------------------------------------
The following lemma gives the estimate for $\|w_k^{(j)}-Ew_k^{(j)}\|_{a,\Omega}$, which is the second term in the right hand
side of (\ref{eqn_split_u_kj-Eu_kj}).
%--------------------------------------------------------------------------------------------------------------
\begin{lemma}\label{lemma_estimate_of_w-Ew}
The following three propositions hold
\begin{enumerate}
\item If $k = n_{\ell}$ for any $1\neq\ell \in \mathbb{N}$, which implies $\bar{f}_k^{(j+1)}=\tilde{u}_{k}^{(j)}$($ j\geq 0$)
is produced by the eigenvalue solving module $\textsf{ESOLVE}$ in step 5 of {\bf Adaptive Algorithm $C$},
we have
\begin{eqnarray}\label{Inequality_1}
\|w_{n_{\ell}}^{(j+1)}-Ew_{n_{\ell}}^{(j+1)}\|_{a,\Omega} \leq \frac{C_{wu}}{2\hat{\lambda}^1\|K\tilde{u}_{n_{\ell}}^{(j)}\|_{a,\Omega}}
\eta_a(V_H)\|u_{n_{\ell}}^{(j)}-Eu_{n_{\ell}}^{(j)}\|_{a,\Omega},
\end{eqnarray}
\revise{where} %  $\hat{\lambda}^2$ denotes the second exact eigenvalue and
\begin{eqnarray*}
C_{wu}:=\frac{ 4C_{e0}C_{ea}\hat{\lambda}^1}{c_a}.
\end{eqnarray*}
\item If $k=n_{\ell}+1$ for any $1\neq \ell \in \mathbb{N}$, which implies $\bar{f}_k^{(0)}=\tilde{u}_{k-1}$
is produced by the eigenvalue solving module $\textsf{ESOLVE}$ in step 6 of {\bf Adaptive Algorithm $C$},
we have
\begin{eqnarray}\label{Inequality_1_1}
\|w_{n_{\ell}+1}^{(0)}-Ew_{n_{\ell}+1}^{(0)}\|_{a,\Omega} \leq \frac{C_{wu}}{2\hat{\lambda}^1\|K\tilde{u}_{n_{\ell}}\|_{a,\Omega}}
\eta_a(V_H)\|u_{n_{\ell}}-Eu_{n_{\ell}}\|_{a,\Omega}.
\end{eqnarray}
\item If $n_{\ell}+1<k \leq n_{\ell+1}$ for any $\ell \in \mathbb{N}$, which implies $\bar{f}_k^{(0)}=u_{k-1}$ is
produced by the boundary value problem solving module $\textsf{LSOLVE}$ in step 3 of {\bf Adaptive Algorithm $C$},
we have
\begin{eqnarray}\label{eqn_estimate_of_w-Ew_LSOLVE}
\|w_k^{(0)}-Ew_k^{(0)}\|_{a,\Omega}\leq \frac{C_{ww}\eta_a(V_{k-1} )\|(I-R_{k-1})w_{k-1} \|_{a,\Omega}}{\revise{2\hat{\lambda}^1\|Kw_{k-1}\|_{a,\Omega}}}
+\revise{\|w_{k-1} -Ew_{k-1} \|_{a,\Omega}},
\end{eqnarray}
where
\begin{eqnarray*}
\revise{C_{ww}:=\frac{2C_{an}\hat{\lambda}^1}{c_a}}.
\end{eqnarray*}
\end{enumerate}
\end{lemma}
%------------------------------------------------------------------------------------------------------------
\begin{proof}
If $k = n_{\ell}$ for any $1\neq\ell \in \mathbb{N}$, from the algorithm definition, we have
\begin{eqnarray*}
w_k^{(j+1)}=Kf_k^{(j+1)}=\frac{K\tilde{u}_{k}^{(j)}}{\|K\tilde{u}_{k}^{(j)}\|_{a,\Omega}}\ \ \ \ {\rm for}\ j\geq 0.
\end{eqnarray*}
By the optimality of spectral projection $E$ in (\ref{Spectral_Projection}), the following inequality holds
\begin{eqnarray}\label{eqn_proof_estimate_w-Ew_1}
\|w_k^{(j+1)}-Ew_k^{(j+1)}\|_{a,\Omega} \leq \left\|w_k^{(j+1)}-\frac{E\tilde{u}_{k}^{(j)}}{\hat{\lambda}^1
\|K\tilde{u}_{k}^{(j)}\|_{a,\Omega}}\right\|_{a,\Omega}.
\end{eqnarray}
Since $E\tilde{u}_{k}^{(j)}\in M(\hat{\lambda}^1)$, we have $E\tilde{u}_{k}^{(j)} = \hat{\lambda}^1KE\tilde{u}_{k}^{(j)}$.
Combining (\ref{eqn_K_norm}), (\ref{definition_w_k_f_k}), (\ref{eqn_proof_estimate_w-Ew_1}),
Corollary \ref{Error_Estimate_Crude_Lemma} and Lemma \ref{lemma_equality_of_2kind_error}
leads to the following estimates
\begin{eqnarray}\label{eqn_estimate_of_w-Ew_1}
&&\|w_k^{(j+1)}-Ew_k^{(j+1)}\|_{a,\Omega} \leq \frac{\|K(\tilde{u}_{k}^{(j)}-E\tilde{u}_{k}^{(j)})\|_{a,\Omega}}{\|K\tilde{u}_{k}^{(j)}\|_{a,\Omega}}
\leq \frac{\|\tilde{u}_{k}^{(j)}-E\tilde{u}_{k}^{(j)}\|_{0,\Omega}}{c_a\|K\tilde{u}_{k}^{(j)}\|_{a,\Omega}}\nonumber\\
&\leq&\frac{C_{e0}}{c_a\|K\tilde{u}_{k}^{(j)}\|_{a,\Omega}}\eta_a(V_H)\|\tilde{u}_{k}^{(j)}-E\tilde{u}_{k}^{(j)}\|_{a,\Omega}
\leq \frac{C_{e0}C_{ea}}{c_a\|K\tilde{u}_{k}^{(j)}\|_{a,\Omega}}\eta_a(V_H)\|\hat u^1-u_{k}^{(j)}\|_{a,\Omega}\nonumber\\
&\leq&\frac{2C_{e0}C_{ea}}{c_a\|K\tilde{u}_{k}^{(j)}\|_{a,\Omega}}\eta_a(V_H)\|u_{k}^{(j)}-Eu_{k}^{(j)}\|_{a,\Omega},
\end{eqnarray}
\comm{which is the desired result (\ref{Inequality_1})}.
%Applying Corollary \ref{Error_Estimate_Crude_Lemma} and Lemma \ref{lemma_equality_of_2kind_error} leads to
%the desired result (\ref{Inequality_1}).
%\begin{eqnarray*}
%\|w_k-Ew_k\|_{a,\Omega} \leq
%\frac{C_{wu}}{\hat{\lambda}^2\|K\tilde{u}_{k-1}\|_{a,\Omega}}{\color{red}\eta_a(V_H)}\|u_{k-1}-Eu_{k-1}\|_{a,\Omega},
%\end{eqnarray*}
%which is the desired result

If $k=n_{\ell}+1$ for any $1\neq \ell \in \mathbb{N}$, we can prove (\ref{Inequality_1_1}) with a similar procedure as above.

If $n_{\ell}+1<k \leq n_{\ell+1}$ for any $\ell \in \mathbb{N}$, from the algorithm definition, we have
\begin{eqnarray*}
w_k^{(0)}=\frac{Ku_{k-1}}{\|Ku_{k-1}\|_{a,\Omega}}.
\end{eqnarray*}
Let us introduce an auxiliary function
\begin{eqnarray}\label{Definition_Tilde_w_k}
\tilde{w}_k = \frac{Kw_{k-1}}{\|Kw_{k-1}\|_{a,\Omega}}.
\end{eqnarray}
From Lemma \ref{lemma_general_v-Ev_split}, the following inequality holds
\begin{eqnarray}\label{eqn_w2w_1}
\|w_k^{(0)} -Ew_k^{(0)} \|_{a,\Omega} \leq \inf_{\tau \in \mathbb{R}}\|\tilde{w}_k -\tau w_k^{(0)}\|_{a,\Omega}
+ \|\tilde{w}_k -E\tilde{w}_k \|_{a,\Omega}.
\end{eqnarray}
Let us define $\tau_2$ as follows
$$\tau_2 = \frac{\|Ku_{k-1}\|_{a,\Omega}\|R_{k-1}w_{k-1}\|_{a,\Omega}}{\|Kw_{k-1}\|_{a,\Omega}}.$$
Since $$u_{k-1} = \frac{R_{k-1}w_{k-1}}{\|R_{k-1}w_{k-1}\|_{a,\Omega}},$$
the following estimate holds
\begin{eqnarray}\label{w2w_6}
\inf_{\tau \in \mathbb{R}}\|\tilde{w}_k -\tau w_k^{(0)}\|_{a,\Omega} \leq \|\tilde{w}_k -\tau_2 w_k^{(0)}\|_{a,\Omega}
=  \frac{\|K(I-R_{k-1})w_{k-1}\|_{a,\Omega}}{\|Kw_{k-1}\|_{a,\Omega}}.
\end{eqnarray}
From (\ref{eqn_K_norm}) and Lemma \ref{Aubin-Nitsche}, we can derive
\begin{eqnarray}\label{w2w_7}
\|K(I-R_{k-1})w_{k-1}\|_{a,\Omega}
\leq\frac{\|(I-R_{k-1})w_{k-1}\|_{0,\Omega}}{c_a}\leq \frac{C_{an}}{c_a}\eta_a(V_{k-1} )\|(I-R_{k-1})w_{k-1} \|_{a,\Omega}.
\end{eqnarray}
Then combining (\ref{w2w_6}) and (\ref{w2w_7}) leads to
\begin{eqnarray}\label{w2w_5}
\inf_{\tau \in \mathbb{R}}\|\tilde{w}_k -\tau w_k^{(0)}\|_{a,\Omega}
\leq \frac{C_{an}}{c_a\|Kw_{k-1}\|_{a,\Omega}}\eta_a(V_{k-1} )\|(I-R_{k-1})w_{k-1} \|_{a,\Omega}.
\end{eqnarray}
\revise{We assume $w_{k-1}\in V$ has the expansion $w_{k-1} = \sum_{i=1}^{\infty} \beta_i \hat{u}^i$. Then}
\begin{eqnarray}
Kw_{k-1} = \sum_{i=1}^{\infty} \frac{\beta_i}{\hat{\lambda}^i} \hat{u}^i.
\end{eqnarray}
By the definition of spectral projection $E$ in (\ref{Spectral_Projection}), (\ref{Definition_Tilde_w_k}) and the fact $\|w_{k-1}\|_{a,\Omega}^2=\sum_{i=1}^{\infty} \beta_i^2=1$, we have
\begin{eqnarray}\label{w2w_44}
&&\|\tilde{w}_k-E\tilde{w}_k\|_{a,\Omega}^2= \|\tilde{w}_k\|_{a,\Omega}^2-\|E\tilde{w}_k\|_{a,\Omega}^2
= 1-\frac{\big(\frac{\beta_1}{\hat{\lambda}^1}\big)^2}{\sum_{i=1}^{\infty} \big(\frac{\beta_i}{\hat{\lambda}^i}\big)^2}\nonumber\\
&\leq& 1-\frac{\big(\frac{\beta_1}{\hat{\lambda}^1}\big)^2}{\big(\frac{1}{\hat{\lambda}^1}\big)^2\sum_{i=1}^{\infty} \beta_i^2}
= 1 - \beta_1^2 = \|w_{k-1}-Ew_{k-1}\|_{a,\Omega}^2.
\end{eqnarray}
%By the definition of spectral projection $E$ in (\ref{Spectral_Projection}) and (\ref{Definition_Tilde_w_k}), we have
%\begin{eqnarray}\label{w2w_1}
%&&\|\tilde{w}_k -E\tilde{w}_k \|_{a,\Omega}^2 = a(\tilde{w}_k -E\tilde{w}_k ,\tilde{w}_k -E\tilde{w}_k )
%=a\Big(\tilde{w}_k -\frac{Ew_{k-1}}{\hat{\lambda}^1\|Kw_{k-1}\|_{a,\Omega}} ,\tilde{w}_k -E\tilde{w}_k \Big)\nonumber\\
%&=&\frac{(w_{k-1} -Ew_{k-1} ,\tilde{w}_k -E\tilde{w}_k )}{\|Kw_{k-1}\|_{a,\Omega}}
%\leq\frac{\|w_{k-1} -Ew_{k-1} \|_{0,\Omega}\|\tilde{w}_k -E\tilde{w}_k \|_{0,\Omega}}{\|Kw_{k-1}\|_{a,\Omega}}.
%\end{eqnarray}
%It is obvious that the properties
%$w_{k-1} -Ew_{k-1}\in M(\hat{\lambda}^1)^{\bot}$ and $\tilde{w}_k -E\tilde{w}_k \in M(\hat{\lambda}^1)^{\bot}$
%imply following estimates
%\begin{eqnarray}\label{w2w_2}
%\|w_{k-1} -Ew_{k-1} \|_{0,\Omega} \leq \frac{1}{\sqrt{\hat{\lambda}^2}}\|w_{k-1} -Ew_{k-1} \|_{a,\Omega}
%\end{eqnarray}
%and
%\begin{eqnarray}\label{w2w_3}
%\|\tilde{w}_k -E\tilde{w}_k \|_{0,\Omega} \leq \frac{1}{\sqrt{\hat{\lambda}^2}}
%\|\tilde{w}_k -E\tilde{w}_k \|_{a,\Omega}.
%\end{eqnarray}
%Combining (\ref{w2w_1}), (\ref{w2w_2}) and (\ref{w2w_3}) leads to
%\begin{eqnarray}\label{w2w_4}
%\|\tilde{w}_k -E\tilde{w}_k \|_{a,\Omega}
%\leq \frac{\|w_{k-1} -Ew_{k-1} \|_{a,\Omega} }{\hat{\lambda}^2\|Kw_{k-1}\|_{a,\Omega}}.
%\end{eqnarray}
From (\ref{eqn_w2w_1}), (\ref{w2w_5}) and (\ref{w2w_44}), we can obtain
the desired result (\ref{eqn_estimate_of_w-Ew_LSOLVE}) and the proof is complete.
\end{proof}
%-----------------------------------------------------------------------------------------------------------
From Lemma \ref{lemma_estimate_of_w-Ew}, it is required to obtain the lower bounds of $\|K\tilde{u}_{n_{\ell}}^{(j)}\|_{a,\Omega}$,
$\|K\tilde{u}_{n_{\ell}}\|_{a,\Omega}$ and $\|Kw_k\|_{a,\Omega}$.
For this aim, we first state the following lemma.
\begin{lemma}\label{lemma_bound_Kv_anorm}
For any function $v\in V$ with $\|v\|_{a,\Omega}=1$, we have the following inequality
\begin{eqnarray}\label{Inequality_3}
\frac{1}{\hat{\lambda}^1}-\frac{2\|v-Ev\|_{a,\Omega}}{c_a^2}\leq\|Kv\|_{a,\Omega}.
\end{eqnarray}
\end{lemma}
%------------------------------------------------------------------------------------------------------------
\begin{proof}
First, from (\ref{eqn_K_norm}) and Lemma \ref{lemma_equality_of_2kind_error} we have
\begin{eqnarray}\label{eqn_proof_bound_Kv_anorm_1}
\|K\hat{u}^1\|_{a,\Omega}-\|Kv\|_{a,\Omega}\leq \|Kv-K\hat{u}^1\|_{a,\Omega}
\leq\frac{\|v-\hat{u}^1\|_{a,\Omega}}{c_a^2}\leq \frac{2\|v-Ev\|_{a,\Omega}}{c_a^2}.
\end{eqnarray}
The property (\ref{Operator_Eigenvalue_Problem}) for $\hat{u}^1$ leads to the following equalities
\begin{eqnarray}\label{Inequality_2}
\|K\hat{u}^1\|_{a,\Omega} = \frac{\|\hat{u}^1\|_{a,\Omega}}{\hat{\lambda}^1}=\frac{1}{\hat{\lambda}^1}.
\end{eqnarray}
Then the desired result (\ref{Inequality_3}) can be deduced from (\ref{eqn_proof_bound_Kv_anorm_1})
and (\ref{Inequality_2}).
\end{proof}

%\begin{remark}
%\comm{
%From (\ref{Inequality_3}), we note that $\|K\tilde{u}_{n_{\ell}}^{(j)}\|_{a,\Omega}$,
%$\|K\tilde{u}_{n_{\ell}}\|_{a,\Omega}$ and $\|Kw_k\|_{a,\Omega}$ have an approximate lower bound $1/\hat{\lambda}^1$.
%And as a result, the ratios $1/(\hat{\lambda}^2\|K\tilde{u}_{n_{\ell}}^{(j)}\|_{a,\Omega})$,
%$1/(\hat{\lambda}^2\|K\tilde{u}_{n_{\ell}}\|_{a,\Omega})$ and $1/(\hat{\lambda}^2\|Kw_k\|_{a,\Omega})$ in
%(\ref{Inequality_1}), (\ref{Inequality_1_1}) and (\ref{eqn_estimate_of_w-Ew_LSOLVE}) of Lemma \ref{lemma_estimate_of_w-Ew},
%respectively, have an approximate upper bound $\hat{\lambda}^1/\hat{\lambda}^2$, which is the converging factor of
%inverse power method and is dependent on the gap between eigenvalues $\hat{\lambda}^1$ and $\hat{\lambda}^2$.
%However, in this paper we give our analysis by estimate these ratios with a crude upper bound $1$.}
%\end{remark}
%-------------------------------------------------------------------------------------------------------------
\begin{corollary}\label{corollary_lowerbound_Kunl}
When $H$ is small enough, for $\ell\in\mathbb{N}$ and $j\geq 0$, we have
\begin{eqnarray}\label{eqn_lowerbound_Kunl}
\revise{\frac{1}{2\hat{\lambda}^1}}\leq \|K\tilde{u}_{n_{\ell}}\|_{a,\Omega}\ \ \text{and}\ \  \revise{\frac{1}{2\hat{\lambda}^1}}\leq\|K\tilde{u}_{n_{\ell}}^{(j)}\|_{a,\Omega}.
\end{eqnarray}
\end{corollary}
\begin{proof}
We only prove the first inequality and the second one can be proved by a similar procedure. For any $\ell \in \mathbb{N}$,
 Lemma \ref{lemma_bound_Kv_anorm} implies
\begin{eqnarray}\label{eqn_proof_lowerbound_Kunl_1}
\frac{1}{\hat{\lambda}^1}-\frac{2\|\tilde{u}_{n_{\ell}}-E\tilde{u}_{n_{\ell}}\|_{a,\Omega}}{c_a^2}
\leq\|K\tilde{u}_{n_{\ell}}\|_{a,\Omega}.
\end{eqnarray}
From Corollary \ref{Error_Estimate_Crude_Lemma}, when $H$ is small enough, the following inequalities hold
\begin{eqnarray}\label{eqn_the_first_two_step_3}
\|\tilde{u}_{n_{\ell}}-E\tilde{u}_{n_{\ell}}\|_{a,\Omega}
\leq C_{ea} \delta_{V_H}(\hat{\lambda}^1)\leq \revise{\frac{c_a^2}{4\hat{\lambda}^1}}.
\end{eqnarray}
Then from (\ref{eqn_proof_lowerbound_Kunl_1}) and (\ref{eqn_the_first_two_step_3}), we can obtain the
desired result (\ref{eqn_lowerbound_Kunl}) and the proof is complete.
\end{proof}
%-----------------------------------------------------------------------------------------------------
%----------------------------------------------------------------------------------------------------------
\begin{lemma}\label{lemma_break_in_finite_steps}
When $H$ is small enough, if $\bar{u}_{n_{\ell}}^{(i)}$ $(0\leq i\leq j)$ produced by {\bf Adaptive Algorithm $C$} satisfies
\begin{eqnarray}\label{eqn_decision_condition_before_break}
\eta_{n_{\ell}}(\bar{f}_{n_{\ell}}^{(i)},\bar{u}_{n_{\ell}}^{(i)};\mathcal{T}_{n_{\ell}}) \leq \theta_2^{i+1} \eta_{\ell-1},\ \ 0\leq i \leq j,
\end{eqnarray}
we have the following property
\begin{eqnarray}\label{eqn_break_in_finite_steps}
\|u_{n_{\ell}}^{(j)}-Eu_{n_{\ell}}^{(j)}\|_{a,\Omega}\leq \revise{4\theta_2^{j+1}\hat{\lambda}^1}C_{\rm up}\eta_{\ell-1}+\big(C_{wu}\eta_a(V_H)\big)^{j}\|u_{n_{\ell}}^{(0)}-Eu_{n_{\ell}}^{(0)}\|_{a,\Omega}.
\end{eqnarray}
\end{lemma}
\begin{proof}
We prove (\ref{eqn_break_in_finite_steps}) by the induction method.
It is obvious that (\ref{eqn_break_in_finite_steps}) holds for $j=0$. Now we suppose (\ref{eqn_break_in_finite_steps})
holds for the case $j-1$ and consider the case $j$.
From Theorem \ref{theorem_split_u_k-Eu_k}, Lemmas \ref{lemma_upper_lower_bound_boundaryvalue}, \ref{lemma_estimate_of_w-Ew} and Corollary \ref{corollary_lowerbound_Kunl}, we have
\begin{eqnarray}\label{eqn_proof_break_in_finite_steps_1}
&&\|u_{n_{\ell}}^{(j)}-Eu_{n_{\ell}}^{(j)}\|_{a,\Omega} \leq \|(I-R_{n_{\ell}})w_{n_{\ell}}^{(j)}\|_{a,\Omega}+\|w_{n_{\ell}}^{(j)}-Ew_{n_{\ell}}^{(j)}\|_{a,\Omega}\nonumber\\
&\leq& C_{\rm up} \eta(f_{n_{\ell}}^{(j)};\mathcal{T}_{n_{\ell}})
+ C_{wu}\eta_a(V_H) \|u_{n_{\ell}}^{(j-1)}-Eu_{n_{\ell}}^{(j-1)}\|_{a,\Omega}.
\end{eqnarray}
Note that $\bar{f}_{n_{\ell}}^{(j)} = \tilde{u}_{n_{\ell}}^{(j-1)}$ and $\|\tilde{u}_{n_{\ell}}^{(j-1)}\|_{a,\Omega}=1$.
Then from (\ref{eqn_lowerbound_Kunl}) and (\ref{eqn_decision_condition_before_break}), we can derive
\begin{eqnarray}\label{eqn_update_decision_condition}
&&\eta(f_{n_{\ell}}^{(j)};\mathcal{T}_{n_{\ell}}) = \frac{1}{\|K\bar{f}_{n_{\ell}}^{(j)}\|_{a,\Omega}}
\eta(\bar{f}_{n_{\ell}}^{(j)};\mathcal{T}_{n_{\ell}})
\leq \revise{2\hat{\lambda}^1}\theta_2^{j+1}\eta_{\ell-1}.
\end{eqnarray}
Since $H$ is small enough, the inequality $2C_{wu}\eta_a(V_H)\leq \theta_2$ holds.
Combining (\ref{eqn_proof_break_in_finite_steps_1}) and (\ref{eqn_update_decision_condition}) leads to
\begin{eqnarray*}
&&\|u_{n_{\ell}}^{(j)}-Eu_{n_{\ell}}^{(j)}\|_{a,\Omega}\nonumber\\
 &\leq& \revise{2C_{\rm up}\hat{\lambda}^1}\theta_2^{j+1}\eta_{\ell-1}+C_{wu}\eta_a(V_H)\Big(\revise{4\theta_2^{j}\hat{\lambda}^1}C_{\rm up}\eta_{\ell-1}
 +\big(C_{wu}\eta_a(V_H)\big)^{j-1}\|u_{n_{\ell}}^{(0)}-Eu_{n_{\ell}}^{(0)}\|_{a,\Omega}\Big)\nonumber\\
&=&\revise{\big(2\theta_2+4C_{wu}\eta_a(V_H)\big)\theta_2^{j}\hat{\lambda}^1}C_{\rm up}\eta_{\ell-1}+\big(C_{wu}\eta_a(V_H)\big)^{j}\|u_{n_{\ell}}^{(0)}-Eu_{n_{\ell}}^{(0)}\|_{a,\Omega}\nonumber\\
&\leq&\revise{4\theta_2^{j+1}\hat{\lambda}^1}C_{\rm up}\eta_{\ell-1} +\big(C_{wu}\eta_a(V_H)\big)^{j}\|u_{n_{\ell}}^{(0)}-Eu_{n_{\ell}}^{(0)}\|_{a,\Omega},
\end{eqnarray*}
which means (\ref{eqn_break_in_finite_steps}) holds for the case $j$ and the proof is complete.
\end{proof}
%--------------------------------------------------------------------------------------------
In this paper, we assume $\hat{u}^1$ is not a piecewise polynomial.
Then from Lemma \ref{lemma_break_in_finite_steps}, we can conclude that there exists $j_{n_{\ell}}\in \mathbb{N}$ such that
\begin{eqnarray}\label{eqn_decision_condition_after_break}
\eta(\bar{f}_{n_{\ell}}^{(j_{n_{\ell}})};\mathcal{T}_{n_{\ell}})> \theta_2^{j_{n_{\ell}}+1}\eta_{\ell-1}.
\end{eqnarray}
Otherwise, from (\ref{eqn_break_in_finite_steps}), we have
\begin{eqnarray*}
\lim_{j \rightarrow +\infty}\|u_{n_{\ell}}^{(j)}-Eu_{n_{\ell}}^{(j)}\|_{a,\Omega} = 0,
\end{eqnarray*}
which implies $\hat{u}^1\in V_{n_{\ell}}$.

As for the $k$-th level mesh satisfying $n_{\ell}<k<n_{\ell+1}$ for any $\ell \in \mathbb{N}$,
it is easy to know that  $\bar{f}_k=\bar{f}_k^{(0)}$ and
\begin{eqnarray*}
\eta(\bar{f}_k^{(0)};\mathcal{T}_k)> \theta_2 \eta_{\ell}.
\end{eqnarray*}
Thus we define $j_k = 0$ in this case. Then from {\bf Adaptive Algorithm $C$}, we know
\begin{eqnarray*}
\bar{f}_k=\bar{f}_k^{(j_k)},\ \  \bar{u}_k=\bar{u}_k^{(j_k)},\ \ f_k=f_k^{(j_k)},\ \ w_k=w_k^{(j_k)} \ \ \text{and}\ \  u_k=u_k^{(j_k)}.
\end{eqnarray*}
Moreover, on mesh $\mathcal{T}_{n_{\ell}}$ we have
\begin{eqnarray*}
\eta(\bar{f}_{n_{\ell}}^{(j_{n_{\ell}}-1)};\mathcal{T}_{n_{\ell}})\leq \theta_2^{j_{n_{\ell}}}\eta_{\ell-1}\leq \theta_2^{j_{n_{\ell}}-1}\eta(\bar{f}_{n_{\ell}-1};\mathcal{T}_{n_{\ell}-1}),
\end{eqnarray*}
which suggests that the number of iterations depends on the ratio of discretization errors on \revise{meshes}
$\mathcal{T}_{n_{\ell}}$ and $\mathcal{T}_{n_{\ell}-1}$. In general case, we assume these ratios are bounded by
a constant and the numbers of iterations are bounded too.
%------------------------------------------------------------------------------------------------------
\section{Convergence analysis}
In this section, we deduce the convergence of {\bf Adaptive Algorithm $C$} with the help of results introduced in Section 2.
%-------------------------------------------------------------------------------------------------------
\subsection{Properties of the first two levels }\label{section_property_first_two_levels}
This subsection is only concerned with the first two levels of {\bf Adaptive Algorithm $C$}
and shows some properties. From the definition of {\bf Adaptive Algorithm $C$}, we have
\begin{eqnarray*}\label{eqn_the_first_level_eigen_prob}
a(u_1,v_1) = \lambda_1 (u_1,v_1),\ \ \ \forall v_1 \in V_1
\end{eqnarray*}
and
\begin{eqnarray*}
a(\bar{u}_2^{(0)},v_2) = (u_1,v_2),\ \ \ \forall v_2 \in V_2.
\end{eqnarray*}
Since $\bar{f}_1=u_1=\tilde{u}_1=\bar{f}_2^{(0)}$ and the definition (\ref{definition_w_k_f_k}),
the following properties hold
\begin{eqnarray*}\label{Equality_1}
f_1=f_2^{(0)}=\frac{u_1}{\|Ku_1\|_{a,\Omega}}\ \ \text{and}\ \  w_1=w_2^{(0)}=\frac{Ku_1}{\|Ku_1\|_{a,\Omega}}.
\end{eqnarray*}
%-------------------------------------------------------------------------------------------------------
%-------------------------------------------------------------------------------------------------------

First, we state error estimates for $\|w_1-Ew_1\|_{a,\Omega}$ and $\|w_2^{(0)}-Ew_2^{(0)}\|_{a,\Omega}$ in the following
lemma.
\begin{lemma}\label{lemma_estimate_of_w-Ew_by_etal_for_k=2}
When the mesh size $H$ is small enough, we have following estimates for the first two levels
\begin{eqnarray}
\|w_1-Ew_1\|_{a,\Omega} \leq C_{w,1}\eta_a(V_H)\eta(f_1;\mathcal{T}_1),\label{eqn_the_first_two_step_10}\\
\|w_2^{(0)}-Ew_2^{(0)}\|_{a,\Omega}\leq C_{w,1}\eta_a(V_H)\eta(f_1;\mathcal{T}_1),\label{eqn_the_first_two_step_8}
\end{eqnarray}
where $$C_{w,1} = \frac{4\hat{\lambda}^1C_{e0}C_{\rm up}}{c_a}.$$
\end{lemma}
%--------------------------------------------------------------------------------------------------------
\begin{proof}
From (\ref{eqn_lowerbound_Kunl}) for $\ell=1$ and analogous to (\ref{eqn_estimate_of_w-Ew_1}),
we can derive
\begin{eqnarray}\label{eqn_the_first_two_step_1}
&&\|w_1-Ew_1\|_{a,\Omega} \leq \frac{ C_{e0}}{c_a\|Ku_1\|_{a,\Omega}}
\eta_a(V_1)\|u_1-Eu_1\|_{a,\Omega}\nonumber\\
&\leq&\frac{ C_{e0}}{c_a\|Ku_1\|_{a,\Omega}}\eta_a(V_H)\|u_1-Eu_1\|_{a,\Omega}
\leq\frac{\revise{2\hat{\lambda}^1} C_{e0}}{c_a}\eta_a(V_H)\|u_1-Eu_1\|_{a,\Omega}.
\end{eqnarray}
From Theorem \ref{theorem_split_u_k-Eu_k}, the following inequality holds
\begin{eqnarray}\label{eqn_the_first_two_step_5}
\|u_1-Eu_1\|_{a,\Omega} \leq \|w_1-R_1 w_1\|_{a,\Omega} + \|w_1-Ew_1\|_{a,\Omega}.
\end{eqnarray}
\revise{Let us choose $H$ to be small enough such that
\begin{eqnarray}\label{Inequality_31}
2\hat\lambda^1C_{e0}\eta_a(V_H)\leq \frac{c_a}{2}.
\end{eqnarray}
Combining (\ref{eqn_the_first_two_step_1}), (\ref{eqn_the_first_two_step_5}) and (\ref{Inequality_31}) leads to}
\begin{eqnarray}\label{eqn_the_first_two_step_6}
\|w_1-Ew_1\|_{a,\Omega} \leq \frac{2\hat{\lambda}^1C_{e0}}{c_a-\revise{2\hat{\lambda}^1}C_{e0}\eta_a(V_H)}\eta_a(V_H)
\|w_1-R_1w_1\|_{a,\Omega}\leq\frac{4\hat{\lambda}^1C_{e0}}{c_a}\eta_a(V_H)\|w_1-R_1w_1\|_{a,\Omega}.
\end{eqnarray}
From Lemma \ref{lemma_upper_lower_bound_boundaryvalue} and (\ref{eqn_the_first_two_step_6}),
we obtain the desired result (\ref{eqn_the_first_two_step_10}).
Since $w_2^{(0)}=w_1$, $\|w_2^{(0)}-Ew_2^{(0)}\|_{a,\Omega}$ has the same estimate (\ref{eqn_the_first_two_step_8}).
\end{proof}
%--------------------------------------------------------------------------------------------------------
\begin{lemma}\label{lemma_estimate_of_w-Ew_by_etak_for_k=2}
When the mesh size $H$ is small enough, we have the following estimate for the second level
\begin{eqnarray}\label{eqn_estimate_of_w-Ew_by_etak_for_k=2}
\|w_2-Ew_2\|_{a,\Omega}\leq C_{w,2}\eta_a(V_H)\eta(f_2;\mathcal{T}_2),\label{eqn_the_first_two_step_9}
\end{eqnarray}
where $$C_{w,2} = \max\Big\{\frac{C_{w,1}}{\theta_2},\frac{8C_{wu}\hat{\lambda}^1C_{\rm up}}{\theta_2c_a^2}\Big\}.$$
\end{lemma}
%--------------------------------------------------------------------------------------------------------
\begin{proof}
We prove the assertion case by case.

\textbf{Case 1:} If $\eta(\bar{f}_2^{(0)};\mathcal{T}_2) > \theta_2 \eta(\bar{f}_1;\mathcal{T}_1)$,
then from step 5 of {\bf Adaptive Algorithm $C$} and definition (\ref{definition_w_k_f_k}),
we have $\bar{f}_2 = \bar{f}_2^{(0)} = \bar{f}_1$ and $w_2=w_2^{(0)}$.
Thus, the following inequality holds
\begin{eqnarray}\label{eqn_update_decision_condition_k=2}
\eta(f_2;\mathcal{T}_2)=\frac{\eta(\bar{f}_2^{(0)};\mathcal{T}_2)}{\|K\bar{f}_2^{(0)}\|_{a,\Omega}}>\frac{\theta_2 \eta(\bar{f}_1;\mathcal{T}_1)}{\|K\bar{f}_1\|_{a,\Omega}}= \theta_2 \eta(f_1;\mathcal{T}_1).
\end{eqnarray}
\revise{From Lemma \ref{lemma_estimate_of_w-Ew_by_etal_for_k=2} and (\ref{eqn_update_decision_condition_k=2})}, we have
\begin{eqnarray}\label{eqn_proof_estimate_of_w-Ew_by_etak_for_k=2_case1}
\|w_2-Ew_2\|_{a,\Omega}\leq C_{w,1}\eta_a(V_H)\eta(f_1;\mathcal{T}_1)< \frac{C_{w,1}}{\theta_2}\eta_a(V_H)\eta(f_2;\mathcal{T}_2).
\end{eqnarray}

\textbf{Case 2:} If $\eta(\bar{f}_2^{(0)};\mathcal{T}_2) \leq \theta_2 \eta(\bar{f}_1;\mathcal{T}_1)$, then we have $j_2>0$ and $n_2=2$.
From Lemma \ref{lemma_break_in_finite_steps}, inequalities (\ref{eqn_decision_condition_before_break})
and (\ref{eqn_break_in_finite_steps}) hold for $j=j_2-1$. Then from (\ref{Inequality_1}),
(\ref{eqn_lowerbound_Kunl}), (\ref{eqn_break_in_finite_steps}) and (\ref{eqn_decision_condition_after_break}), we can derive
\begin{eqnarray}\label{eqn_estimate_w2-Ew2_by_eta2_1}
&&\|w_2^{(j_2)}-Ew_2^{(j_2)}\|_{a,\Omega}\leq C_{wu}\eta_a(V_H)\|u_2^{(j_2-1)}-Eu_2^{(j_2-1)}\|_{a,\Omega}\nonumber\\
&\leq& C_{wu}\eta_a(V_H)\Big(4\theta_2^{j_2}\hat{\lambda}^1C_{\rm up}\eta_1+\big(C_{wu}\eta_a(V_H)\big)^{j_2-1}\|u_2^{(0)}-Eu_2^{(0)}\|_{a,\Omega}\Big)\nonumber\\
&\leq&\frac{4C_{wu}\hat{\lambda}^1C_{\rm up}}{\theta_2}\eta_a(V_H)\eta(\bar{f}_2^{(j_2)};\mathcal{T}_2) +\big(C_{wu}\eta_a(V_H)\big)^{j_2}\|u_2^{(0)}-Eu_2^{(0)}\|_{a,\Omega}.
\end{eqnarray}
Note that $\|\bar{f}_1\|_{a,\Omega}=\|u_1\|_{a,\Omega}=1$ and
$\|\bar{f}_2^{(j_k)}\|_{a,\Omega}=\|\tilde{u}_2^{(j_k-1)}\|_{a,\Omega}=1$.
Combining (\ref{eqn_K_norm}) and (\ref{eqn_lowerbound_Kunl})
leads to the following bounds for $\|K\bar{f}_1\|_{a,\Omega}$ and $\|K\bar{f}_2^{(j_k)}\|_{a,\Omega}$
\begin{eqnarray}\label{eqn_upper_lower_bound_for_Kbarfn+1}
\frac{1}{2\hat{\lambda}^1}\leq \|K\bar{f}_1\|_{a,\Omega}
\leq \frac{\|\bar{f}_1\|_{a,\Omega}}{c_a^2}= \frac{1}{c_a^2},\ \ \
\frac{1}{2\hat{\lambda}^1} \leq \|K\bar{f}_2^{(j_k)}\|_{a,\Omega}\leq\frac{1}{c_a^2}.
\end{eqnarray}
By a similar procedure as (\ref{eqn_update_decision_condition_k=2}), we can derive
\begin{eqnarray}\label{eqn__update_decision_condition_k=2_inverse}
\eta(f_2^{(0)};\mathcal{T}_2)\leq \theta_2 \eta(f_1;\mathcal{T}_1).
\end{eqnarray}
From (\ref{eqn_split_u_kj-Eu_kj}), (\ref{eqn_the_first_two_step_8}), (\ref{eqn_upper_lower_bound_for_Kbarfn+1}), (\ref{eqn__update_decision_condition_k=2_inverse}) and Lemma \ref{lemma_upper_lower_bound_boundaryvalue}, we have
\begin{eqnarray}\label{eqn_estimate_w2-Ew2_by_eta2_10}
&&\|u_2^{(0)}-Eu_2^{(0)}\|_{a,\Omega} \leq \|(I-R_2)w_2^{(0)}\|_{a,\Omega}+\|w_2^{(0)}-Ew_2^{(0)}\|_{a,\Omega}\nonumber\\
&\leq& C_{\rm up} \eta(f_2^{(0)};\mathcal{T}_2) + C_{w,1}\eta_a(V_H) \eta(f_1;\mathcal{T}_1)
\leq \big(C_{\rm up}\theta_2+C_{w,1}\eta_a(V_H)\big)\eta(f_1;\mathcal{T}_1)\nonumber\\
&\leq& 2\big(C_{\rm up}\theta_2+C_{w,1}\eta_a(V_H)\big)\hat{\lambda}^1\eta(\bar{f}_1;\mathcal{T}_1).
\end{eqnarray}
Thus, combining (\ref{eqn_decision_condition_after_break}) and (\ref{eqn_estimate_w2-Ew2_by_eta2_10}) leads to following inequalities
\begin{eqnarray}\label{eqn_estimate_w2-Ew2_by_eta2_2}
&&\big(C_{wu}\eta_a(V_H)\big)^{j_2}\|u_2^{(0)}-Eu_2^{(0)}\|_{a,\Omega}\nonumber\\
&\leq&2\big(C_{\rm up}\theta_2+C_{w,1}\eta_a(V_H)\big)\hat{\lambda}^1\frac{C_{wu}\eta_a(V_H)}{\theta_2^2}\Big(\frac{C_{wu}\eta_a(V_H)}{
\theta_2}\Big)^{j_2-1}\eta(\bar{f}_2^{(j_2)};\mathcal{T}_2)\nonumber\\
&\leq&\frac{4C_{\rm up}\hat{\lambda}^1C_{wu}}{\theta_2}\eta_a(V_H)\eta(\bar{f}_2^{(j_2)};\mathcal{T}_2),
\end{eqnarray}
where we use the fact that the mesh size $H$ is  small enough such that $C_{w,1}\eta_a(V_H)\leq C_{\rm up}\theta_2$
and $2C_{wu}\eta_a(V_H)\leq \theta_2$.
From (\ref{eqn_estimate_w2-Ew2_by_eta2_1}), (\ref{eqn_upper_lower_bound_for_Kbarfn+1}) and (\ref{eqn_estimate_w2-Ew2_by_eta2_2}),
the following estimates hold
\begin{eqnarray}\label{eqn_proof_estimate_of_w-Ew_by_etak_for_k=2_case2}
\|w_2^{(j_2)}-Ew_2^{(j_2)}\|_{a,\Omega}\leq \frac{8C_{wu}\hat{\lambda}^1C_{\rm up}}{\theta_2}\eta_a(V_H)\eta(\bar{f}_2^{(j_2)};\mathcal{T}_2)
\leq\frac{8C_{wu}\hat{\lambda}^1C_{\rm up}}{\theta_2c_a^2}\eta_a(V_H)\eta(f_2^{(j_2)};\mathcal{T}_2).
\end{eqnarray}
The combination of (\ref{eqn_proof_estimate_of_w-Ew_by_etak_for_k=2_case1}) and (\ref{eqn_proof_estimate_of_w-Ew_by_etak_for_k=2_case2})
leads to the desired result (\ref{eqn_estimate_of_w-Ew_by_etak_for_k=2}).
\end{proof}
%--------------------------------------------------------------------------------------------------------
\begin{remark}\label{remark_estimate_wk-Ewk_with_etanell}
The above proof indicates that the keys to deduce the result (\ref{eqn_estimate_of_w-Ew_by_etak_for_k=2}) are
the property (\ref{eqn_the_first_two_step_8}) and the bounds of $\|K\bar{f}_1\|_{a,\Omega}$, $\|K\bar{f}_2^{(0)}\|_{a,\Omega}$ and $\|K\bar{f}_2^{(j_2)}\|_{a,\Omega}$. In what follows, we will prove a similar property as (\ref{eqn_the_first_two_step_8})
for general \revise{cases}, namely, the estimate of $\|w_{k}^{(0)}-Ew_{k}^{(0)}\|_{a,\Omega}$ is controlled by the estimate $\eta_a(V_H)\eta(f_{n_{\ell}};\mathcal{T}_{n_{\ell}})$  for $n_{\ell} < k \leq n_{\ell+1}$. The role of parameter $\theta_2$ in the
{\bf Adaptive Algorithm $C$} is to give an uniform lower bound $\theta_2\eta(f_{n_{\ell}};\mathcal{T}_{n_{\ell}})$ of the error estimators
$\eta(f_{k};\mathcal{T}_{k})$,
whenever $n_{\ell}<k< n_{\ell+1}$. As for $n_{\ell+1}$-th level,
the eigenvalue problem solving is designed to reduce the ratio of nonlinear error $\|w_{n_{\ell}+1}^{(j)}-Ew_{n_{\ell}+1}^{(j)}\|_{a,\Omega}$
to discretization error $\|w_{n_{\ell}+1}^{(j)}-R_{n_{\ell}+1}w_{n_{\ell}+1}^{(j)}\|_{a,\Omega}$ to a small value and reset the reference indicator as $\eta_{\ell+1}:=\eta(\bar{f}_{n_{\ell+1}};\mathcal{T}_{n_{\ell+1}})$.
As a result,
we will prove that $\|w_{k}-Ew_{k}\|_{a,\Omega}$ is a high order term compared
with $\eta(f_{k};\mathcal{T}_{k})$ for $k\in \mathbb{N}$. This is the key property of {\bf Adaptive Algorithm $C$}, and also the reason why the
adaptive refinement can be designed with respect to the a posterior error estimate of
the corresponding linear boundary value problem.
\end{remark}
%-----------------------------------------------------------------------------------------------------------
%\subsection{Estimate of $\|(I-R_{n+2})w_{n+2}\|_{a,\Omega}$}
The following lemmas will be used in our analysis to estimate $\|(I-R_{2})w_{2}\|_{a,\Omega}$.
\begin{lemma}(\cite{DaiXuZhou})\label{lemma_etimate_osc_anorm}
There exists a constant $C_{oa}$ which depends on $A$, regularity constant $\gamma^*$
and coefficient $\varphi$, such that
\begin{eqnarray*}
{\rm osc}(Lv_k;\mathcal{T}_k) \leq C_{oa} \|v_k\|_{a,\Omega},\ \ \ \forall v_k \in V_k.
\end{eqnarray*}
\end{lemma}
%-----------------------------------------------------------------------------------------------------------
\begin{lemma}\label{lemma_eta_triangle_inequ}
%For $f \in L^2(\Omega)$ and $g \in L^2(\Omega)$,
The following inequality holds
\begin{eqnarray}\label{Inequality_11}
\eta(f;\mathcal{T}_{k})-\eta(g;\mathcal{T}_{k}) \leq\eta(f-g;\mathcal{T}_{k}),\ \ \
\forall \revise{f,\ g\in L^2(\Omega)}.
\end{eqnarray}
\end{lemma}
%-----------------------------------------------------------------------------------------------------------
\begin{proof}
For any element $T\in \mathcal{T}_k$, we have
\begin{eqnarray*}
&&\eta_k(f,R_{k}Kf;T)-\eta_k(g,R_{k}Kg;T)\nonumber\\
&=&\Big(h_T^2\|f-LR_{k}Kf\|_{0,T}^2
 +\sum\limits_{E\in \mathcal{E}_k,E\subset \partial T}h_E\big\|[[A\nabla (R_{k}Kf)]]_E\cdot \nu_E\big\|^2_{0,E}\Big)^{\frac{1}{2}}\nonumber\\
 &&-\Big(h_T^2\|g-LR_{k}Kg\|_{0,T}^2
 +\sum\limits_{E\in \mathcal{E}_k,E\subset \partial T}h_E\big\|[[A\nabla (R_{k}Kg)]]_E\cdot \nu_E\big\|^2_{0,E}\Big)^{\frac{1}{2}}\nonumber\\
 &\leq&\left(h_T^2\Big(\|f-LR_{k}Kf\|_{0,T}-\|g-LR_{k}Kg\|_{0,T}\Big)^2\right.\nonumber\\
 &&\left.+\sum\limits_{E\in \mathcal{E}_k,E\subset \partial T}h_E\Big(\big\|[[A\nabla (R_{k}Kf)]]_E
 \cdot \nu_E\big\|_{0,E}-\big\|[[A\nabla (R_{k}Kg)]]_E\cdot \nu_E\big\|_{0,E}\Big)^2\right)^{\frac{1}{2}}\nonumber\\
 &\leq&\eta_k\big(f-g,R_{k}K(f-g);T\big).
\end{eqnarray*}
The previous inequality leads to the following inequalities
\begin{eqnarray*}
\eta(f;\mathcal{T}_k)-\eta(g;\mathcal{T}_k)
&=&\left(\sum_{T\in \mathcal{T}_k}\eta_k^2(f,R_{k}Kf;T)\right)^{1/2}
-\left(\sum_{T\in \mathcal{T}_k}\eta_k^2(g,R_{k}Kg;T)\right)^{1/2}\nonumber\\
&\leq&\left(\sum_{T\in \mathcal{T}_k}
\Big(\eta_k(f,R_{k}Kf;T)-\eta_k(g,R_{k}Kg;T)\Big)^2\right)^{1/2}\nonumber\\
&\leq&\left(\sum_{T\in \mathcal{T}_k}\eta_k^2\big(f-g,R_{k}K(f-g);T\big)\right)^{1/2}
=\eta(f-g;\mathcal{T}_{k}).
\end{eqnarray*}
This is the desired result (\ref{Inequality_11}) and the proof is complete.
\end{proof}
%----------------------------------------------------------------------------------------------------------------------
\begin{lemma}\label{lemma_relation _eta_2l}
%For any function $f \in V_k$ and $g \in V_k$,
The following inequality holds
\begin{eqnarray}\label{Inequality_12}
\eta(f;\mathcal{T}_{k})-\eta(g;\mathcal{T}_{k})
\leq\frac{1+C_{oa}}{C_{\rm low}}\|K(f-g)\|_{a,\Omega},\ \ \
\revise{\forall f,\ g\in V_k.}
\end{eqnarray}
\end{lemma}
%----------------------------------------------------------------------------------------------------------------------
\begin{proof}
From Lemmas \ref{lemma_upper_lower_bound_boundaryvalue}, \ref{lemma_etimate_osc_anorm},
\ref{lemma_eta_triangle_inequ} and Remark \ref{remark_fk_vanish_in_osc}, we have
\begin{eqnarray*}
\eta(f;\mathcal{T}_k)-\eta(g;\mathcal{T}_k)\leq\eta(f-g;\mathcal{T}_{k})&\leq& \frac{1}{C_{\rm low}}\Big(\|(R_{k}-I)K(f-g)\|_{a,\Omega}
+{\rm osc}\big(LR_{k}K(f-g);\mathcal{T}_{k}\big)\Big)\nonumber\\
&\leq&\frac{1+C_{oa}}{C_{\rm low}}\|K(f-g)\|_{a,\Omega},
\end{eqnarray*}
which is the desired result (\ref{Inequality_12}) and the proof is complete.
\end{proof}
%--------------------------------------------------------------------------------------------------------------
To estimate $\|(I-R_{2})w_{2}\|_{a,\Omega}$, we establish a relation  between the finite element error
estimates  of two boundary value problems with the right hand side terms $f_{1}$ and $f_{2}$, respectively.
%In this paper, we choose $\alpha$ satisfying $0<\hat{\alpha}< \alpha <1$ where $\hat{\alpha}$ is introduced
%in Lemma \ref{Adaptive_Convergence_Source_Theorem}.
%We come to prove the following contraction property.
% being the right hand term.
%--------------------------------------------------------------------------------------------------------------
\begin{lemma}\label{lemma_error_reduction_2}
The following relation  holds
\begin{eqnarray}\label{eqn_error_reduction_2}
\|(I-R_{2})w_{2}\|_{a,\Omega}^2+\gamma \eta^2(f_{2};\mathcal{T}_{2})
\leq \alpha_1^2\big(\|(I-R_{1})w_{1}\|_{a,\Omega}^2+\gamma\eta^2(f_{1};\mathcal{T}_{1})\big),
\end{eqnarray}
where
\begin{eqnarray}\label{definition_alpha_1}
\alpha_1^2 = \frac{(1+\delta)\hat{\alpha}^2+\frac{8C_{w,1}^2C_{u\eta}}{\gamma}\Big(1+\frac{1}{\delta}\Big)
\eta_a^2(V_H)}{1-\frac{8C_{w,2}^2C_{u\eta}}{\gamma}\Big(1+\frac{1}{\delta}\Big)\eta_a^2(V_H)}
\end{eqnarray}
with
\begin{eqnarray}\label{definition_C_ueta}
\revise{\delta > 0,\ \ \ C_{u\eta} =1+\gamma\Big(\frac{1+C_{oa}}{C_{\rm low}}\Big)^2.}
\end{eqnarray}
\end{lemma}
%--------------------------------------------------------------------------------------------------------------
\begin{proof}
First, from the definition of $R_{2}$ and triangle inequality, the following inequalities hold
\begin{eqnarray}\label{error_reduction_proof_1}
\|(I-R_{2})w_{2}\|_{a,\Omega}
&\leq& \|(I-R_{2})w_{1}\|_{a,\Omega}
+\|(I-R_{2})(w_{1}-w_{2})\|_{a,\Omega}\nonumber\\
&\leq& \|(I-R_{2})w_{1}\|_{a,\Omega} + \|w_{1}- w_{2}\|_{a,\Omega}.
\end{eqnarray}
From Lemma \ref{lemma_relation _eta_2l}, we have
\begin{eqnarray}\label{error_reduction_proof_4}
\eta(f_{2};\mathcal{T}_{2})
&\leq&\eta(f_{1};\mathcal{T}_{2})+\frac{1+C_{oa}}{C_{\rm low}}\|w_{1}- w_{2}\|_{a,\Omega}.
\end{eqnarray}
Combining (\ref{error_reduction_proof_1}) and (\ref{error_reduction_proof_4}) leads to
\begin{eqnarray}\label{eqn_reduction_lvl_1to2_proof_1}
\|(I-R_{2})w_{2}\|_{a,\Omega}^2+\gamma \eta^2(f_{2};\mathcal{T}_{2}) &\leq& (1+\delta)\big(\|(I-R_{2})w_{1}\|_{a,\Omega}^2+\gamma\eta^2(f_{1};\mathcal{T}_{2})\big)\nonumber\\
&&\ \ +\Big(1+\frac{1}{\delta}\Big)C_{u\eta}\|w_{1}- w_{2}\|_{a,\Omega}^2,
\end{eqnarray}
where $\delta>0$ and $C_{u\eta}$ is defined by (\ref{definition_C_ueta}).
For the last term of (\ref{eqn_reduction_lvl_1to2_proof_1}), by triangle inequality
and Lemma \ref{lemma_equality_of_2kind_error}, we have
\begin{eqnarray}\label{eqn_error_reduction_proof_5}
\|w_{1}- w_{2}\|_{a,\Omega} &\leq& \|w_{1}-\hat{u}^1\|_{a,\Omega}
+\|\hat{u}^1-w_{2}\|_{a,\Omega}\nonumber\\
&\leq& 2\|w_{1}-Ew_{1}\|_{a,\Omega}+2 \|w_{2}-Ew_{2}\|_{a,\Omega}.
\end{eqnarray}
Then combining (\ref{eqn_the_first_two_step_10}), (\ref{eqn_the_first_two_step_9})
and (\ref{eqn_error_reduction_proof_5}) leads to
\begin{eqnarray}\label{error_reduction_proof_2}
\|w_{1}-w_{2}\|_{a,\Omega}^2&\leq& \big(2C_{w,1}\eta_a(V_H)\eta(f_{1};\mathcal{T}_{1})+2C_{w,2}\eta_a(V_H)\eta(f_{2};\mathcal{T}_{2})\big)^2\nonumber\\
&\leq&8C_{w,1}^2\eta_a^2(V_H)\eta^2(f_{1};\mathcal{T}_{1})+8C_{w,2}^2\eta_a^2(V_H)\eta^2(f_{2};\mathcal{T}_{2}).
\end{eqnarray}
%From (\ref{error_reduction_proof_1}) and (\ref{error_reduction_proof_2}), we obtain
%\begin{eqnarray}\label{error_reduction_proof_3}
%&&\|(I-R_{2})w_{2}\|_{a,\Omega}^2 \leq (1+\delta)\|(I-R_{2})w_{1}\|_{a,\Omega}^2\nonumber\\
%&&\ \ \ \ \ \ \ \ \ \ +8\Big(1+\frac{1}{\delta}\Big)\eta_a^2(V_H)\big(C_{w,1}^2\eta^2(f_{1};\mathcal{T}_{1})
%+C_{w,2}^2\eta^2(f_{2};\mathcal{T}_{2})\big),
%\end{eqnarray}
%
%From (\ref{error_reduction_proof_2}) and (\ref{error_reduction_proof_4}), we obtain
%\begin{eqnarray}\label{error_reduction_proof_6}
%&&\eta^2(f_{2};\mathcal{T}_{2})\leq(1+\delta)\eta^2(f_{1};\mathcal{T}_{2})\nonumber\\
%&&\ \ \ +8\big(\frac{1+C_{oa}}{C_{\rm low}}\big)^2\Big(1+\frac{1}{\delta}\Big)\eta_a^2(V_H)\big(C_{w,1}^2\eta^2(f_{1};\mathcal{T}_{1})
%+C_{w,2}^2\eta^2(f_{2};\mathcal{T}_{2})\big).
%\end{eqnarray}
From (\ref{eqn_reduction_lvl_1to2_proof_1}) and (\ref{error_reduction_proof_2}), the following inequalities hold
\begin{eqnarray}\label{eqn_relation _errorestimate_2boundvalueproblems}
&&\Big(1-\frac{8C_{w,2}^2C_{u\eta}}{\gamma}\Big(1+\frac{1}{\delta}\Big)\eta_a^2(V_H)\Big)\big(\|(I-R_{2})w_{2}\|_{a,\Omega}^2+\gamma \eta^2(f_{2};\mathcal{T}_{2})\big)\nonumber\\
&\leq&\|(I-R_{2})w_{2}\|_{a,\Omega}^2+\gamma \eta^2(f_{2};\mathcal{T}_{2})-8C_{u\eta}\Big(1+\frac{1}{\delta}\Big)\eta_a^2(V_H)C_{w,2}^2\eta^2(f_{2};\mathcal{T}_{2})\nonumber\\
&\leq&(1+\delta)\big(\|(I-R_{2})w_{1}\|_{a,\Omega}^2+\gamma\eta^2(f_{1};\mathcal{T}_{2})\big)
+8C_{u\eta}\Big(1+\frac{1}{\delta}\Big)\eta_a^2(V_H)C_{w,1}^2\eta^2(f_{1};\mathcal{T}_{1}).
\end{eqnarray}
Using Theorem \ref{Adaptive_Convergence_Source_Theorem} and the definitions
\begin{eqnarray*}
f_{1}=\frac{\bar{f}_{1}}{\|K\bar{f}_{1}\|_{a,\Omega}}, \ \ \ \
w_{1} = \frac{K\bar{f}_{1}}{\|K\bar{f}_{1}\|_{a,\Omega}},
\end{eqnarray*}
we have
\begin{eqnarray}\label{error_reduction_proof_7}
\|(I-R_{2})w_{1}\|_{a,\Omega}^2+\gamma\eta^2(f_{1};\mathcal{T}_{2})
\leq \hat{\alpha}^2\big(\|(I-R_{1})w_{1}\|_{a,\Omega}^2+\gamma\eta^2(f_{1};\mathcal{T}_{1})\big),
\end{eqnarray}
where $\hat{\alpha} \in (0,1)$.
Combining (\ref{definition_alpha_1}), (\ref{eqn_relation _errorestimate_2boundvalueproblems})
and (\ref{error_reduction_proof_7}) leads to the desired result (\ref{eqn_error_reduction_2}) %as follows
%% following inequalities
%\begin{eqnarray*}
%\|(I-R_{2})w_{2}\|_{a,\Omega}^2+\gamma \eta^2(f_{2};\mathcal{T}_{2})
%\leq \alpha_1^2\big(\|(I-R_{1})w_{1}\|_{a,\Omega}^2+\gamma\eta^2(f_{1};\mathcal{T}_{1})\big),
%\end{eqnarray*}
%where we choose $\delta$ and $H$ small enough such that
%\begin{eqnarray*}
%\alpha^2 >\frac{(1+\delta)\hat{\alpha}^2+\frac{8C_{w,1}^2C_{u\eta}}{\gamma}\Big(1+\frac{1}{\delta}\Big)\eta_a^2(V_H)}
%{1-\frac{8C_{w,2}^2C_{u\eta}}{\gamma}\Big(1+\frac{1}{\delta}\Big)\eta_a^2(V_H)}.
%\end{eqnarray*}
%which is the desired result (\ref{eqn_error_reduction_2}) and
and the proof is complete.
\end{proof}
\revise{By the definition of $\alpha_1$ in (\ref{eqn_error_reduction_2})}, we can choose $\delta$ and $H$ small enough such that
\begin{eqnarray*}
\hat{\alpha}<\alpha_1 < 1.
\end{eqnarray*}
Furthermore, for a fixed constant $\alpha$ satisfying $\hat{\alpha}<\alpha<1$, we can choose $\delta$ and $H$
small enough such that
\begin{eqnarray*}
\hat{\alpha}<\alpha_1 < \alpha <1.
\end{eqnarray*}
From Lemma \ref{lemma_error_reduction_2}, when $H$ is small enough, we have
\begin{eqnarray}\label{Convergence_Line_2}
\|(I-R_{2})w_{2}\|_{a,\Omega}^2+\gamma \eta^2(f_{2};\mathcal{T}_{2})
\leq \alpha^2\big(\|(I-R_{1})w_{1}\|_{a,\Omega}^2+\gamma\eta^2(f_{1};\mathcal{T}_{1})\big).
\end{eqnarray}
%--------------------------------------------------------------------------------------------------------------
\subsection{Induction procedure}\label{section_induction_procedure}
In this subsection, the convergence for the general levels will be deduced.
These results are provided as the preparation for the analysis of the convergence by the induction
method for {\bf Adaptive Algorithm $C$} proposed in this paper.
For this aim, we define the following notation
%\begin{eqnarray}
%&&C_{w,1}=2C_{wu}C_{\rm up},\ \ \ C_{w,2}=\frac{C_{ww}\alpha(C_{\rm up}^2
%+\gamma)^{\frac{1}{2}}\hat{\lambda}^2}{(1-\alpha)c_a^2\theta_2}
%+\frac{C_{w,1}\hat{\lambda}^2}{c_a^2\theta_2},\\
%&&C_w = \max\{C_{w,0},C_{w,1},C_{w,2}\}, \ \ \ \hat{\alpha}<\alpha<1.
%&&C_{\mathcal{M}}=C_{\mathcal{M},0},\ \ C_c=C_{\#}^{2s}C_{\mathcal{M}}^{2s}
%\Big(1+\frac{2\gamma}{C_{\rm low}^2}\Big)\alpha^2(1-\alpha^{\frac{1}{2}})^{-2s}
%\Big(2+\frac{5}{\gamma}\Big),\\
%&&C_m=\frac{\hat{C}_{\#}^sC_c^{\frac{1}{2}}}{m},\ \ \
%C_{\eta,1} = \frac{4C_m(C_{\rm up}+2)}{C_{\rm low}},\ \ \
%C_{\eta} =\max\{1,C_{\eta,0},C_{\eta,1}\}.
%\end{eqnarray}
\begin{eqnarray}
\bar{C}_{w}^{(0)} &=&\frac{C_{ww}\alpha(C_{\rm up}^2+\gamma)^{\frac{1}{2}}
}{(1-\alpha)}+2C_{wu}C_{\rm up},\label{definition_Cw2}\\
\bar{C}_w &=&\frac{2\bar{C}_{w}^{(0)}\hat{\lambda}^1}{\theta_2c_a^2},\\
C_w &=& \max\Big\{C_{w,1},C_{w,2},\bar{C}_w\Big\}.\label{definition_Cw}
\end{eqnarray}
For the induction procedure, \revise{from (\ref{Convergence_Line_2}),
Lemmas \ref{lemma_estimate_of_w-Ew_by_etal_for_k=2} and \ref{lemma_estimate_of_w-Ew_by_etak_for_k=2}}
in the previous section, it is natural to suppose that
\begin{eqnarray}
\|(I-R_{k+1})w_{k+1}\|_{a,\Omega}^2 + \gamma \eta^2(f_{k+1};\mathcal{T}_{k+1})\leq \alpha^2\Big(\|(I-R_k)w_{k}\|_{a,\Omega}^2
+ \gamma \eta^2(f_{k};\mathcal{T}_{k})\Big)\label{eqn_induction_error_reduction}
%\#\mathcal{M}_{k} \leq C_{\mathcal{M}}\big(\|(I-R_k)w_k\|_{a,\Omega}^2
%+{\rm osc}^2(LR_kw_{k};\mathcal{T}_{k})\big)^{-\frac{1}{2s}}\label{eqn_induction_bound_num_Mk}
\end{eqnarray}
holds for $1 \leq k\leq n$, and
\begin{eqnarray}
&&\|w_k-Ew_k\|_{a,\Omega} \leq C_w \eta_a(V_H)\eta(f_k;\mathcal{T}_k)\label{eqn_induction_wk-Ewk}
%&&\big(\|u_{k}-Eu_k\|_{a,\Omega}^2 + {\rm osc}^2(Lu_{k};\mathcal{T}_{k})\big)
%(\#\mathcal{T}_{k}-\#\mathcal{T}_{1})^{2s} \leq C_{c}\label{eqn_induction_bound_{n+1}um_Tk-T1}
\end{eqnarray}
holds for $1\leq k \leq n+1$.

Since $C_{w,1}\leq C_w$ and $C_{w,2}\leq C_w$, the properties (\ref{eqn_induction_error_reduction})-(\ref{eqn_induction_wk-Ewk}) hold
for $n=1$ from Lemmas \ref{lemma_estimate_of_w-Ew_by_etal_for_k=2} and \ref{lemma_estimate_of_w-Ew_by_etak_for_k=2},
and (\ref{Convergence_Line_2}). Now we come to prove  (\ref{eqn_induction_error_reduction}) holds for $k=n+1$
and (\ref{eqn_induction_wk-Ewk}) holds for $k=n+2$. In this subsection, we suppose $n_{\ell}<n+2\leq n_{\ell+1}$
for some $\ell\in \mathbb{N}$.

%-----------------------------------------------------------------------------------------------
%\subsection{Estimate of $\|w_{n+2}-Ew_{n+2}\|_{a,\Omega}$}
%-----------------------------------------------------------------------------------------------
\begin{lemma}
\revise{When the mesh size $H$ and $h_1$ are small enough}, the following bounds hold
\begin{eqnarray}\label{eqn_lowerbound_of_Kuk_Kwk}
\|Ku_k\|_{a,\Omega} \geq \frac{1}{2\hat{\lambda}^1}\ \ \
\text{and}\ \ \ \|Kw_{k}\|_{a,\Omega} \geq \frac{1}{2\hat{\lambda}^1},
\ \ \ {\rm for}\ k\leq n+1.
\end{eqnarray}
\end{lemma}
%-----------------------------------------------------------------------------------------------
\begin{proof}
By Theorem \ref{theorem_split_u_k-Eu_k}, the following inequality holds
\begin{eqnarray}\label{eqn_proof_C_w_4}
\|u_{k}-Eu_{k}\|_{a,\Omega}\leq \|(I-R_{k})w_{k}\|_{a,\Omega}+\|w_{k}-Ew_{k}\|_{a,\Omega}.
\end{eqnarray}
Since \revise{the mesh size $H$ is small enough, (\ref{eqn_induction_error_reduction}) holds for $k\leq n$ and (\ref{eqn_induction_wk-Ewk})
holds for $k\leq n+1$}, we can derive
\begin{eqnarray}\label{eqn_proof_C_w_5}
&&\|u_{k}-Eu_{k}\|_{a,\Omega}\leq \|(I-R_{k})w_{k}\|_{a,\Omega}
+C_w\eta_a(V_H)\eta(f_{k};\mathcal{T}_{k})\nonumber\\
&&\leq\sqrt{2}\Big(\|(I-R_{k})w_{k}\|_{a,\Omega}^2
+\gamma\eta^2(f_{k};\mathcal{T}_{k})\Big)^{\frac{1}{2}}
\leq\sqrt{2}\Big(\|(I-R_{1})w_1\|_{a,\Omega}^2
+\gamma\eta^2(f_{1};\mathcal{T}_{1})\Big)^{\frac{1}{2}}.
\end{eqnarray}
Since the initial mesh size $h_1$ is small enough,  the following estimate holds
\begin{eqnarray}\label{eqn_assumption_initialmesh_1}
\sqrt{2}\Big(\|(I-R_{1})w_1\|_{a,\Omega}^2+\gamma\eta^2(f_{1};\mathcal{T}_{1})\Big)^{\frac{1}{2}}
\leq \frac{c_a^2}{4\hat{\lambda}^1},
\end{eqnarray}
which combined with (\ref{eqn_proof_C_w_5}) leads to
\begin{eqnarray}\label{Inequality_9}
\|u_{k}-Eu_{k}\|_{a,\Omega}\leq \frac{c_a^2}{4\hat{\lambda}^1}.
\end{eqnarray}
Similarly, the combination of \revise{(\ref{eqn_induction_wk-Ewk})}, (\ref{eqn_proof_C_w_5}) and (\ref{eqn_assumption_initialmesh_1}) implies
\begin{eqnarray}\label{Inequality_10}
\|w_{k}-Ew_{k}\|_{a,\Omega} \leq \frac{c_a^2}{4\hat{\lambda}^1}.
\end{eqnarray}
Then from Lemma \ref{lemma_bound_Kv_anorm}, (\ref{Inequality_9}) and (\ref{Inequality_10}),
we can obtain that the estimate (\ref{eqn_lowerbound_of_Kuk_Kwk}) holds for $k\leq n+1$.
\end{proof}

%------------------------------------------------------------------------------------------------
From Remark \ref{remark_estimate_wk-Ewk_with_etanell}, in order to show (\ref{eqn_induction_wk-Ewk})
holds for $k=n+2$, we need the bounds for $\|K\bar{f}_{n_{\ell}}\|_{a,\Omega}$,
$\|K\bar{f}_{n+2}^{(0)}\|_{a,\Omega}$ and $\|K\bar{f}_{n+2}\|_{a,\Omega}$.

Note that $\bar{f}_{n_{\ell}} = \tilde{u}_{n_{\ell}}^{(j_{n_{\ell}}-1)}$ and
$ \|\tilde{u}_{n_{\ell}}^{(j_{n_{\ell}}-1)}\|_{a,\Omega}=1$.
From (\ref{eqn_K_norm}) and (\ref{eqn_lowerbound_Kunl}), we have the following bound for $\|K\bar{f}_{n_{\ell}}\|_{a,\Omega}$
\begin{eqnarray}\label{eqn_upper_lower_bound_for_Kbarfn+11}
\frac{1}{2\hat{\lambda}^1}\leq \|K\bar{f}_{n_{\ell}}\|_{a,\Omega}
\leq \frac{\|\bar{f}_{n_{\ell}}\|_{a,\Omega}}{c_a^2}= \frac{1}{c_a^2}.
\end{eqnarray}
If $n+2 = n_{\ell}+1$, then $\bar{f}_{n+2}^{(0)} = \tilde{u}_{n_{\ell}}$ and $\|\tilde{u}_{n_{\ell}}\|_{a,\Omega} =1$.
Using (\ref{eqn_K_norm}) and (\ref{eqn_lowerbound_Kunl}) again, we have the following bound for $\|K\bar{f}_{n+2}^{(0)}\|_{a,\Omega}$
\begin{eqnarray}\label{eqn_bound_Kf_n+2}
\frac{1}{2\hat{\lambda}^1}\leq \|K\bar{f}_{n+2}^{(0)}\|_{a,\Omega}
\leq  \frac{1}{c_a^2}.
\end{eqnarray}
Else, if $n_{\ell}+1< n+2 \leq n_{\ell+1}$, we know that $\bar{f}_{n+2}^{(0)} = u_{n+1}$ and $\|u_{n+1}\|_{a,\Omega}=1$.
From (\ref{eqn_K_norm}) and (\ref{eqn_lowerbound_of_Kuk_Kwk}), we still have the result (\ref{eqn_bound_Kf_n+2}).

Now, it remains to show the bound for $\|K\bar{f}_{n+2}\|_{a,\Omega}$. If $n_{\ell}<n+2< n_{\ell+1}$,
we already obtain the same estimate as (\ref{eqn_bound_Kf_n+2}) since $f_{n+2} = f_{n+2}^{(0)}$ in this case.
Else, in the case that $n+2 = n_{\ell+1}$, we know $f_{n+2} = \tilde{u}_{n_{\ell+1}}^{(j_{n_{\ell+1}}-1)}$ and $\|\tilde{u}_{n_{\ell+1}}^{(j_{n_{\ell+1}}-1)}\|_{a,\Omega}=1$.
The following inequality can be deduced by using (\ref{eqn_K_norm}) and (\ref{eqn_lowerbound_Kunl}) again
\begin{eqnarray}\label{New_Notation_2}
\frac{1}{2\hat{\lambda}^1}\leq \|K\bar{f}_{n+2}\|_{a,\Omega}
\leq  \frac{1}{c_a^2}.
\end{eqnarray}

%-----------------------------------------------------------------------------------------------------
\begin{lemma}
When $H$ is small enough, the following estimate holds
\begin{eqnarray}\label{eqn_wn+2-Ewn+2_by_eta_nl}
\|w_{n+2}^{(0)}-Ew_{n+2}^{(0)}\|_{a,\Omega} \leq \bar{C}_{w}^{(0)}\eta_a(V_H)\eta(f_{n_{\ell}};\mathcal{T}_{n_{\ell}}),
\end{eqnarray}
where constant $\bar{C}_w^{(0)}$ is defined by (\ref{definition_Cw2}).
\end{lemma}
\begin{proof}
Since (\ref{eqn_induction_wk-Ewk})
holds for $k=n_{\ell}$, using Lemma \ref{lemma_upper_lower_bound_boundaryvalue} and
Theorem \ref{theorem_split_u_k-Eu_k}, we have
\begin{eqnarray}\label{eqn_proof_C_w_6}
\|u_{n_{\ell}}-Eu_{n_{\ell}}\|_{a,\Omega}
\leq\|(I-R_{n_{\ell}})w_{n_{\ell}}\|_{a,\Omega}+\|w_{n_{\ell}}-Ew_{n_{\ell}}\|_{a,\Omega}\leq\Big(C_{\rm up}+C_w\eta_a(V_H)\Big)\eta(f_{n_{\ell}};\mathcal{T}_{n_{\ell}}).
\end{eqnarray}
According to Lemma \ref{lemma_estimate_of_w-Ew}, (\ref{eqn_lowerbound_Kunl}) and (\ref{eqn_proof_C_w_6}),
the following inequalities hold
\begin{eqnarray}\label{eqn_proof_C_w_7}
&&\|w_{n_{\ell}+1}^{(0)}-Ew_{n_{\ell}+1}^{(0)}\|_{a,\Omega}
\leq C_{wu}\eta_a(V_H)\|u_{n_{\ell}}-Eu_{n_{\ell}}\|_{a,\Omega}\nonumber\\
&\leq&C_{wu}\Big(C_{\rm up}+C_w\eta_a(V_H)\Big)\eta_a(V_H) \eta(f_{n_{\ell}};\mathcal{T}_{n_{\ell}})
\leq 2C_{wu}C_{\rm up}\eta_a(V_H) \eta(f_{n_{\ell}};\mathcal{T}_{n_{\ell}}),
\end{eqnarray}
when $H$ is small enough. From (\ref{eqn_proof_C_w_7}), we already prove (\ref{eqn_wn+2-Ewn+2_by_eta_nl}) for case $n+2=n_{\ell}+1$.

Now we consider the case $n+2>n_{\ell}+1$. For $n_{\ell}+1\leq i<n+2\leq n_{\ell+1}$, we observe that $j_i=0$ and $w_i=w_i^{(0)}$.
Thus, from Lemma \ref{lemma_estimate_of_w-Ew}, (\ref{eqn_lowerbound_of_Kuk_Kwk}) and recursive argument,
$\|w_{n+2}^{(0)}-Ew_{n+2}^{(0)}\|_{a,\Omega}$ has following estimates
\begin{eqnarray}\label{eqn_proof_C_w_3}
&&\|w_{n+2}^{(0)}-Ew_{n+2}^{(0)}\|_{a,\Omega}
\leq C_{ww}\eta_a(V_{n+1} )\|(I-R_{n+1})w_{n+1} \|_{a,\Omega}
+\|w_{n+1}^{(0)} -Ew_{n+1}^{(0)} \|_{a,\Omega}\nonumber\\
&\leq&C_{ww}\sum_{i=n_{\ell}+1}^{n+1}\eta_a(V_{i} )\|(I-R_{i})w_{i}
\|_{a,\Omega}+\|w_{n_{\ell}+1}^{(0)}-Ew_{n_{\ell}+1}^{(0)}\|_{a,\Omega}.
\end{eqnarray}
By using (\ref{eqn_proof_C_w_3}) and the following fact
\begin{eqnarray*}
\eta_a(V_i) \leq \eta_a(V_H) ,\ \ \ \forall i\in\mathbb{N},
\end{eqnarray*}
the following inequalities hold
\begin{eqnarray}\label{eqn_prove_C_w_eta_10}
&&\sum_{i=n_{\ell}+1}^{n+1}\eta_a(V_{i} )\|(I-R_{i})w_{i} \|_{a,\Omega}
\leq \eta_a(V_H)\sum_{i=n_{\ell}+1}^{n+1}\|(I-R_{i})w_{i} \|_{a,\Omega}\nonumber\\
&&\leq\eta_a(V_H)\sum_{i=n_{\ell}+1}^{n+1}\Big(\|(I-R_{i})w_{i} \|_{a,\Omega}^2
+\gamma \eta^2(f_i;\mathcal{T}_i)\Big)^{\frac{1}{2}}.
\end{eqnarray}
Since (\ref{eqn_induction_error_reduction}) holds for $k\leq n$,
from Lemma \ref{lemma_upper_lower_bound_boundaryvalue}, (\ref{definition_w_kj_f_kj}) and (\ref{eqn_prove_C_w_eta_10}),
we have
\begin{eqnarray}\label{prove_C_w_eta_1}
&&\sum_{i=n_{\ell}+1}^{n+1}\eta_a(V_{i} )\|(I-R_{i})w_{i} \|_{a,\Omega}
\leq \eta_a(V_H)\Big(\|(I-R_{n_{\ell}})w_{n_{\ell}} \|_{a,\Omega}^2
+\gamma \eta^2(f_{n_{\ell}};\mathcal{T}_{n_{\ell}})\Big)^{\frac{1}{2}}
\sum_{i=1}^{n+1-n_{\ell}}\alpha^i\nonumber\\
&&\leq \frac{\alpha (C_{\rm up}^2+\gamma)^{\frac{1}{2}}}{1-\alpha}\eta_a(V_H)
\eta(f_{n_{\ell}};\mathcal{T}_{n_{\ell}}).
\end{eqnarray}
Combining (\ref{eqn_proof_C_w_7}), (\ref{eqn_proof_C_w_3}), (\ref{prove_C_w_eta_1}) and the definition (\ref{definition_Cw2})
leads to the desired result (\ref{eqn_wn+2-Ewn+2_by_eta_nl}).
%\begin{eqnarray}\label{prove_C_w_eta_2}
%\|w_{n+2}^{(0)}-Ew_{n+2}^{(0)}\|_{a,\Omega}
%\leq \bar{C}_{w}^{(0)}\eta_a(V_H)\eta(f_{n_{\ell}};\mathcal{T}_{n_{\ell}}).
%%\left(\frac{C_{ww}\alpha(C_{\rm up}^2+\gamma)^{\frac{1}{2}}
%%}{(1-\alpha)}+2C_{wu}C_{\rm up}
%%\right)\eta_a(V_H) \eta(f_{n_{\ell}};\mathcal{T}_{n_{\ell}}).
%\end{eqnarray}
\end{proof}
%------------------------------------------------------------------------------------------------------
So far, we already have  the bounds for $\|K\bar{f}_{n_{\ell}}\|_{a,\Omega}$, $\|K\bar{f}_{n+2}^{(0)}\|_{a,\Omega}$
and $\|K\bar{f}_{n+2}\|_{a,\Omega}$ in (\ref{eqn_upper_lower_bound_for_Kbarfn+11}), (\ref{eqn_bound_Kf_n+2}) and (\ref{New_Notation_2}),
respectively, as well as the estimate (\ref{eqn_wn+2-Ewn+2_by_eta_nl}). Recall Remark \ref{remark_estimate_wk-Ewk_with_etanell}
and the proofs of Lemmas \ref{lemma_estimate_of_w-Ew_by_etak_for_k=2}
and \ref{lemma_error_reduction_2}, the following two lemmas can be deduced by a similar procedure.
%---------------------------------------------------------------------------------------------
\begin{lemma}
When $H$ is small enough, the following estimate holds
\begin{eqnarray}\label{eqn_estimate_wn+2-Ewn+2}
\|w_{n+2}-Ew_{n+2}\|_{a,\Omega} \leq C_w\eta_a(V_H)\eta(f_{n+2};\mathcal{T}_{n+2}),
\end{eqnarray}
where the constant $C_w$ is defined by (\ref{definition_Cw}).
\end{lemma}
%---------------------------------------------------------------------------------------------
\begin{lemma}
When $H$ is small enough, the following relation  holds
\begin{eqnarray}\label{eqn_error_reduction_n+2}
\|(I-R_{n+2})w_{n+2}\|_{a,\Omega}^2+\gamma \eta^2(f_{n+2};\mathcal{T}_{n+2})
\leq \alpha^2\big(\|(I-R_{n+1})w_{n+1}\|_{a,\Omega}^2+\gamma\eta^2(f_{n+1};\mathcal{T}_{n+1})\big).
\end{eqnarray}
\end{lemma}
\subsection{Main results of convergence}
From the discussion in the last two subsections, we can conclude the convergence result of {\bf Adaptive Algorithm $C$}.
%the following convergence results.
\begin{theorem}\label{theorem_main_result}
Let $\{u_k\}_{k\in \mathbb{N}}$ and $\{\mathcal{T}_k\}_{k\in \mathbb{N}}$
be produced by {\bf Adaptive Algorithm $C$},  $w_k$ and $f_k$ be defined by
(\ref{definition_w_k_f_k}). When $H$ is sufficiently small, for $k\in \mathbb{N}$, the following inequalities hold
\begin{eqnarray}
\|w_k-Ew_k\|_{a,\Omega} &\leq& C_w\eta_a(V_H) \eta(f_k;\mathcal{T}_k),\label{eqn_mainresult_estimate_eta}\\
\|(I-R_{k+1})w_{k+1}\|_{a,\Omega}^2+\gamma \eta^2(f_{k+1};\mathcal{T}_{k+1}) &\leq & \alpha^2\big(\|(I-R_k)w_k\|_{a,\Omega}^2+\gamma\eta^2(f_k;\mathcal{T}_{k})\big),
\label{eqn_mainresult_contraction}
%&&\#\mathcal{M}_{k}\leq C_{\mathcal{M}}\big(\|(I-R_{k})w_{k}\|_{a,\Omega}^2
%+{\rm osc}^2(LR_{k}w_{k};\mathcal{T}_{k})\big)^{-\frac{1}{2s}},\label{eqn_mainresult_num_M}\\
%&&\big(\|u_{k}-Eu_{k}\|_{a,\Omega}^2 + {\rm osc}^2(Lu_{k};\mathcal{T}_{k})\big)
%(\#\mathcal{T}_{k}-\#\mathcal{T}_{1})^{2s} \leq  C_{c}.\label{eqn_mainresult_complexity}
\end{eqnarray}
where $\alpha$ satisfies $\hat{\alpha} < \alpha <1$.
\end{theorem}
%-------------------------------------------------------------------------------------------------------
\begin{proof}
Based on the results in Sections \ref{section_property_first_two_levels}
and \ref{section_induction_procedure},  we give the proof by the induction method.

First, combining (\ref{eqn_the_first_two_step_10}), (\ref{eqn_estimate_of_w-Ew_by_etak_for_k=2}) and (\ref{Convergence_Line_2})
in Section \ref{section_property_first_two_levels} leads to that (\ref{eqn_mainresult_estimate_eta})
holds for $k \leq 2$, (\ref{eqn_mainresult_contraction}) holds for $k=1$.

Second, for $n\in \mathbb{N}$, we suppose (\ref{eqn_mainresult_estimate_eta}) holds for $k \leq n+1$ and (\ref{eqn_mainresult_contraction}) holds for $k\leq n$. Then from (\ref{eqn_estimate_wn+2-Ewn+2}), the desired result (\ref{eqn_mainresult_estimate_eta})
holds for $k = n+2$. From (\ref{eqn_error_reduction_n+2}), we obtain
(\ref{eqn_mainresult_contraction}) holds for $k=n+1$.

Finally, we can deduce assertions (\ref{eqn_mainresult_estimate_eta}) and (\ref{eqn_mainresult_contraction}) hold for all $k\in\mathbb{N}$ by the induction method and
the proof is complete.
\end{proof}
%-------------------------------------------------------------------------------------------------------
Now we state a useful Rayleigh quotient expansion for the eigenvalue which is expressed by
the error of the eigenfunction approximation (see \cite{BabuskaOsborn_1989,LinXie}).
%-------------------------------------------------------------------------------------------------------
\begin{lemma}\label{Expansion_Eigenvalue_Lemma}
For any
$w \in V$, we have
\begin{eqnarray*}\label{Expansion_Eigenvalue}
\frac{a(w,w)}{(w,w)}-\hat{\lambda}^1 =\frac{a(w-\hat{u}^1,w-\hat{u}^1)}{(w,w)}
-\frac{\hat{\lambda}^1(w-\hat{u}^1,w-\hat{u}^1)}{(w,w)}.
\end{eqnarray*}
\end{lemma}
%-------------------------------------------------------------------------------------------------------
Using Theorem \ref{theorem_main_result} and Lemma \ref{Expansion_Eigenvalue_Lemma},
we state the convergence property of eigenpair approximation $(\lambda_k,u_k)$ produced
by {\bf Adaptive Algorithm $C$} as follows.
%-------------------------------------------------------------------------------------------------------
\begin{theorem}\label{theorem_convergence_eigenpair}
When $H$ is small enough, there exist constants $C_{\rm con}^u$ and $C_{\rm con}^{\lambda}$ such that
the following inequalities hold
\begin{eqnarray}
\|u_{k}-Eu_{k}\|_{a,\Omega}^2&\leq&  C_{\rm con}^u\alpha^{2(k-1)},\label{eqn_convergence_eigenvector}\\
|\lambda_k-\hat{\lambda}^1| &\leq& C_{\rm con}^{\lambda}\alpha^{2(k-1)},\label{eqn_convergence_eigenvalue}
\end{eqnarray}
where
\begin{eqnarray}\label{New_Definitions}
C_{\rm con}^u = 3\Big(\|w_1-R_1 w_1\|_{a,\Omega}^2+\gamma\eta^2(f_1;\mathcal{T}_1)\Big)\ \ \
{\rm and}\ \ \ C_{\rm con}^{\lambda}=8\hat{\lambda}^1C_{\rm con}^u.
\end{eqnarray}
\end{theorem}
%-------------------------------------------------------------------------------------------------------
\begin{proof}
From Theorems \ref{theorem_split_u_k-Eu_k} and \ref{theorem_main_result}, we have
\begin{eqnarray*}
\|u_{k}-Eu_k\|_{a,\Omega}^2
&\leq& 2\|w_k-R_k w_k\|_{a,\Omega}^2+ 2\|w_k-Ew_k\|_{a,\Omega}^2\nonumber\\
&\leq&2\|w_k-R_k w_k\|_{a,\Omega}^2+2C_{w}^2\eta_a^2(V_H)\eta^2(f_k;\mathcal{T}_k)\nonumber\\
&\leq&\Big(2+\frac{2C_{w}^2\eta_a^2(V_H)}{\gamma}\Big)\Big(\|w_k-R_k w_k\|_{a,\Omega}^2
+\gamma\eta^2(f_k;\mathcal{T}_k)\Big)\nonumber\\
&\leq& 3\alpha^{2(k-1)}\Big(\|w_1-R_1 w_1\|_{a,\Omega}^2+\gamma\eta^2(f_1;\mathcal{T}_1)\Big),
\end{eqnarray*}
when $H$ is small enough. Then the desired assertion (\ref{eqn_convergence_eigenvector}) can be deduced  by
the definition of $C_{\rm con}^u$ in  (\ref{New_Definitions}).
%\begin{eqnarray*}
%C_{\rm con}^u = 3\Big(\|w_1-R_1 w_1\|_{a,\Omega}^2+\gamma\eta^2(f_1;\mathcal{T}_1)\Big).
%\end{eqnarray*}

\revise{From min-max principle, Lemmas \ref{lemma_equality_of_2kind_error} and \ref{Expansion_Eigenvalue_Lemma}}, we have
\begin{eqnarray}\label{eqn_estimate_eigenvalue_1}
0\leq\lambda_k-\hat{\lambda}^1 = \frac{\|u_k-\hat{u}^1\|_{a,\Omega}^2}{\|u_k\|_{0,\Omega}^2}-\frac{\hat{\lambda}^1
\|u_k-\hat{u}^1\|_{0,\Omega}^2}{\|u_k\|_{0,\Omega}^2}\leq \frac{4\|u_k-Eu_k\|_{a,\Omega}^2}{\|u_k\|_{0,\Omega}^2}.
\end{eqnarray}
The mesh size $H$ is chosen to be small enough such that $\|u_k-Eu_k\|_{a,\Omega}^2\leq C_{\rm con}^u\leq 1/8$
holds for all $k\in \mathbb{N}$.
From (\ref{eqn_estimate_eigenvalue_1}) and the fact $\lambda_k\|u_k\|_{0,\Omega}^2 =\|u_k\|_{a,\Omega}^2 = 1$,
the following estimates hold
\begin{eqnarray}\label{New_Notation_3}
\lambda_k-\hat{\lambda}^1 \leq \frac{\lambda_k}{2}\ \ \ \text{or}\ \ \ \lambda_k \leq 2\hat{\lambda}^1 .
\end{eqnarray}
By using (\ref{eqn_estimate_eigenvalue_1}) and (\ref{New_Notation_3}), we can obtain
\begin{eqnarray*}
\lambda_k-\hat{\lambda}^1 \leq 4\lambda_k\|u_k-Eu_k\|_{a,\Omega}^2 \leq 8\hat{\lambda}^1C_{\rm con}^u\alpha^{2(k-1)}.
\end{eqnarray*}
Then the desired assertion (\ref{eqn_convergence_eigenvalue}) can be deduced by defining
$C_{\rm con}^{\lambda}:=8\hat{\lambda}^1C_{\rm con}^u$ and the proof is complete.
\end{proof}
%---------------------------------------------------------------------------------------------------------

\section{Complexity analysis}
%----------------------------------------------------------------------------------------------------------------
In this section, the complexity analysis of {\bf Adaptive Algorithm $C$} will be provided.
In order to state the results of complexity estimate, we need the following assumption.
\begin{assumption}\label{assumption_function_class_u}
There exist an order $s>0$ and a constant $M>0$ such that
\begin{eqnarray*}
M=\sup\limits_{\varepsilon > 0}\varepsilon\inf\limits_{\{\mathcal{T}_{1}\leq\mathcal{T}_{\varepsilon}:
\inf\limits_{u_{\varepsilon}\in V_{\varepsilon},\|u_{\varepsilon}\|_{a,\Omega}=1}
\big(\|u_{\varepsilon}-Eu_{\varepsilon}\|_{a,\Omega}^2+{\rm osc}^2(Lu_{\varepsilon};\mathcal{T}_{\varepsilon})\big)
\leq \varepsilon^2\}}(\#\mathcal{T}_{\varepsilon}-\#\mathcal{T}_{1})^s < \infty,
%\inf\limits_{\varepsilon > 0}\varepsilon\inf\limits_{\{\mathcal{T}_{1}\leq\mathcal{T}_{\varepsilon}:
%\inf\limits_{u_{\varepsilon}\in V_{\varepsilon},\|u_{\varepsilon}\|_{a,\Omega}=1}
%\|u_{\varepsilon}-Eu_{\varepsilon}\|_{a,\Omega}^2+{\rm osc}^2(Lu_{\varepsilon};\mathcal{T}_{\varepsilon})
%\leq \varepsilon^2\}}(\#\mathcal{T}_{\varepsilon}-\#\mathcal{T}_{1})^s = m >0,
\end{eqnarray*}
where $\mathcal{T}_{1}\leq\mathcal{T}_{\varepsilon}$ means $\mathcal{T}_{\varepsilon}$
is a conforming refinement of $\mathcal{T}_{1}$ and $\#\mathcal T$ denotes the number of elements
in the mesh $\mathcal T$.
\end{assumption}
%------------------------------------------------------------------------------------------------------------
%\begin{assumption}\label{assumption_num_refine_mesh}
%There exists a constant $C_{\#}$ such that we have
%\begin{eqnarray}\label{eqn__num_refine_mesh}
%\#\mathcal{T}_{k+1}-\#\mathcal{T}_{k} \leq C_{\#} \#\mathcal{M}_k, \ \ \ {\rm for}\ k\in \mathbb{N}.
%\end{eqnarray}
%%for $k \in \mathbb{N}$.
%\end{assumption}
%----------------------------------------------------------
Now, it is time to consider the complexity of {\bf Adaptive Algorithm $C$}.
First, the number of refined elements has the following upper bound.
\begin{lemma}\label{lemma_bound_num_M}
When $H$ is small enough, the number of elements in $\mathcal{M}_k$ has the following estimate
\begin{eqnarray}\label{eqn_bound_num_M}
\#\mathcal{M}_k\leq C_{\mathcal{M}}\Big(\|(I-R_k)w_k\|_{a,\Omega}^2
+{\rm osc}^2(LR_kw_k;\mathcal{T}_k)\Big)^{-\frac{1}{2s}},
\end{eqnarray}
where
\begin{eqnarray*}
C_{\mathcal{M}}=\Big(\frac{\tilde{\xi}_0^2}{16\hat{C}_D}\Big)^{-\frac{1}{2s}}M^{\frac{1}{s}}.
\end{eqnarray*}
\end{lemma}
%------------------------------------------------------------------------------------------------------
\begin{proof}
We choose $\tilde{\xi}_0$ in Lemma \ref{Error_estimate_Lower_Bound_Lemma}
small enough such that $\tilde{\theta}=\theta_*\sqrt{1-2\tilde{\xi}_0^2}\geq\theta_1$
and $\varepsilon$ to be
\begin{eqnarray}\label{definition_varepsilon_for_k=1}
\varepsilon := \Big(\frac{\tilde{\xi}_0^2}{16\hat{C}_D}\Big)^{\frac{1}{2}}\Big(\|(I-R_k)w_k\|_{a,\Omega}^2
+{\rm osc}^2(LR_kw_k;\mathcal{T}_k)\Big)^{\frac{1}{2}}.
\end{eqnarray}
Let $\mathcal{T}_{{\varepsilon}}$ be a conforming refinement of $\mathcal{T}_1$ with
minimum degrees of freedom satisfying
\begin{eqnarray}\label{Def_u_k_varepsilon}
\|u_{{\varepsilon}}-Eu_{{\varepsilon}}\|_{a,\Omega}^2+
{\rm osc}^2(Lu_{\varepsilon};\mathcal{T}_{{\varepsilon}})
\leq \varepsilon^2,
\end{eqnarray}
where $u_{\varepsilon}\in V_{\varepsilon}$ satisfies $\|u_{\varepsilon}\|_{a,\Omega}=1$.
From Assumption \ref{assumption_function_class_u}, we can get that
\begin{eqnarray}\label{Inequality_k4}
\hskip-0.5cm\#\mathcal{T}_{{\varepsilon}}-\#\mathcal{T}_{1}\leq \varepsilon^{-\frac{1}{s}}M^{\frac{1}{s}}.
\end{eqnarray}
Let $\mathcal{T}_{{k,+}}$ be the smallest common
refinement of $\mathcal{T}_k$ and $\mathcal{T}_{{\varepsilon}}$.
Since both $\mathcal{T}_k$ and $\mathcal{T}_{{\varepsilon}}$
are conforming, $\mathcal{T}_{{k,+}}$ is conforming, and the number of elements in
$\mathcal{T}_{k,+}$ that are not in $\mathcal{T}_k$ is less
than the number of elements that must be added from
$\mathcal{T}_{1}$ to $\mathcal{T}_{\varepsilon}$, i.e.,
\begin{eqnarray*}
\#\mathcal{T}_{k,+}-\#\mathcal{T}_k
\leq \#\mathcal{T}_{\varepsilon}- \#\mathcal{T}_{1}.
\end{eqnarray*}
We conclude from Lemma \ref{quasi optimality of the total error},
there exists a constant $\hat{C}_D$ such that
\begin{eqnarray}\label{eqn_mintotalerror_property_R_k}
&&\|(I-R_{k,+})w_k\|_{a,\Omega}^2
+{\rm osc}^2(LR_{{k,+}}w_k;\mathcal{T}_{k,+})
 \leq \hat{C}_D\Big(\|w_k-u_{{\varepsilon}}\|_{a,\Omega}^2
+{\rm osc}^2(Lu_{\varepsilon};\mathcal{T}_{k,+})\Big)\nonumber\\
&& \leq \hat{C}_D\Big(\|w_k-u_{{\varepsilon}}\|_{a,\Omega}^2
+{\rm osc}^2(Lu_{\varepsilon};\mathcal{T}_{\varepsilon})\Big),
\end{eqnarray}
where $R_{k,+}$ denotes the Galerkin projection onto the finite element space $V_{k,+}$ over $\mathcal{T}_{k,+}$.
By applying Lemma \ref{lemma_equality_of_2kind_error}, Theorem \ref{theorem_main_result}
and triangle inequality, we have the following estimates for the first term of (\ref{eqn_mintotalerror_property_R_k})
\begin{eqnarray}\label{eqn_prove_bound_num_M_k_k}
&&\big\|w_k-u_{{\varepsilon}}\big\|_{a,\Omega}^2\leq 2\|w_k-\hat{u}^1\|_{a,\Omega}^2
+2\|\hat{u}^1-u_{{\varepsilon}}\|_{a,\Omega}^2 \nonumber\\
&\leq& 8\|w_k-Ew_k\|_{a,\Omega}^2+8\|u_{{\varepsilon}}-Eu_{{\varepsilon}}\|_{a,\Omega}^2
\leq 8C_w^2\eta_a^2(V_H)\eta^2(f_k;\mathcal{T}_k)+8\|u_{{\varepsilon}}-Eu_{{\varepsilon}}\|_{a,\Omega}^2.
\end{eqnarray}
The lower bound (\ref{Lower_Bound_Posteriori_Error}) in Lemma \ref{lemma_upper_lower_bound_boundaryvalue}
implies
\begin{eqnarray}\label{eqn_prove_bound_num_M_k_2}
\eta^2(f_k;\mathcal{T}_k) &\leq& \frac{1}{C_{\rm low}^2}\Big(\|(I-R_k)w_k\|_{a,\Omega}
+{\rm osc}(LR_kw_k;\mathcal{T}_k)\Big)^2\nonumber\\
&\leq&\frac{2}{C_{\rm low}^2}\Big(\|(I-R_k)w_k\|_{a,\Omega}^2
+{\rm osc}^2(LR_kw_k;\mathcal{T}_k)\Big).
\end{eqnarray}
From (\ref{definition_varepsilon_for_k=1}), (\ref{Def_u_k_varepsilon}), (\ref{eqn_mintotalerror_property_R_k}),
(\ref{eqn_prove_bound_num_M_k_k}) and (\ref{eqn_prove_bound_num_M_k_2}), when $H$ is small enough,
we obtain
\begin{eqnarray}\label{eqn_prove_bound_num_M_k_3}
&&\|(I-R_{k,+})w_k\|_{a,\Omega}^2
+{\rm osc}^2(LR_{k,+}w_k;\mathcal{T}_{k,+})\nonumber\\
&\leq& \hat{C}_D \Big(8\|u_{{\varepsilon}}-Eu_{{\varepsilon}}\|_{a,\Omega}^2+{\rm osc}^2(Lu_{\varepsilon};\mathcal{T}_{{\varepsilon}})+8C_w^2\eta_a^2(V_H)\eta^2(f_k;\mathcal{T}_k)\Big)\nonumber\\
&\leq&\hat{C}_D \Big(8\varepsilon^2 +\frac{16C_w^2\eta_a^2(V_H)}{C_{\rm low}^2}\big(\|(I-R_k)w_k\|_{a,\Omega}^2
+{\rm osc}^2(LR_kw_k;\mathcal{T}_k)\big)\Big)\nonumber\\
&=&\Big(\frac{\tilde\xi_0^2}{2}+\frac{16\hat{C}_DC_w^2\eta_a^2(V_H)}{C_{\rm low}^2}\Big)\Big(\|(I-R_k)w_k\|_{a,\Omega}^2+{\rm osc}^2(LR_kw_k;\mathcal{T}_k)\Big)\nonumber\\
&\leq& \tilde{\xi}_0^2\Big(\|(I-R_k)w_k\|_{a,\Omega}^2+{\rm osc}^2(LR_kw_k;\mathcal{T}_k)\Big).
\end{eqnarray}
Thus combining Lemma \ref{Error_estimate_Lower_Bound_Lemma} and (\ref{eqn_prove_bound_num_M_k_3})
leads to the following estimate
\begin{eqnarray*}\label{Newspacedorfler}
\eta_k^2\big(f_k,R_kKf_k;\mathcal{T}_k\backslash(\mathcal{T}_{k,+}\cap\mathcal{T}_k)\big)
\geq \theta_1^2\eta_k^2(f_k,R_kKf_k;\mathcal{T}_k).
\end{eqnarray*}
Since $\bar{u}_k=R_kK\bar{f}_k$ and
\begin{eqnarray*}
f_k=\frac{\bar{f}_k}{\|K\bar{f}_k\|_{a,\Omega}},
\end{eqnarray*}
\revise{the following inequality holds}
\begin{eqnarray*}\label{Newspacedorfler2}
\eta_k^2\big(\bar{f}_k,\bar{u}_k;\mathcal{T}_k\backslash(\mathcal{T}_{k,+}
\cap\mathcal{T}_k)\big)
\geq \theta_1^2\eta_k^2(\bar{f}_k,\bar{u}_k;\mathcal{T}_k).
\end{eqnarray*}
Since (\ref{definition_varepsilon_for_k=1}), (\ref{Inequality_k4}) and the marking step selects a minimum set $\mathcal{M}_k$ satisfying
\begin{eqnarray*}
\eta_k^2(\bar{f}_k,\bar{u}_k;\mathcal{M}_k)\geq \theta_1^2\eta_k^2
(\bar{f}_k,\bar{u}_k;\mathcal{T}_k),
\end{eqnarray*}
we know that the value $\#\mathcal M_k$ satisfies
\begin{eqnarray*}\label{eqn_bound_num_M1}
&&\#\mathcal{M}_k  \leq  \#\big(\mathcal{T}_k\backslash(\mathcal{T}_{k,+}\cap\mathcal{T}_k)\big)
\leq \#\mathcal{T}_{k,+}
-\#\mathcal{T}_k  \leq \# \mathcal{T}_{{\varepsilon}} -\#\mathcal{T}_{1}\nonumber\\
&&\leq \Big(\frac{\tilde{\xi}_0^2}{16\hat{C}_D}\Big)^{-\frac{1}{2s}}M^{\frac{1}{s}}
\Big(\|(I-R_k)w_k\|_{a,\Omega}^2
+{\rm osc}^2 (LR_kw_k;\mathcal{T}_k)\Big)^{-\frac{1}{2s}}.
\end{eqnarray*}
This is the desired result (\ref{eqn_bound_num_M}) and the proof is complete.
\end{proof}
%%------------------------------------------------------------------------------------------------------
%%\subsection{Estimate of $\#\mathcal{T}_2-\#\mathcal{T}_1$}
%%In this subsection,
%
%%------------------------------------------------------------------------------------------------
%\begin{proof}

%\end{proof}
%----------------------------------------------------------------------------------------------
%\subsection{Estimate of $\#\mathcal{T}_{n+2}-\#\mathcal{T}_1$}
In our analysis, we also need the following result (see, e.g.,
\cite{CasconKreuzerNochettoSiebert,DaiXuZhou,Nochetto,Stevenson_2007,Stevson_2008}).
\begin{proposition}(\cite[Lemma 2.3]{CasconKreuzerNochettoSiebert})\label{lemma_num_upbound}
 (Complexity of \textsf{REFINE}) Let $\{\mathcal{T}_{k}\}_{k\in \mathbb{N}}$ be a sequence of
conforming nested partitions generated by \textsf{REFINE} starting from $\mathcal{T}_{1}$,
$\mathcal{M}_{k}$ be the set of elements of $\mathcal{T}_{k}$ which
is marked for refinement.
Then there exists a constant $C_{\#}$ depending solely on $\mathcal{T}_{1}$ such that
\begin{eqnarray*}\label{Refinement_bounded}
\#\mathcal{T}_{k+1}- \#\mathcal{T}_{1}
\leq C_{\#}\sum\limits_{i=1}^k\#\mathcal{M}_{i}.
\end{eqnarray*}
\end{proposition}
%----------------------------------------------------------------------------------------------
Now, it is time to show the quasi-optimality by estimating $\#\mathcal{T}_k-\#\mathcal{T}_1$ for $k\in \mathbb{N}$.
\begin{theorem}\label{theorem_complexity}
When $H$ is small enough, for $k\in \mathbb{N}$, the following inequality holds
\begin{eqnarray}\label{eqn_bound_num_Tn+2-Tn+1}
\|u_k-Eu_k\|_{a,\Omega}^2 + {\rm osc}^2(Lu_k;\mathcal{T}_k)\leq  C_{c}(\#\mathcal{T}_k-\#\mathcal{T}_{1})^{-2s},
\end{eqnarray}
where
\begin{eqnarray*}
C_c=C_{\#}^{2s}C_{\mathcal{M}}^{2s}
\Big(1+\frac{2\gamma}{C_{\rm low}^2}\Big)\alpha^2(1-\alpha^{\frac{1}{2}})^{-2s}
\Big(2+\frac{5}{\gamma}\Big).
\end{eqnarray*}
\end{theorem}
%----------------------------------------------------------------------------------------------
\begin{proof}
From (\ref{eqn_prove_bound_num_M_k_2}), we have the following inequality
for $k\in \mathbb{N}$
\begin{eqnarray*}
\|(I-R_k)w_k\|_{a,\Omega}^2+\gamma \eta^2(f_k;\mathcal{T}_k)
\leq\Big(1+\frac{2\gamma}{C_{\rm low}^2}\Big)\Big(\|(I-R_k)w_k\|_{a,\Omega}^2
+{\rm osc}^2(LR_kw_k;\mathcal{T}_k)\Big).
\end{eqnarray*}
Combining Lemma \ref{lemma_bound_num_M} and Proposition \ref{lemma_num_upbound} leads to following estimates
\begin{eqnarray}\label{eqn_complexity_proof_2}
&&\#\mathcal{T}_k-\#\mathcal{T}_{1} \leq C_{\#}\sum_{i=1}^{k-1}\#\mathcal{M}_i
\leq C_{\#}C_{\mathcal{M}}\sum_{i=1}^{k-1}\Big(\|(I-R_i)w_i\|_{a,\Omega}^2
+{\rm osc}^2(LR_iw_i;\mathcal{T}_i)\Big)^{-\frac{1}{2s}}\nonumber\\
&&\leq C_{\#}C_{\mathcal{M}}\Big(1+\frac{2\gamma}{C_{\rm low}^2}\Big)^{\frac{1}{2s}}
\sum_{i=1}^{k-1}\Big(\|(I-R_i)w_i\|_{a,\Omega}^2+\gamma \eta^2(f_i;\mathcal{T}_i)\Big)^{-\frac{1}{2s}}.
\end{eqnarray}
From Theorem \ref{theorem_main_result}, the following inequalities hold
\begin{eqnarray}\label{eqn_complexity_proof_3}
&&\sum_{i=1}^{k-1}\Big(\|(I-R_i)w_i\|_{a,\Omega}^2
+\gamma \eta^2(f_i;\mathcal{T}_i)\Big)^{-\frac{1}{2s}}\leq\Big(\|(I-R_k)w_k\|_{a,\Omega}^2
+\gamma \eta^2(f_k;\mathcal{T}_k)\Big)^{-\frac{1}{2s}}
\sum_{i=1}^{k-1}\alpha^{\frac{i}{s}}\nonumber\\
&&\leq\frac{\alpha^{\frac{1}{s}}}{1-\alpha^{\frac{1}{s}}}
\Big(\|(I-R_k)w_k\|_{a,\Omega}^2
+\gamma \eta^2(f_k;\mathcal{T}_k)\Big)^{-\frac{1}{2s}}.
\end{eqnarray}
The following estimates can be deduced from Theorems \ref{theorem_split_u_k-Eu_k}
and \ref{theorem_main_result},
\begin{eqnarray}\label{eqn_proof_C_c_k=2_1}
\|u_k-Eu_k\|_{a,\Omega}^2
&\leq& 2\|w_k-R_k w_k\|_{a,\Omega}^2+ 2\|w_k-Ew_k\|_{a,\Omega}^2\nonumber\\
&\leq&2\|w_k-R_k w_k\|_{a,\Omega}^2+2C_w^2\eta_a^2(V_H)\eta^2(f_k;\mathcal{T}_k).
\end{eqnarray}
By applying Theorem \ref{theorem_main_result} again, we have the lower bound of $\|R_kw_k\|_{a,\Omega}$
\begin{eqnarray*}\label{eqn_lower_bound_Rkwk}
\|R_kw_k\|_{a,\Omega} \geq \|w_k\|_{a,\Omega} -\|w_k-R_kw_k\|_{a,\Omega}
\geq 1 - \Big(\|w_1-R_1w_1\|_{a,\Omega}^2+\gamma \eta^2(f_1;\mathcal{T}_1)\Big)^{\frac{1}{2}}\geq \frac{1}{2},
\end{eqnarray*}
which implies the following estimate
\begin{eqnarray}\label{eqn_proof_C_c_k=2_2}
{\rm osc}(Lu_k;\mathcal{T}_k)=\frac{{\rm osc}(LR_kw_k;\mathcal{T}_k)}{\|R_kw_k\|_{a,\Omega}}
 \leq 2{\rm osc}(LR_kw_k;\mathcal{T}_k)\leq 2\eta(f_k;\mathcal{T}_k).
\end{eqnarray}
When $H$ is small enough, combining (\ref{eqn_proof_C_c_k=2_1}) and (\ref{eqn_proof_C_c_k=2_2})
leads to the following inequalities
\begin{eqnarray}\label{eqn_complexity_proof_4}
&&\|u_k-Eu_k\|_{a,\Omega}^2 + {\rm osc}^2(Lu_k;\mathcal{T}_k)
\leq 2\|w_k-R_k w_k\|_{a,\Omega}^2+\big(2C_w^2\eta_a^2(V_H)+4\big)\eta^2(f_k;\mathcal{T}_k)\nonumber\\
&&\leq\Big(2+\frac{2C_w^2\eta_a^2(V_H)+4}{\gamma}\Big)\big(\|w_k-R_k w_k\|_{a,\Omega}^2
+\gamma\eta^2(f_k;\mathcal{T}_k)\big)\nonumber\\
&&\leq \Big(2+\frac{5}{\gamma}\Big)\big(\|w_k-R_k w_k\|_{a,\Omega}^2+\gamma\eta^2(f_k;\mathcal{T}_k)\big).
\end{eqnarray}
Then from (\ref{eqn_complexity_proof_2}), (\ref{eqn_complexity_proof_3})
and (\ref{eqn_complexity_proof_4}), we have
\begin{eqnarray*}
(\#\mathcal{T}_k-\#\mathcal{T}_{1}) \leq C_c^{\frac{1}{2s}}
\big(\|u_k-Eu_k\|_{a,\Omega}^2
+ {\rm osc}^2(Lu_k;\mathcal{T}_k)\big) ^{-\frac{1}{2s}}.
\end{eqnarray*}
This is equivalent to the desired result (\ref{eqn_bound_num_Tn+2-Tn+1}) and the proof is complete.
\end{proof}
%-----------------------------------------------------------------------------------------------------
\section{Numerical experiments}
In this section, we investigate the numerical performance of {\bf Adaptive Algorithm $C$} for
the second order elliptic eigenvalue problems by four numerical examples.
Here,  the well known implicitly restarted Lanczos method,
which is included in the package ARPACK, is adopted for solving the concerned generalized eigenvalue problems and the
geometric multigrid (GMG) method is selected as the linear solver for the boundary value problems.
All through these papers, we choose $\mathcal T_H=\mathcal T_1$.

%-----------------------------------------------------------------------------------
{\bf Example 1.}\ In this example, we consider the following eigenvalue problem
\begin{equation}\label{eigenproblem_Exam_1}
\left\{
\begin{array}{rcl}
-\frac{1}{2}\Delta u +\frac{1}{2}|x|^2u&=&\lambda u\ \ \ {\rm in}\ \Omega,\\
u&=&0\ \ \ \ \ {\rm on}\ \partial\Omega,\\
\|u\|_{a,\Omega}&=&1,
\end{array}
\right.
\end{equation}
where $\Omega \subset \mathbb{R}^2$ and $|x|=\sqrt{x_1^2+x_2^2}$.
The first eigenvalue of (\ref{eigenproblem_Exam_1}) is $\lambda=1$ and the associated eigenfunction
is $u=\kappa e^{-|x|^2/2}$ with any nonzero constant $\kappa$.
In our computation, we set $\Omega=(-5, 5)\times (-5, 5)$.

\begin{figure}[ht]
\centering
\includegraphics[width=5.5cm,height=5.5cm]{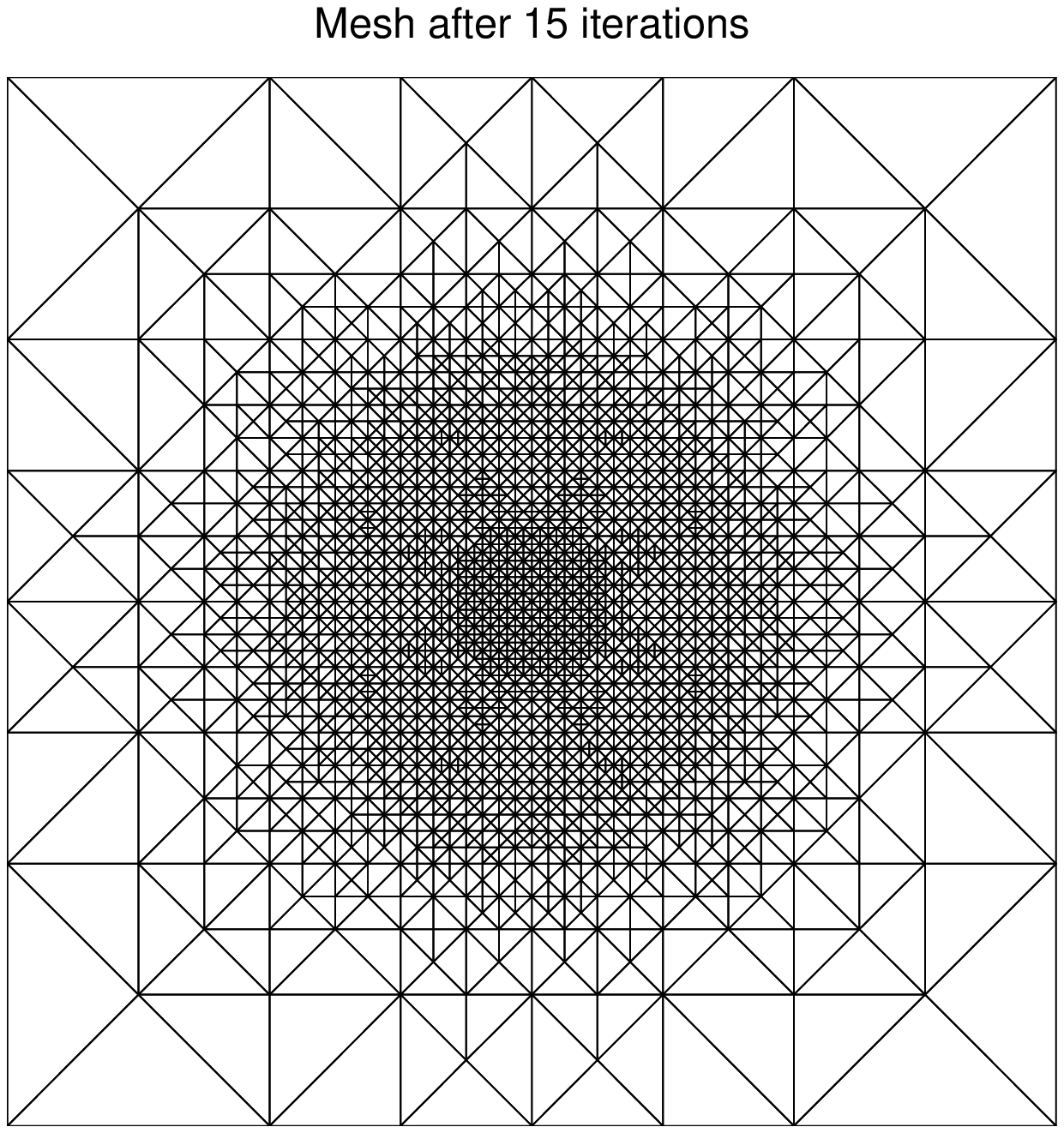}
\includegraphics[width=5.5cm,height=5.5cm]{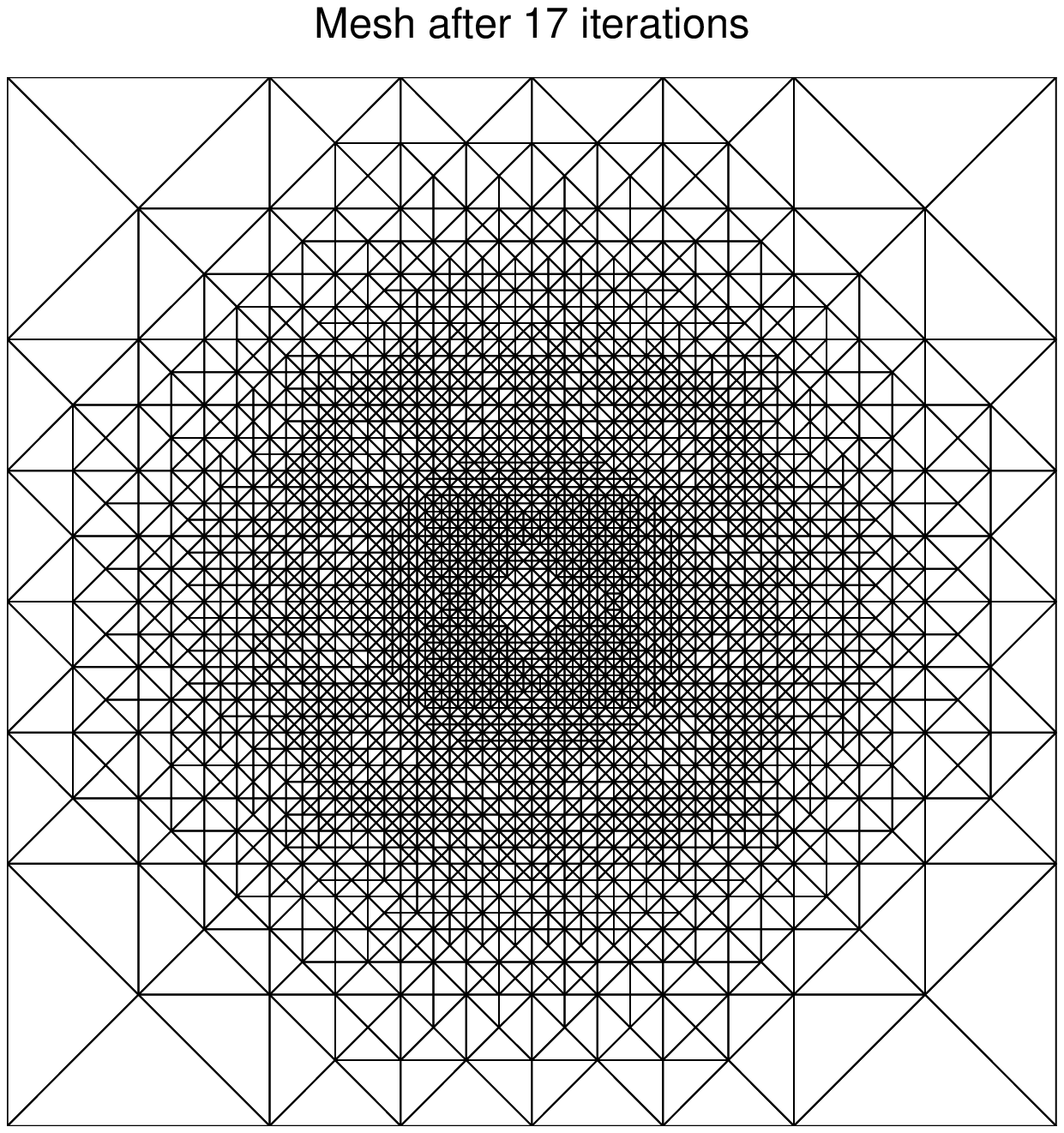}
\caption{The triangulations after adaptive iterations for
Example 1 by the linear element (left) and the quadratic element (right)}\label{Mesh_AFEM_Exam_1}
\end{figure}

\begin{figure}[ht]
\centering
\includegraphics[width=5.5cm,height=5.5cm]{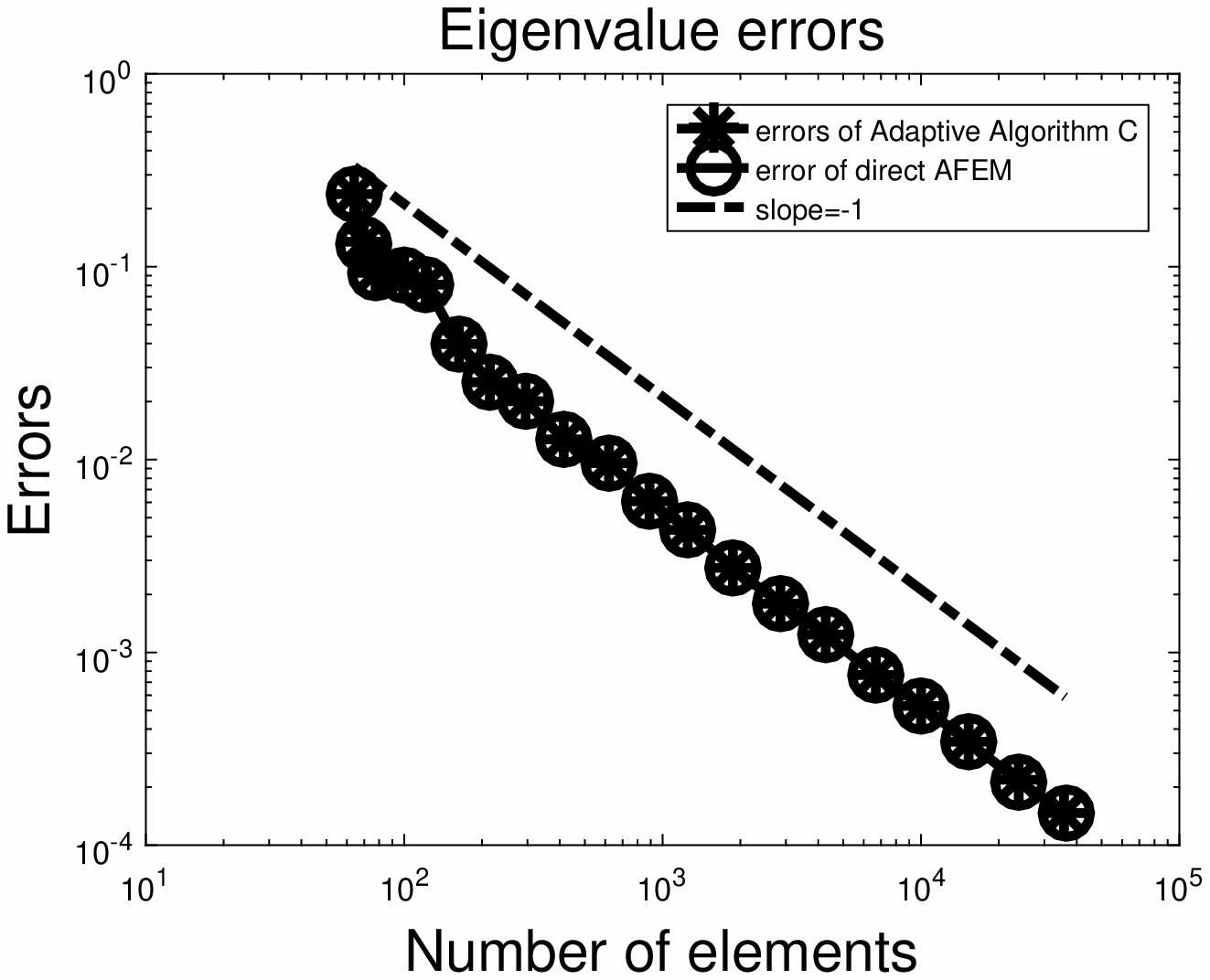}
\includegraphics[width=5.5cm,height=5.5cm]{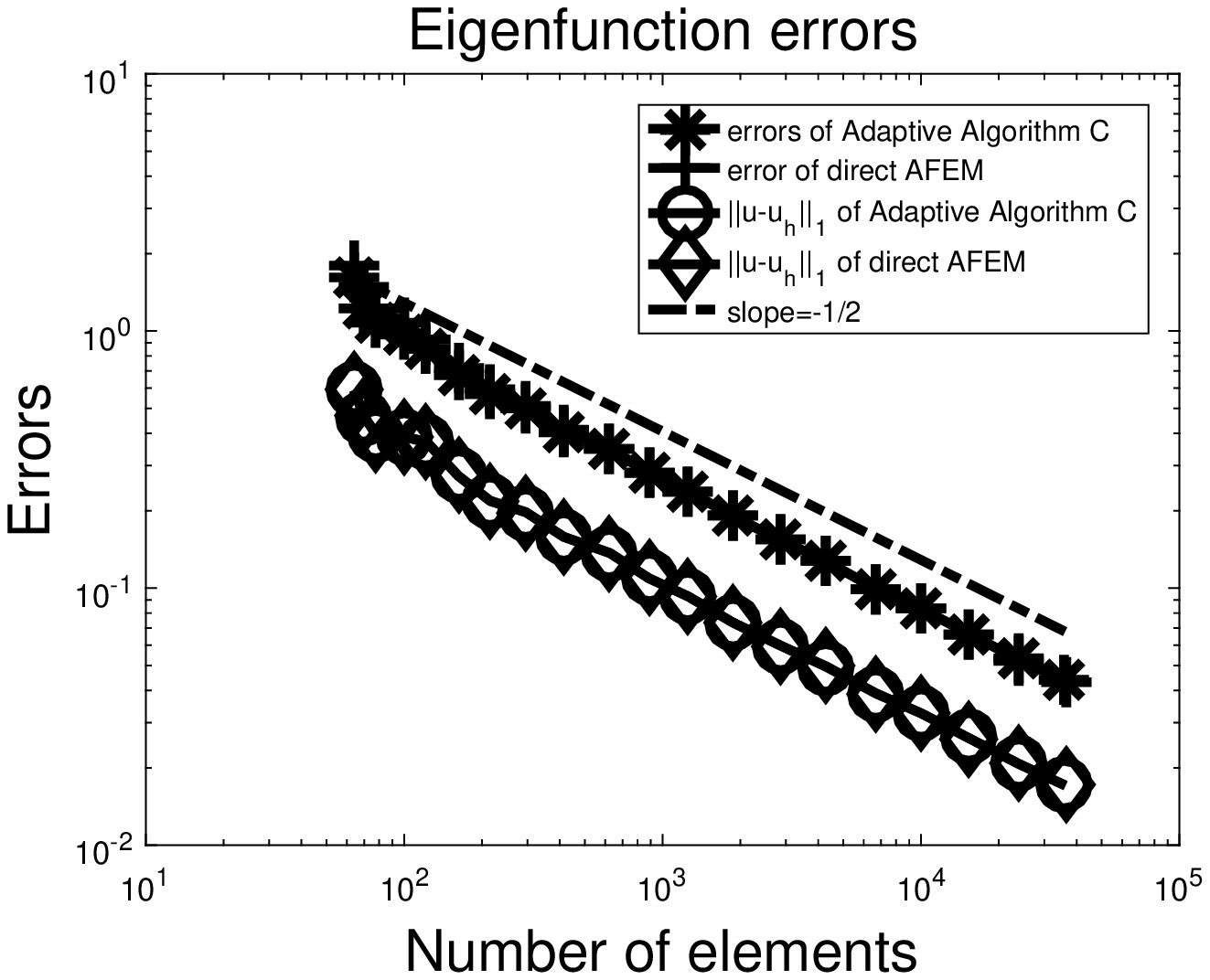}
\caption{The errors of the eigenvalue and the associated eigenfunction
approximations by {\bf Adaptive Algorithm $C$}
and direct AFEM for Example 1 with the linear element}\label{Convergence_AFEM_Exam_1_First}
\end{figure}

\begin{figure}[ht]
\centering
\includegraphics[width=5.5cm,height=5.5cm]{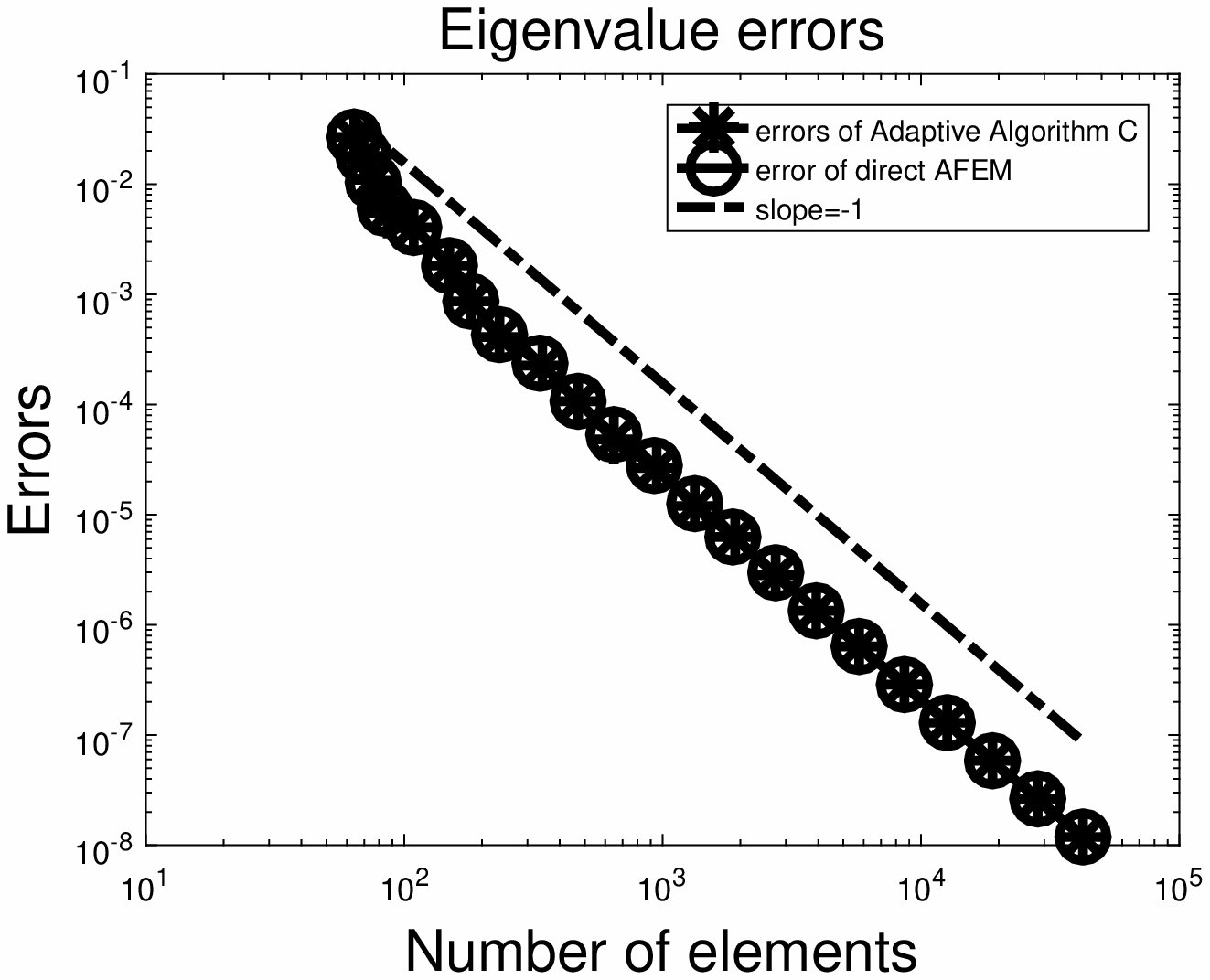}
\includegraphics[width=5.5cm,height=5.5cm]{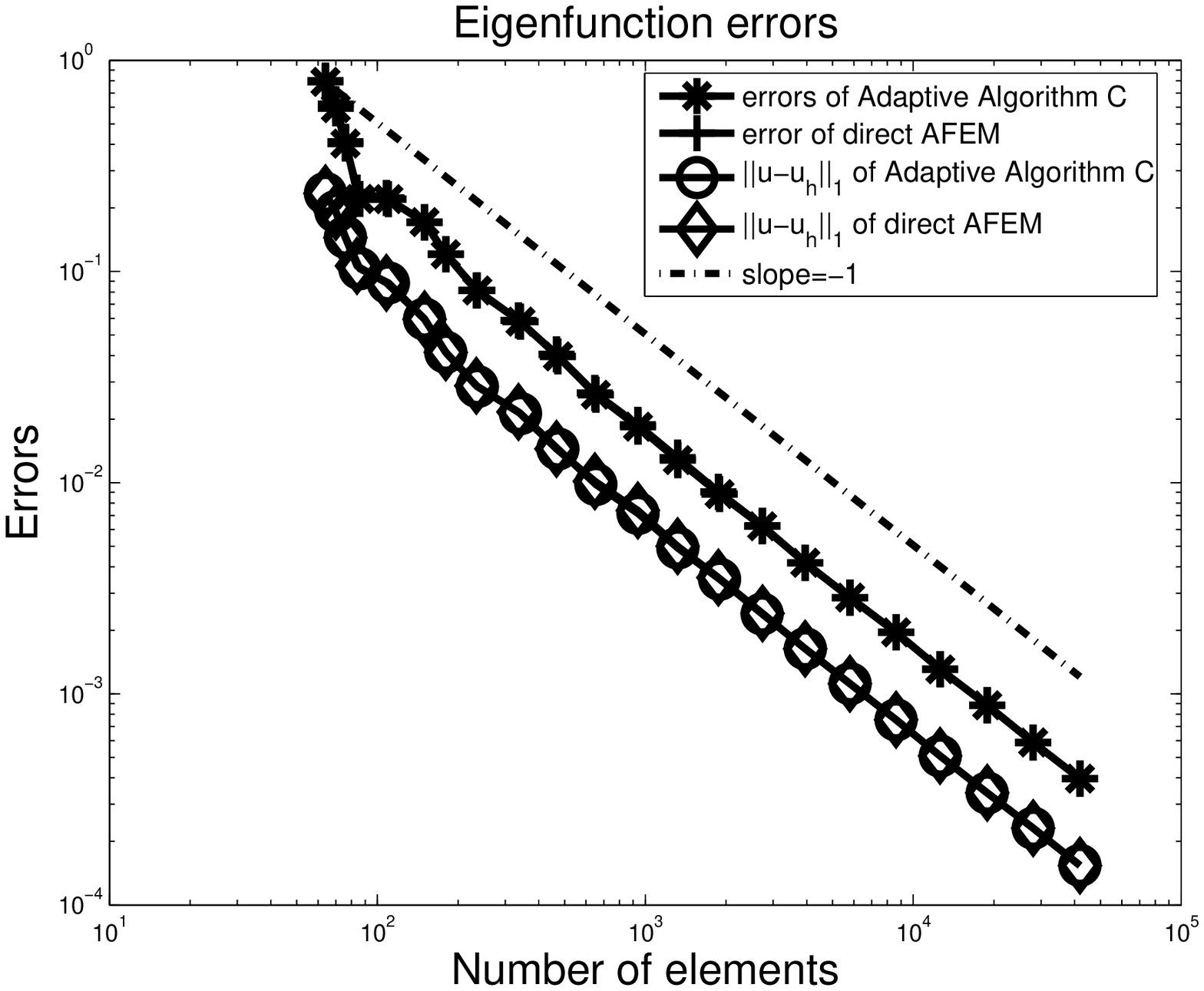}
\caption{The errors of the eigenvalue and the associated eigenfunction
approximations by {\bf Adaptive Algorithm $C$}
and direct AFEM for Example 1 with the quadratic element}\label{Convergence_AFEM_Exam_1_First_Quadratic}
\end{figure}

The eigenvalue problem is solved by {\bf Adaptive Algorithm $C$} with the parameters $\theta_1 = 0.4$
and $\theta_2=0.6$ for linear element and $\theta_1 = 0.4$ and $\theta_2=0.4$ for quadratic element.
Here, we check the numerical results for the first eigenvalue approximations.
Figures \ref{Mesh_AFEM_Exam_1} shows the triangulations by  {\bf Adaptive Algorithm $C$}
with the linear and quadratic finite element methods, respectively.
Figure \ref{Convergence_AFEM_Exam_1_First} gives the corresponding numerical results for the
first $20$ adaptive iterations with linear finite element method.
In order to show the efficiency of {\bf Adaptive Algorithm $C$} more clearly, we compare
the results with those obtained with direct AFEM. Similarly,
Figure \ref{Convergence_AFEM_Exam_1_First_Quadratic} shows the numerical results for
the first $22$ adaptive iterations with quadratic finite element method.

It is observed from Figures \ref{Convergence_AFEM_Exam_1_First} and \ref{Convergence_AFEM_Exam_1_First_Quadratic},
the approximations of eigenvalue as well as eigenfunction have the optimal convergence rate
which coincides with our theory. With {\bf Adaptive Algorithm $C$},
we only need to solve boundary value problems on each adaptively refined meshes and
small scale eigenvalue problems on the low dimensional space $V_H+{\rm span}\{u_k\}$
when $k=1, 5, 8, 11, 14, 17, 20$ ($\#\mathcal T_k=[64,122,296,884,2846,10044,36070]$ and $j_k=1,1,1,1,1,1,1$) and
$k=1, 4, 8, 11, 14, 17, 20$ ($\#\mathcal T_k=[64,86,236,654,1894,5810,18876]$ and $j_k=1,1,1,1,1,1,1$)
for linear element and quadratic element, respectively.
%on the meshes with the number of elements:
%$[64,122,296,884,2846,10044,36070]$ (${\rm when}\ k=1, 5, 8, 11, 14, 17, 20$, $\revise{{\rm and}\ j_k=1,1,1,1,1,1,1}$)
%for linear element and $[64,86,236,654,1894,5810,18876]$
%(${\rm when}\ k=1,4,8,11,14,17,20$, $\revise{{\rm and}\ j_k=1,1,1,1,1,1,1}$) for quadratic element, respectively.
But the accuracy obtained by {\bf Adaptive Algorithm $C$} is almost the same as the standard AFEM which
validate the efficiency of {\bf Adaptive Algorithm $C$}.

%%-----------------------------------------------------------------------------------
{\bf Example 2.}\
In the second example, we consider the Laplace eigenvalue problem on the $L$-shape domain
\begin{equation}\label{eigenproblem_Exam_2}
\left\{
\begin{array}{rcl}
-\Delta u &=&\lambda u\ \ \ \ {\rm in}\ \Omega,\\
u&=&0\ \ \ \ \ \ {\rm on}\ \partial\Omega,\\
\|u\|_{a,\Omega}&=&1,
\end{array}
\right.
\end{equation}
where $\Omega=(-1,1)\times(-1,1)\backslash[0, 1)\times (-1, 0]$.
Since $\Omega$ has a reentrant corner, eigenfunctions with singularities are expected. The
convergence order for eigenvalue approximations is less than $2$ by the linear finite element method
 which is the order predicted by the theory for regular eigenfunctions.

\begin{figure}[ht]
\centering
\includegraphics[width=6cm,height=6cm]{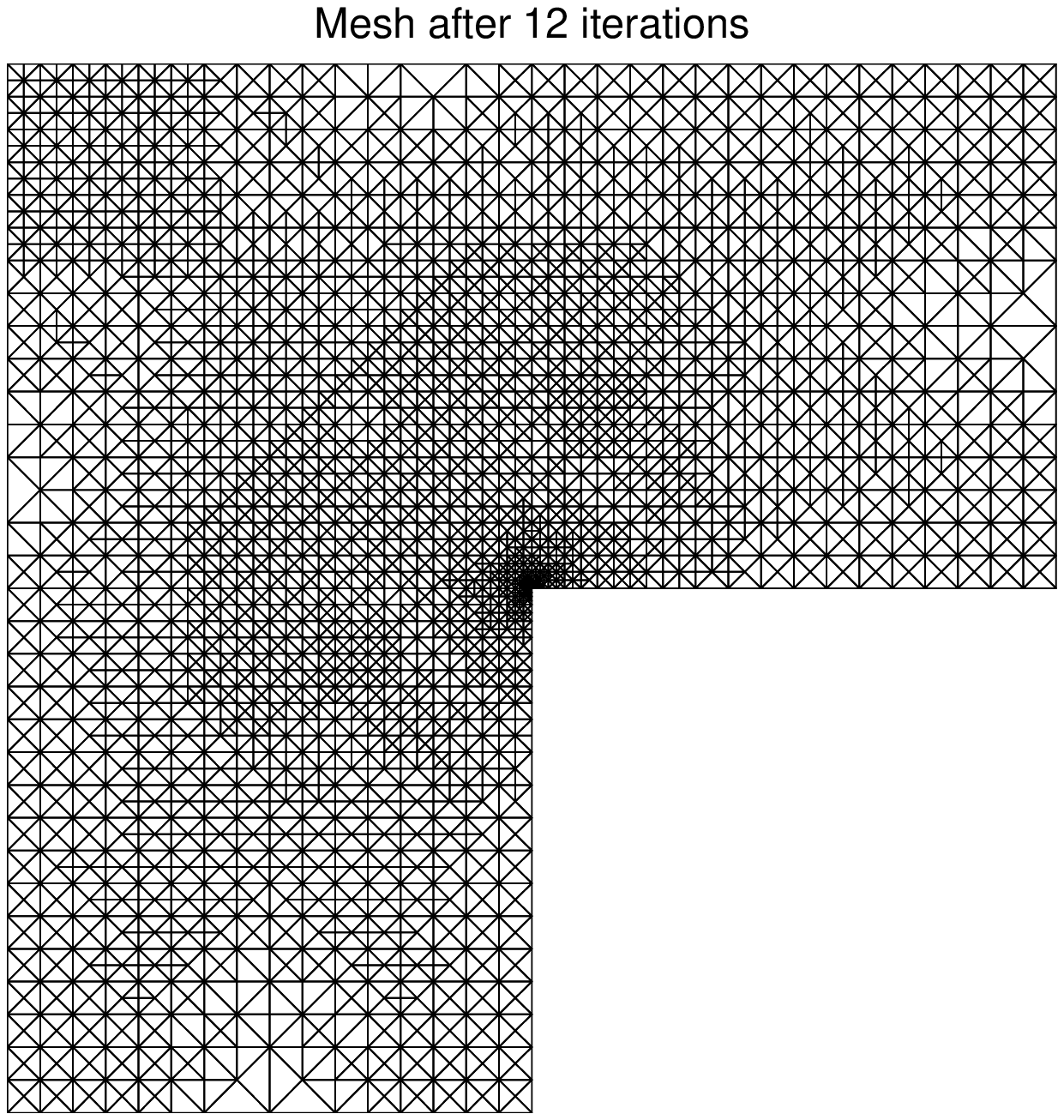}
\includegraphics[width=6cm,height=6cm]{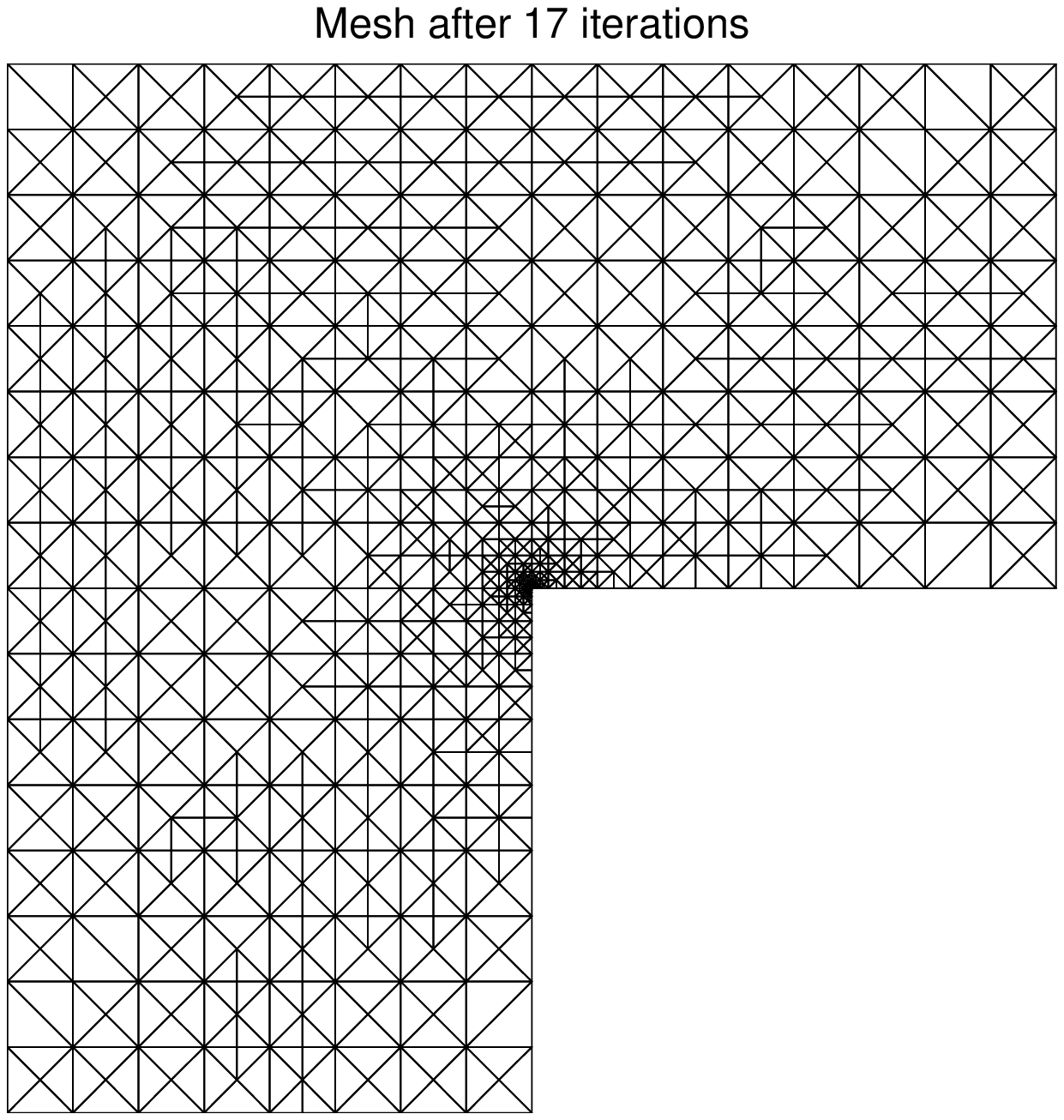}
\caption{The triangulations after adaptive iterations for
Example 2 by the linear element (left) and the quadratic element (right)}\label{Mesh_AFEM_Exam_2}
\end{figure}

\begin{figure}[ht]
\centering
\includegraphics[width=5.5cm,height=5.5cm]{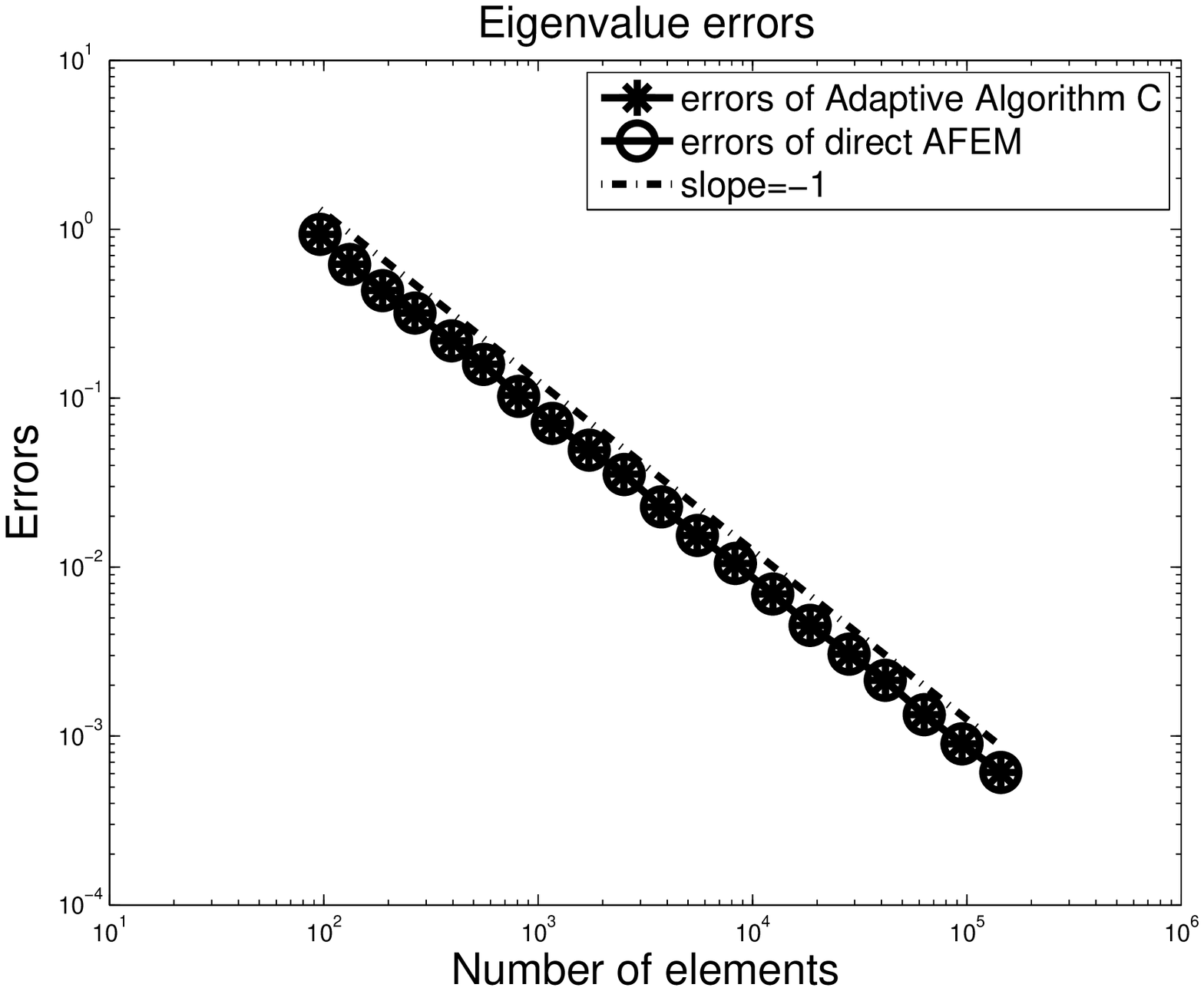}
\includegraphics[width=5.5cm,height=5.5cm]{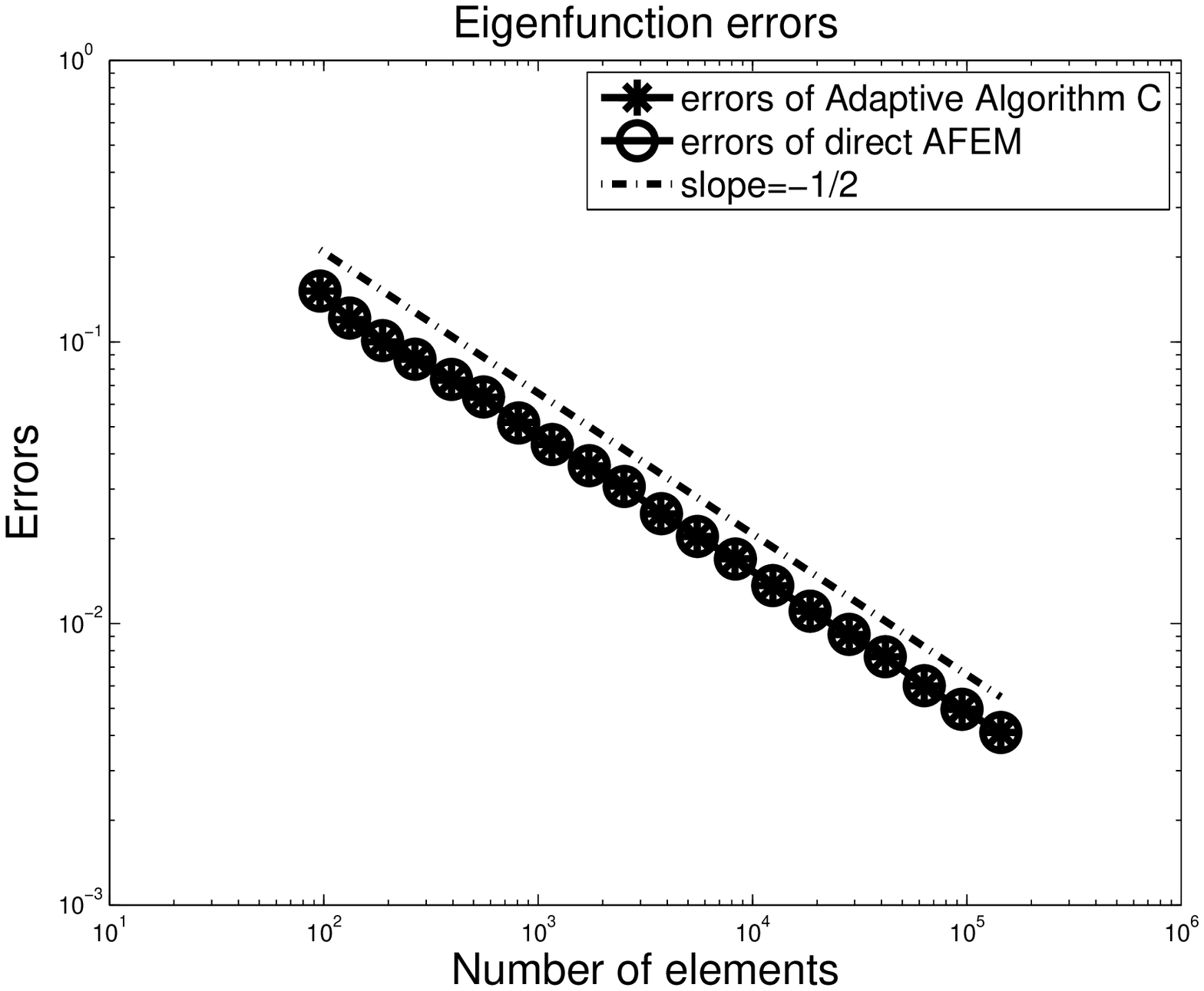}
\caption{The errors of the smallest eigenvalue approximations and the a
posteriori errors of the associated eigenfunction approximations by
{\bf Adaptive Algorithm $C$}
and direct AFEM for Example 2 with the linear element}\label{Convergence_AFEM_Exam_2_First}
\end{figure}

\begin{figure}[ht]
\centering
\includegraphics[width=5.5cm,height=5.5cm]{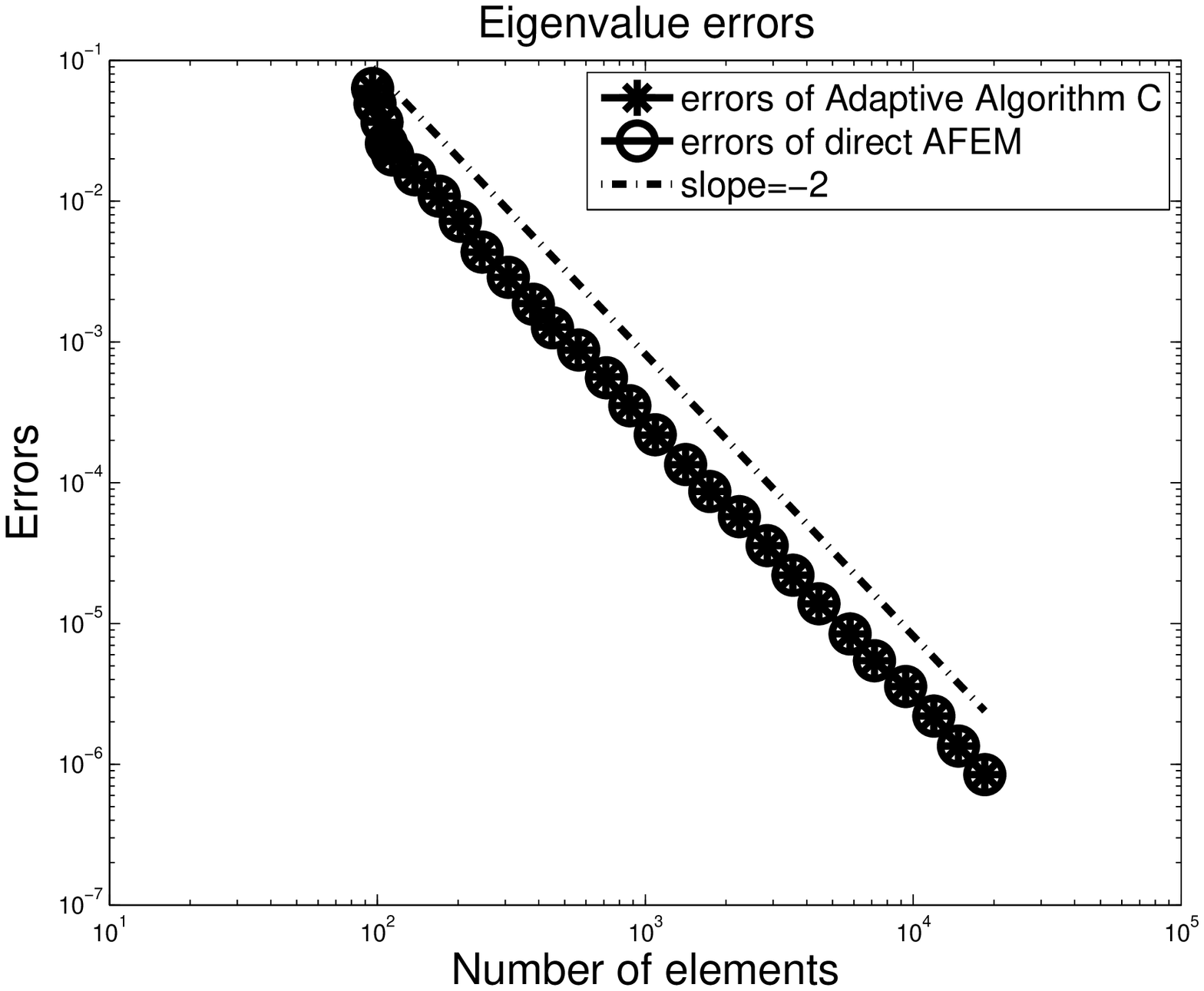}
\includegraphics[width=5.5cm,height=5.5cm]{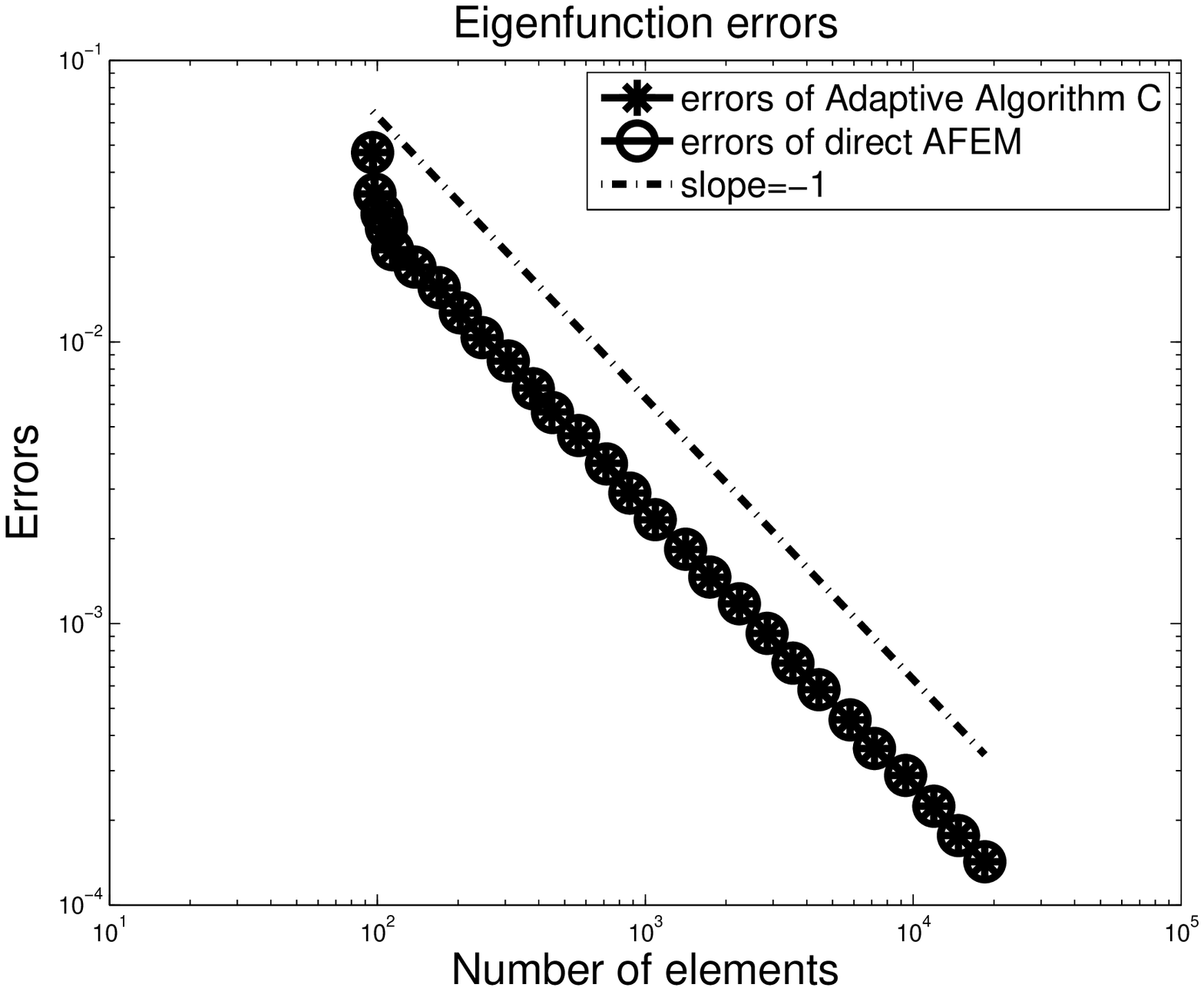}
\caption{The errors of the smallest eigenvalue approximations and the a
posteriori errors of the associated eigenfunction approximations by
{\bf Adaptive Algorithm $C$}
and direct AFEM for Example 2 with the quadratic element}\label{Convergence_AFEM_Exam_2_First_Quadratic}
\end{figure}

Here, we give the numerical results of {\bf Adaptive Algorithm $C$} with parameters $\theta_1 = 0.4$
and $\theta_2=0.6$ for linear element and $\theta_1=0.4$ and $\theta_2=0.4$ for quadratic element, respectively.
First, we investigate the numerical results for the first eigenvalue approximations.
Since the exact eigenvalue is not known, an adequately accurate approximation
$\lambda = 9.6397238440219$ is chosen as the exact first eigenvalue for our numerical tests.
Figure \ref{Mesh_AFEM_Exam_2} shows the triangulations after
adaptive iterations with the linear and quadratic finite element methods, respectively.
Figures \ref{Convergence_AFEM_Exam_2_First} and \ref{Convergence_AFEM_Exam_2_First_Quadratic}
give the corresponding numerical results. In order to show the efficiency of {\bf Adaptive Algorithm $C$}
more clearly, we  also compare the results with those obtained by direct AFEM.
\revise{With {\bf Adaptive Algorithm $C$}, it is only required to solve boundary value problems on the
adaptively refined triangulations and small scale eigenvalue problems on the low dimensional space $V_H+{\rm span}\{u_k\}$
when the numbers of elements of the meshes are $[96,266,808,2508,8292 ,28276,94922]$
($k=1, 4, 7, 10, 13, 16, 19$ and $j_k=1,1,1,1,1,1,1$) for linear element and
$[96,138,380,1090 ,2852,7164 ,18484]$ ($k=1,6,11,16,20,24,28$ and $j_k=1,1,1,1,1,1,1$)
for quadratic element, respectively. }
%%-----------------------------------------------------------------------------------
%%%-----------------------------------------------------------------------------------

\begin{figure}[ht]
\centering
\includegraphics[width=5.5cm,height=5.5cm]{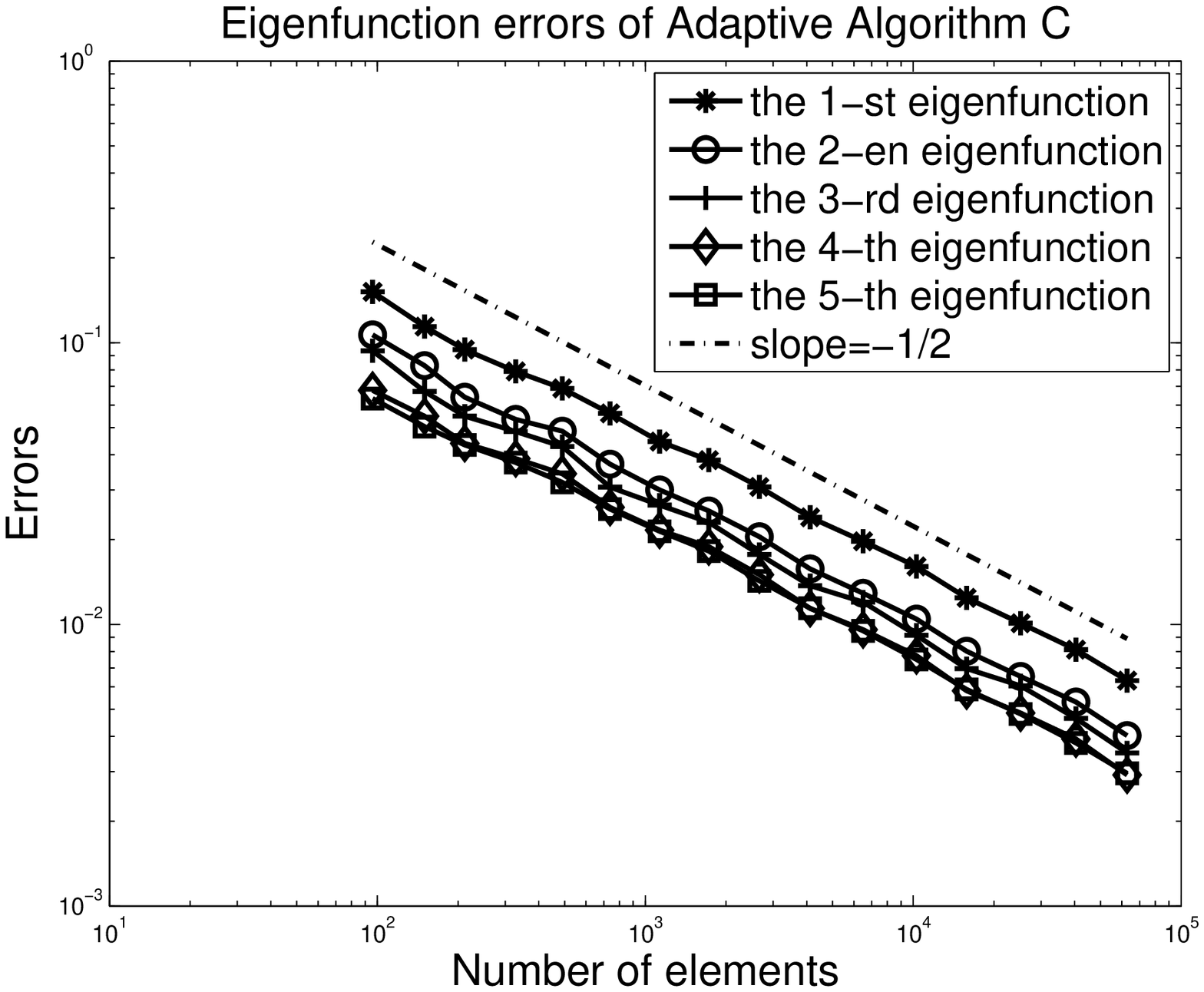}
\includegraphics[width=5.5cm,height=5.5cm]{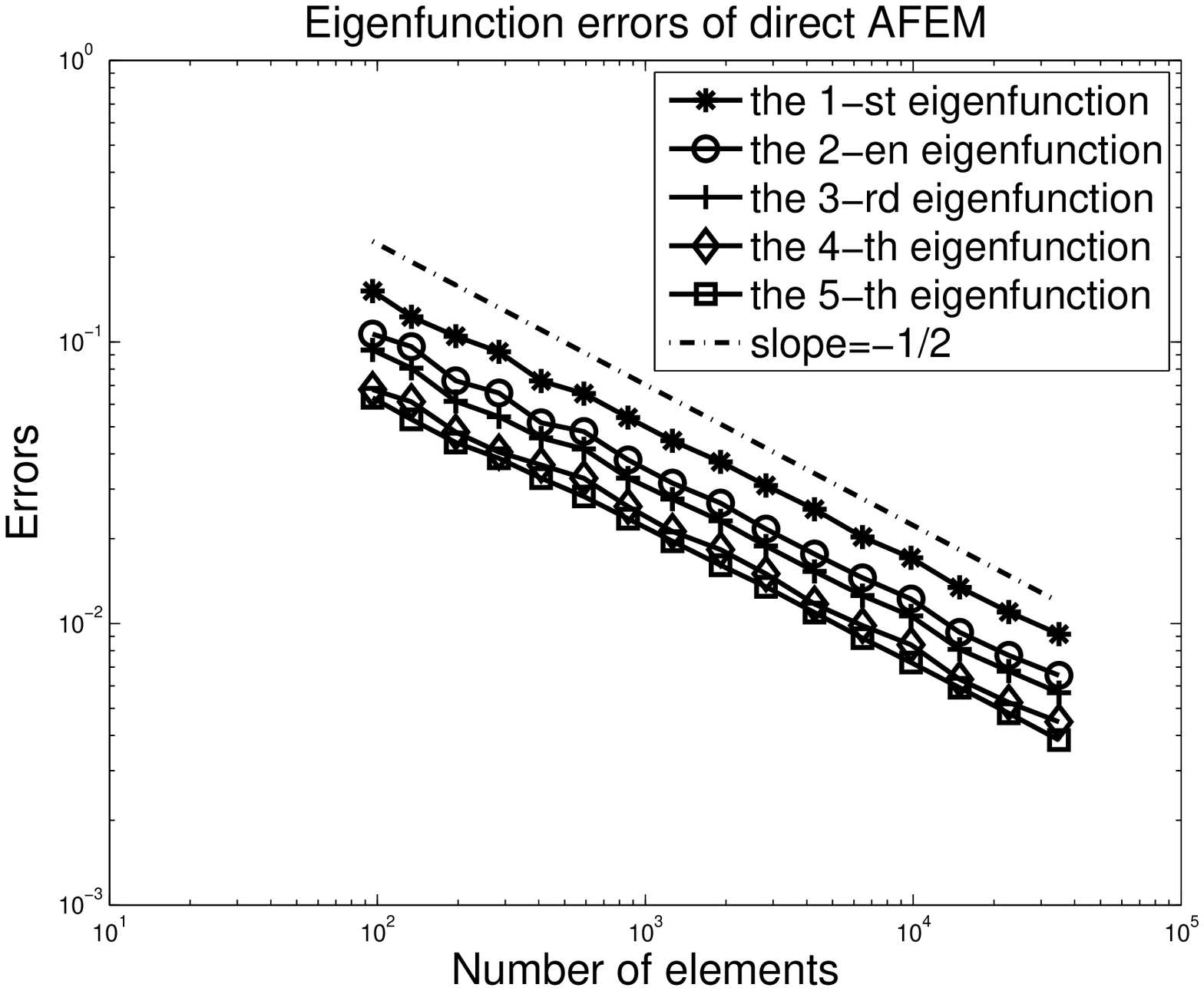}
\caption{The a posteriori error estimates of the eigenfunction
 approximations by {\bf Adaptive Algorithm $C$}
and direct AFEM for Example 2 with the linear element}\label{Convergence_AFEM_Exam_2_5_Small}
\end{figure}

\begin{figure}[ht]
\centering
\includegraphics[width=5.5cm,height=5.5cm]{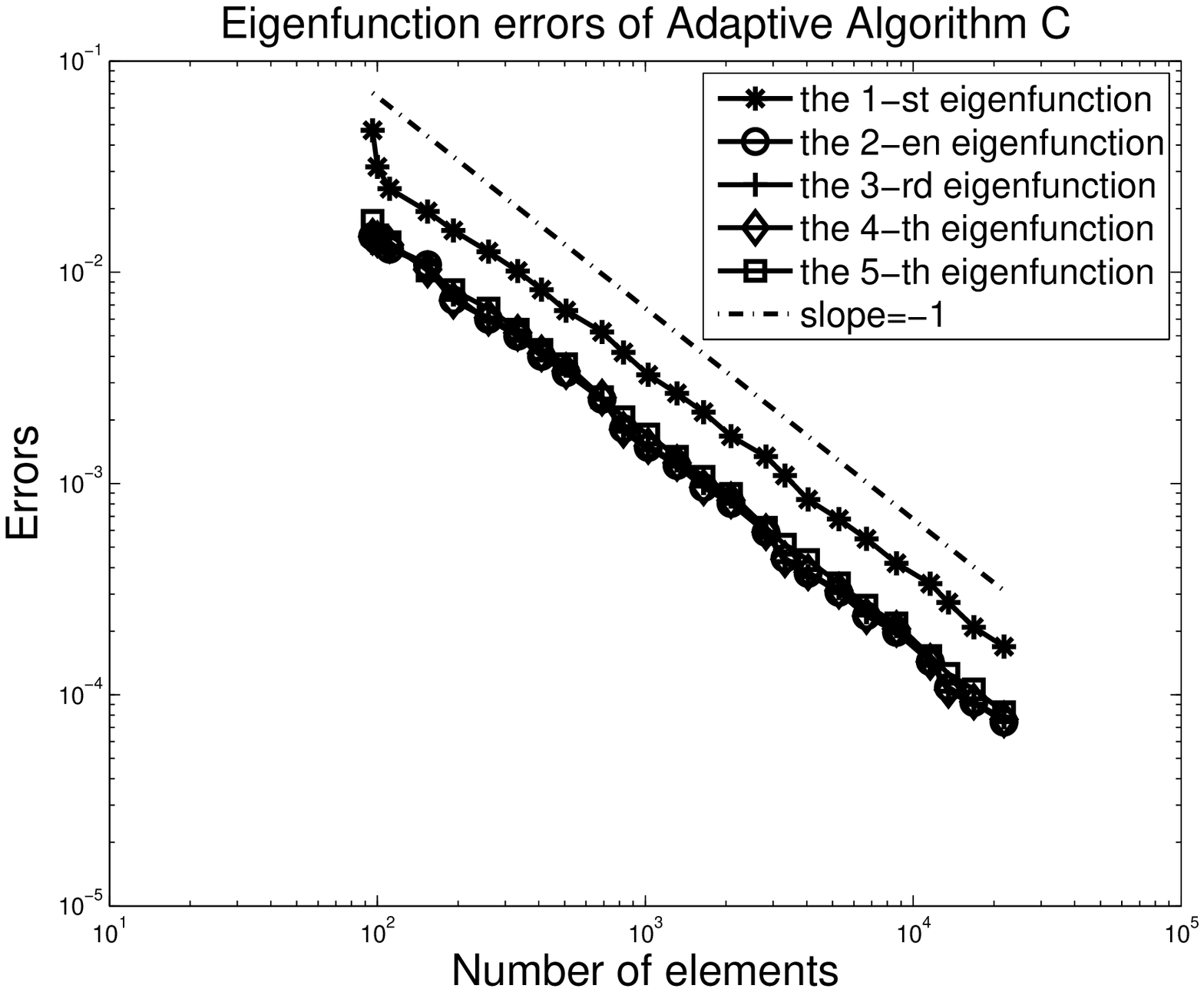}
\includegraphics[width=5.5cm,height=5.5cm]{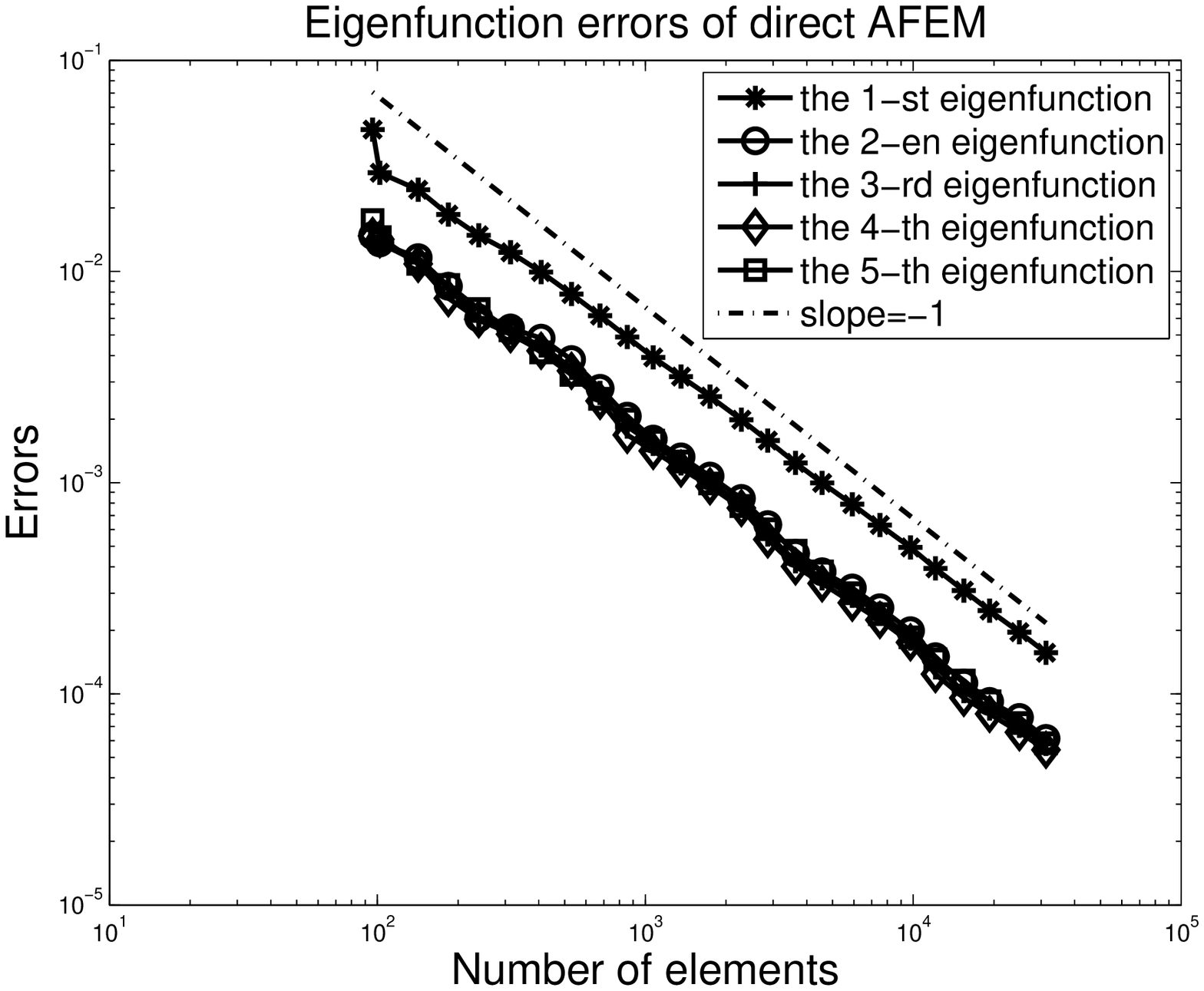}
\caption{The a posteriori error estimates of the eigenfunction
approximations by {\bf Adaptive Algorithm $C$}
and direct AFEM for Example 2 with the quadratic element}\label{Convergence_AFEM_Exam_2_5_Small_Quadratic}
\end{figure}

We also test {\bf Adaptive Algorithm $C$} for the first $5$  eigenvalue approximations and their associated
eigenfunction approximations.  Figures \ref{Convergence_AFEM_Exam_2_5_Small} and
\ref{Convergence_AFEM_Exam_2_5_Small_Quadratic} show the corresponding a posteriori
error estimator produced by {\bf Adaptive Algorithm $C$}
and direct AFEM with the linear and quadratic finite element methods, respectively.
In these cases, {\bf Adaptive Algorithm $C$} only need to solve the small scale eigenvalue problems
on the low dimensional space $V_H+{\rm span}\{u_k\}$ when the numbers of elements of the meshes are
$[96,328,1130,4124,15832,62802]$
($k=1, 4, 7, 10, 13, 16$ and $j_k=1,1,1,1,1,1$) for linear element and
$[ 96,192,690,2086,6686,16832]$ ($k=1, 5, 10, 15, 20, 24$ and $j_k=1,1,1,1,1,1$) for quadratic element, respectively.
%%%-----------------------------------------------------------------------------------

\revise{
From Figures \ref{Convergence_AFEM_Exam_2_First}, \ref{Convergence_AFEM_Exam_2_First_Quadratic},
\ref{Convergence_AFEM_Exam_2_5_Small} and \ref{Convergence_AFEM_Exam_2_5_Small_Quadratic},
we can find the approximations of eigenvalues as well as eigenfunctions have the optimal convergence
rate as the direct AFEM which coincides with our theory.}

%%-----------------------------------------------------------------------------------
{\bf Example 3.}\ In this example, we consider the following second order elliptic eigenvalue problem
\begin{equation}
\left\{
\begin{array}{rcl}
-\nabla\cdot(\mathcal{A}\nabla u) +\varphi u&=&\lambda u\ \ \ {\rm in}\ \Omega,\\
u&=&0\ \ \ \ \ {\rm on}\ \partial\Omega,\\
\|u\|_{a,\Omega}&=&1,
\end{array}
\right.
\end{equation}
with
\begin{equation*}
\mathcal{A} =
\left(
\begin{array}{cc}
1+(x_1-\frac{1}{2})^2 & (x_1-\frac{1}{2})(x_2-\frac{1}{2})\\
(x_1-\frac{1}{2})(x_2-\frac{1}{2}) & 1+(x_2-\frac{1}{2})^2
\end{array}
\right),
\end{equation*}
$\varphi=e^{(x_1-\frac{1}{2})(x_2-\frac{1}{2})}$ and
$\Omega=(-1,1)\times(-1,1)\backslash[0, 1)\times (-1, 0]$.

\begin{figure}[ht]
\centering
\includegraphics[width=6cm,height=6cm]{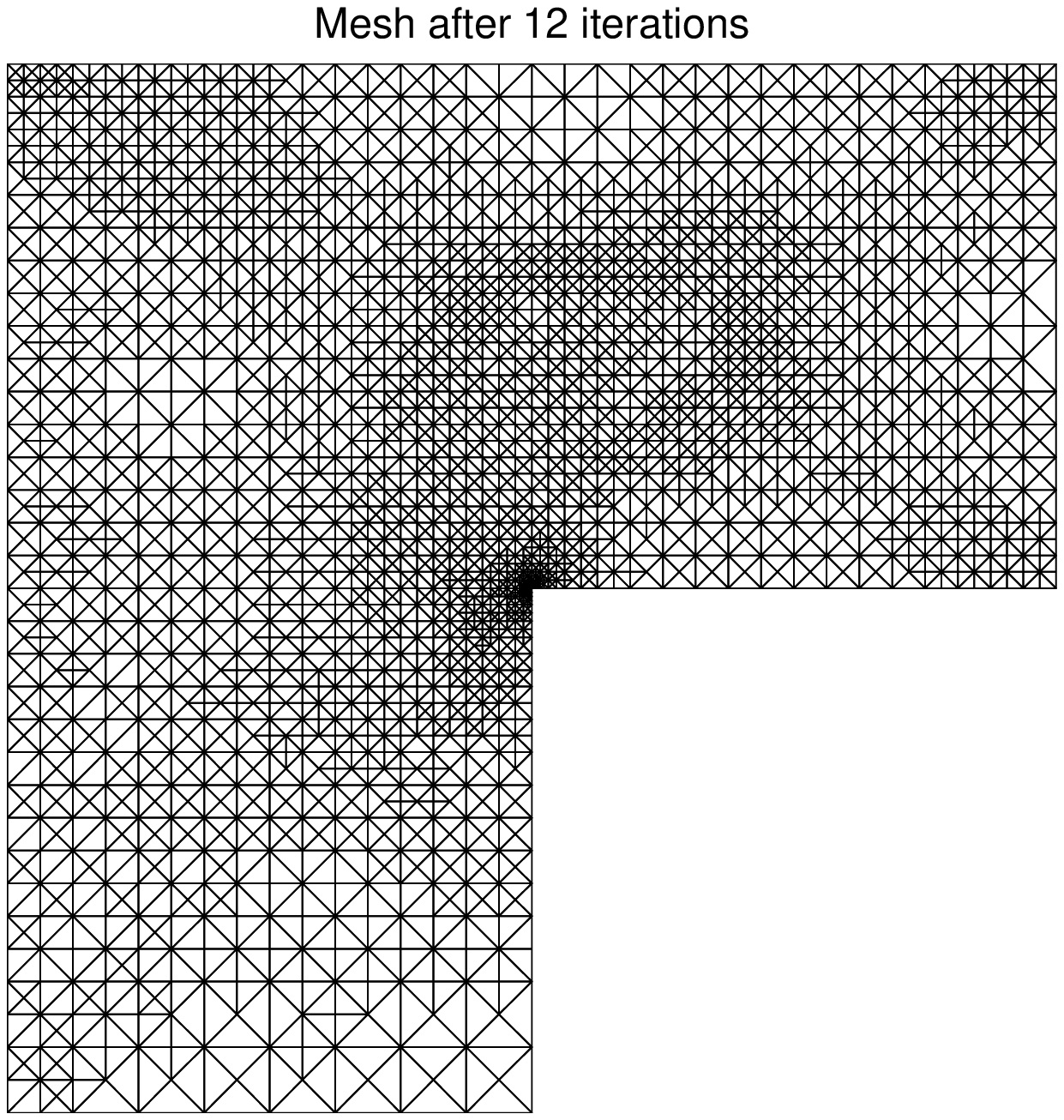}
\includegraphics[width=6cm,height=6cm]{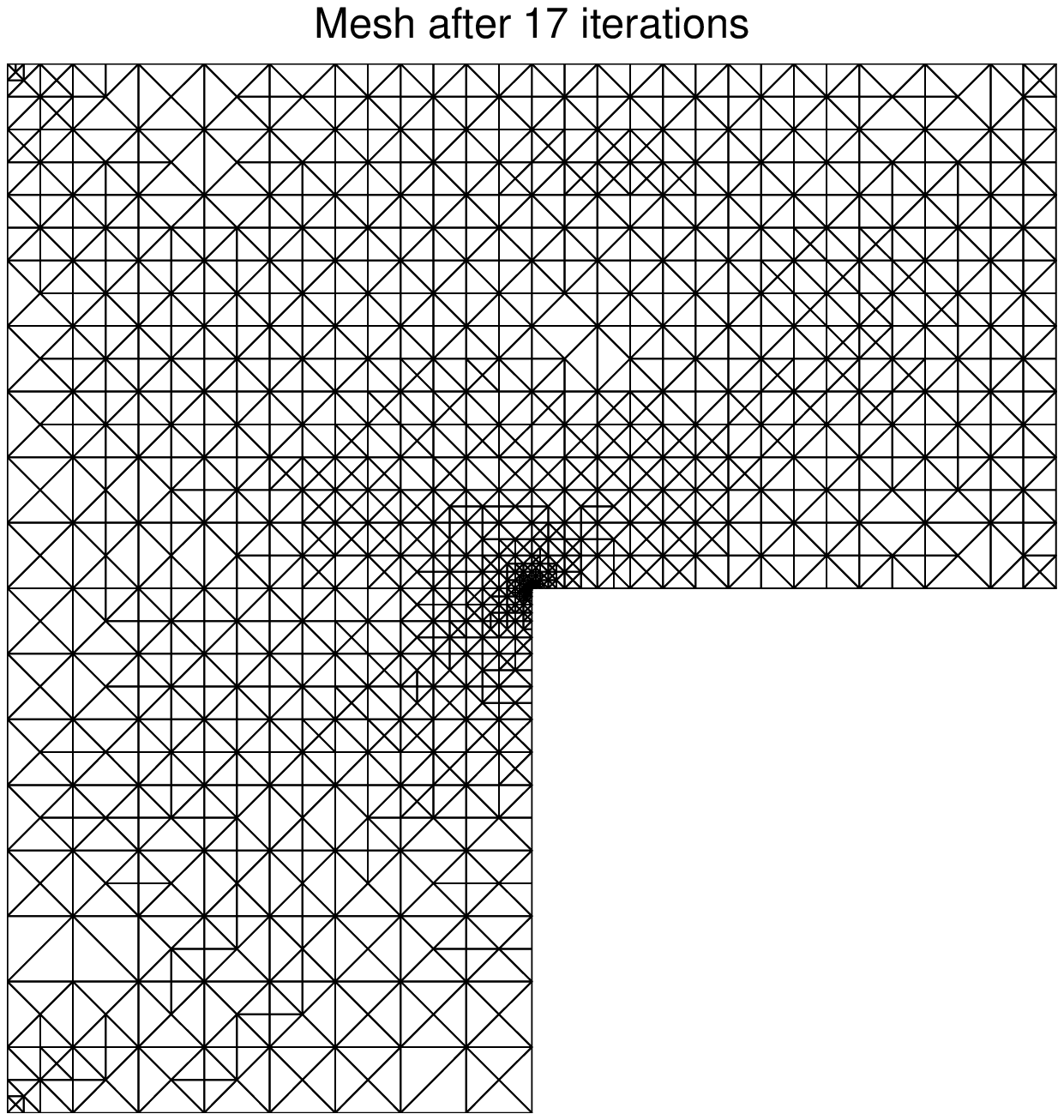}
\caption{The triangulations after adaptive iterations for
Example 3 by the linear element (left) and the quadratic element (right)} \label{Mesh_AFEM_Exam_3}
\end{figure}

\begin{figure}[ht]
\centering
\includegraphics[width=5.5cm,height=5.5cm]{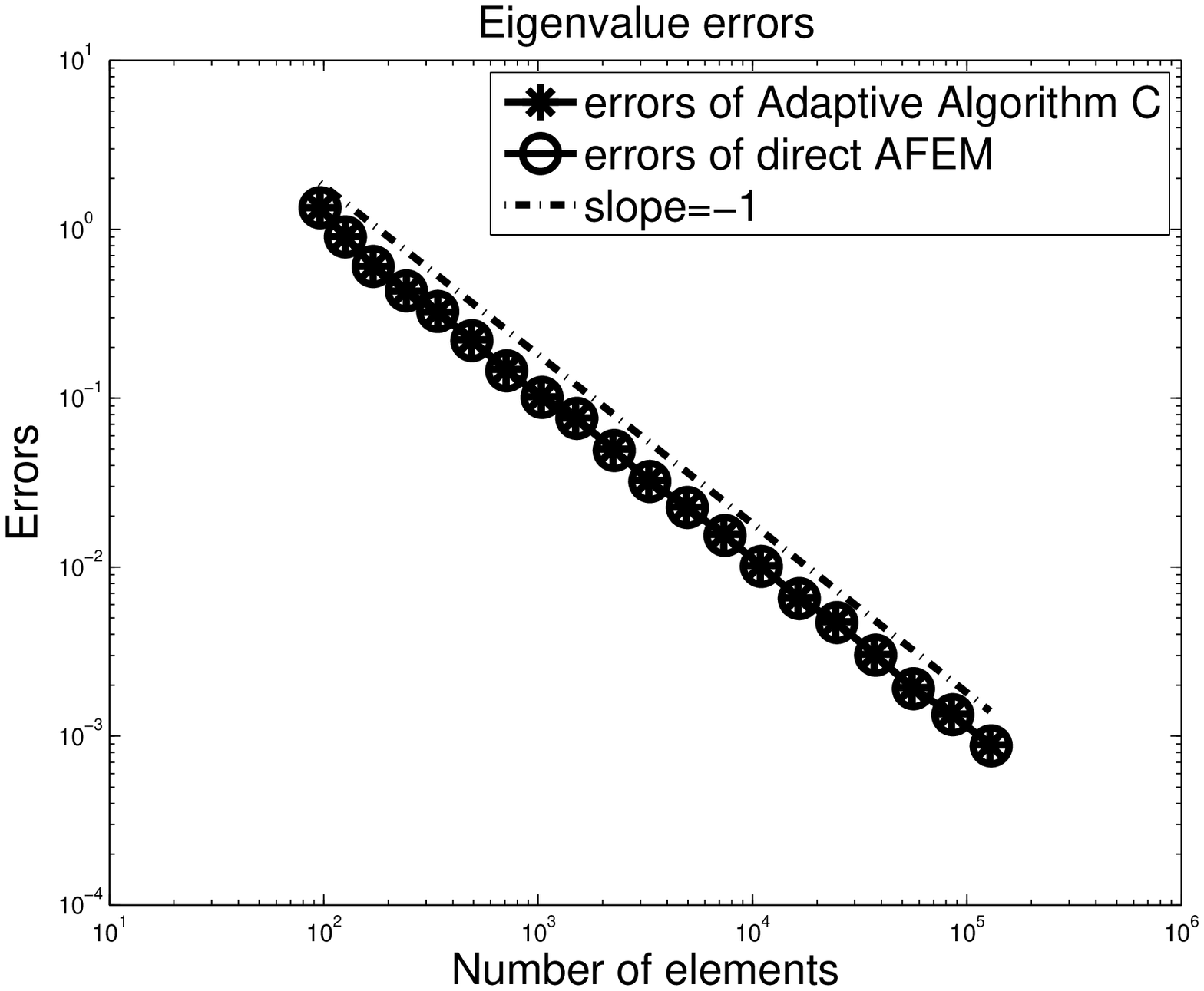}
\includegraphics[width=5.5cm,height=5.5cm]{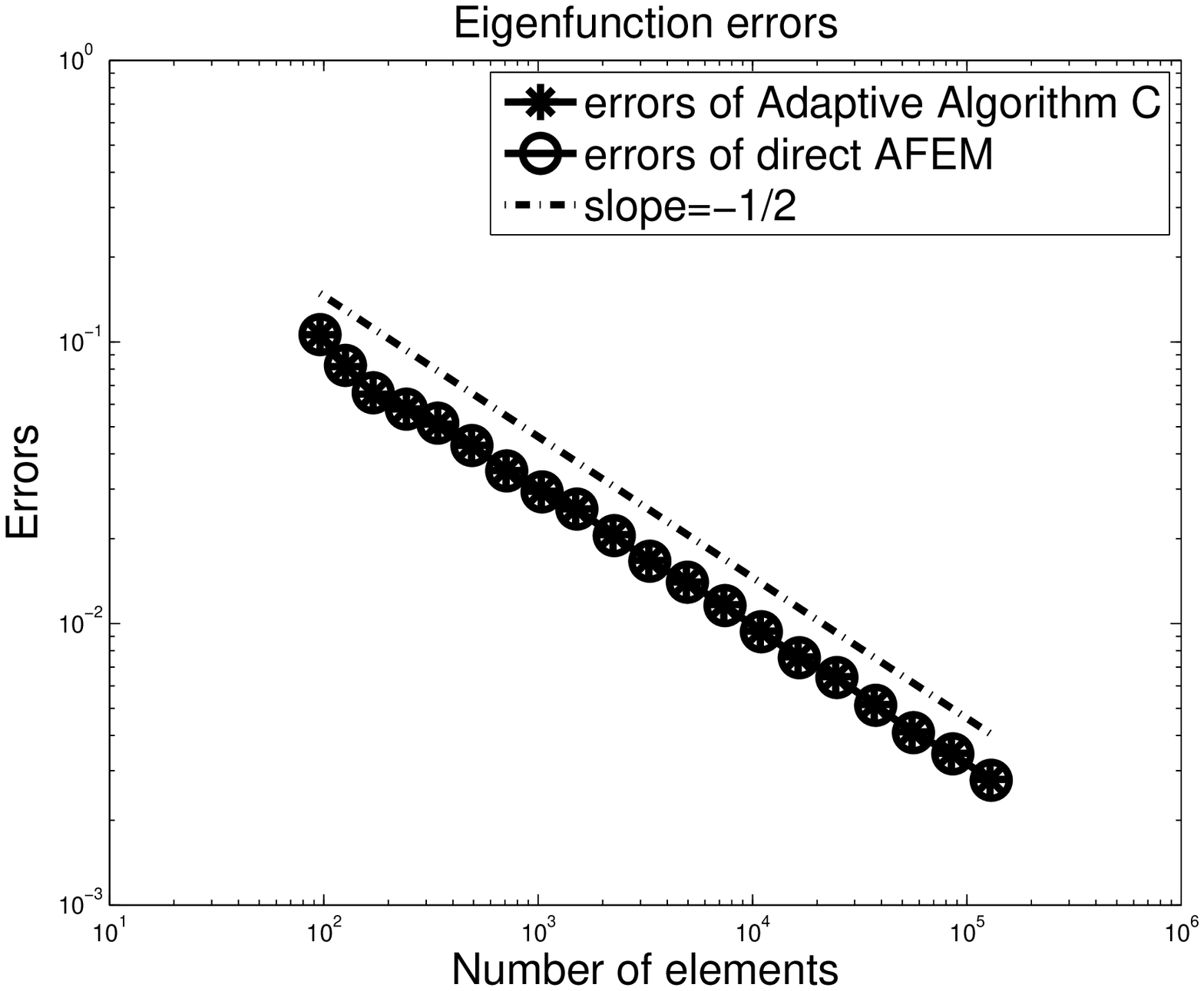}
\caption{The errors of the smallest eigenvalue and the associated
eigenfunction approximations by {\bf Adaptive Algorithm $C$}
and direct AFEM for Example 3 with the linear element} \label{Convergence_AFEM_Exam_3_First}
\end{figure}

\begin{figure}[ht]
\centering
\includegraphics[width=5.5cm,height=5.5cm]{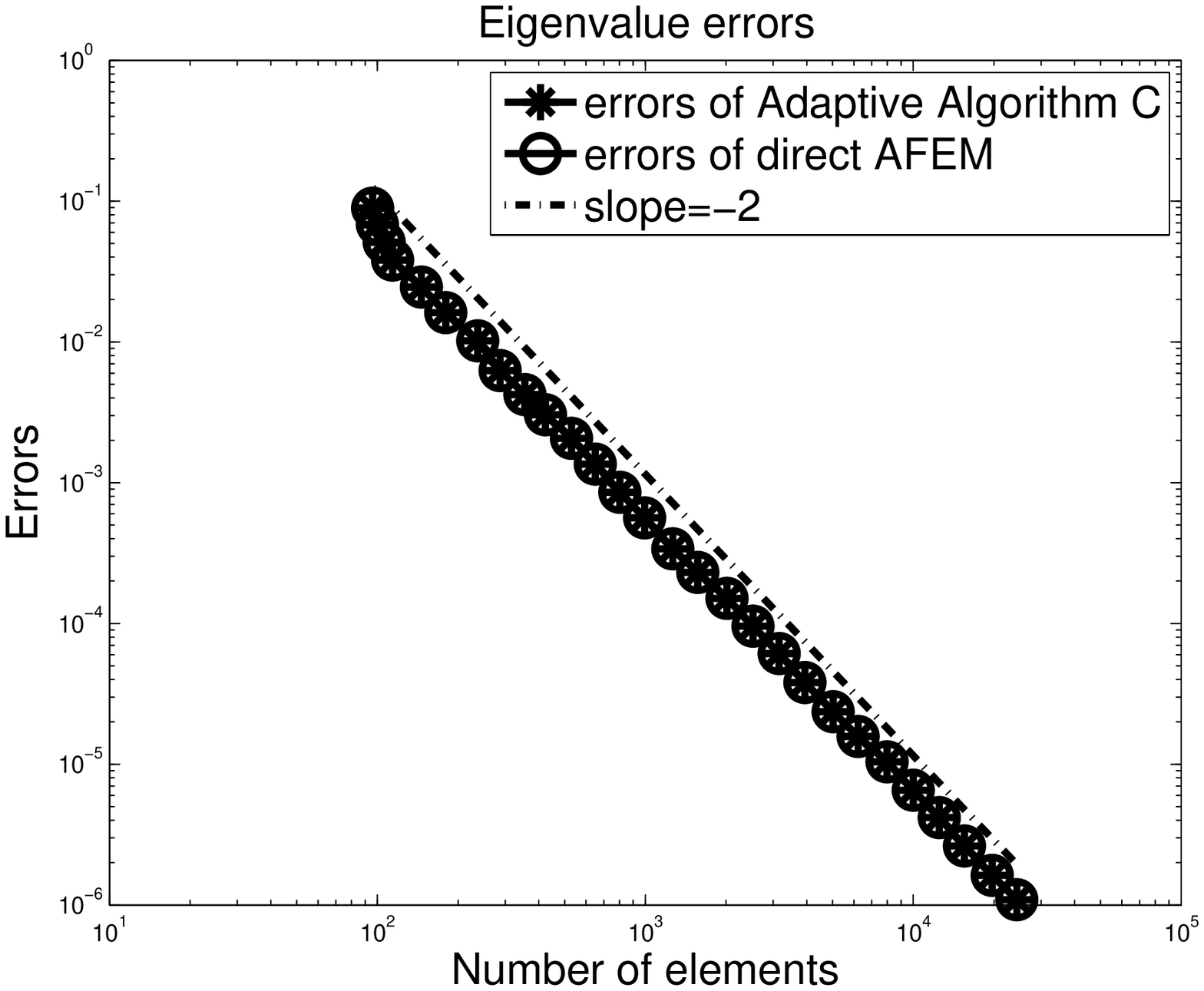}
\includegraphics[width=5.5cm,height=5.5cm]{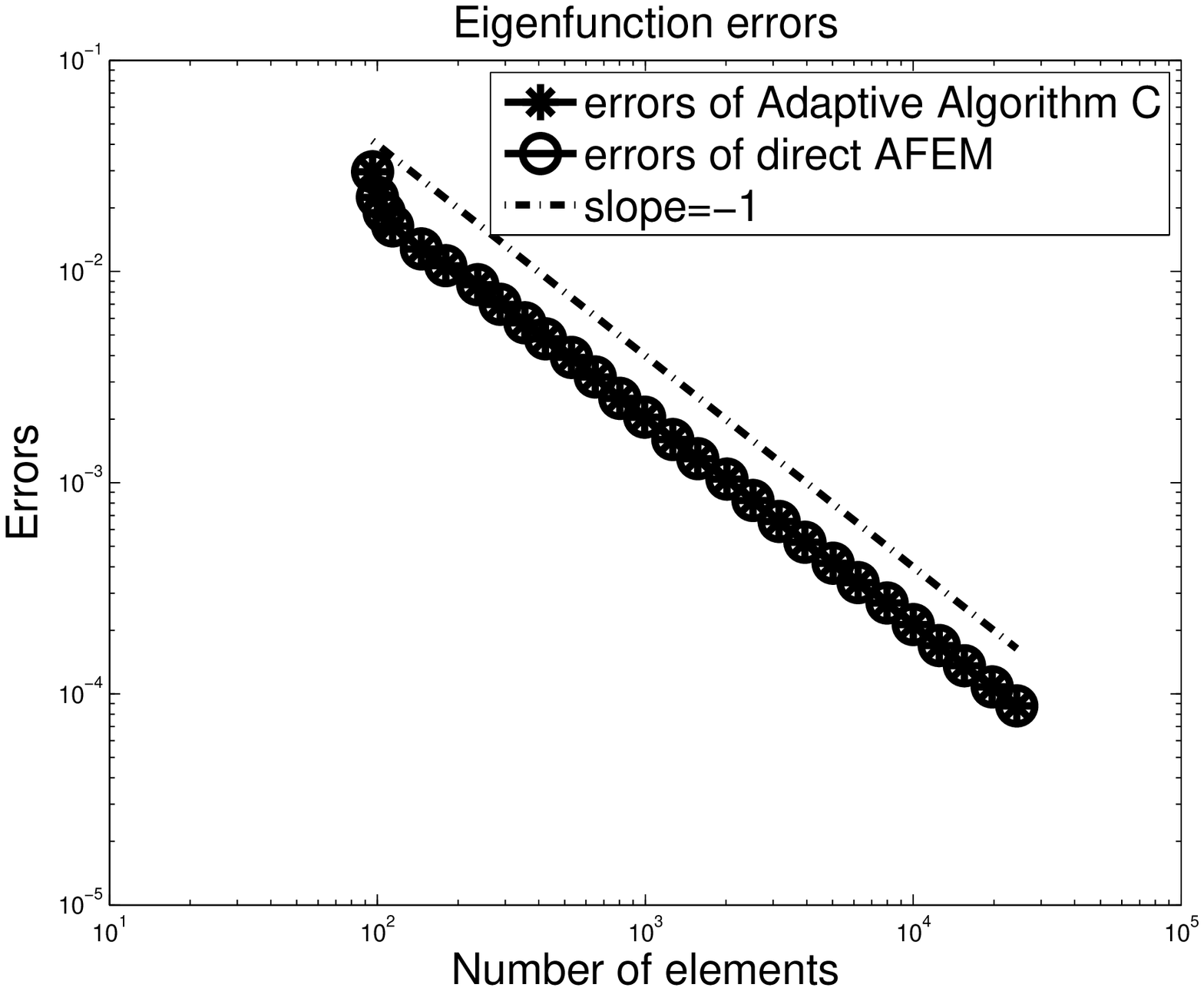}
\caption{The errors of the smallest eigenvalue and the associated
eigenfunction approximations by {\bf Adaptive Algorithm $C$}
and direct AFEM for Example 3 with the quadratic element}\label{Convergence_AFEM_Exam_3_First_Quadratic}
\end{figure}

We give the numerical results  by {\bf Adaptive Algorithm $C$} with parameters $\theta_1 = 0.4$ and $\theta_2=0.6$ for linear element
and $\theta_1=0.4$ and $\theta_2=0.4$ for quadratic element, respectively.
We first investigate the numerical results for the first eigenvalue approximations.
Since the exact eigenvalue is not known neither, we choose an adequately accurate
approximation $\lambda= 15.134144021256400$ as the exact eigenvalue for our numerical tests.
Figure \ref{Mesh_AFEM_Exam_3} shows the triangulations after adaptive iterations by the
linear and quadratic finite element methods, respectively.
Figures \ref{Convergence_AFEM_Exam_3_First} and \ref{Convergence_AFEM_Exam_3_First_Quadratic}
give the corresponding numerical results by the linear and quadratic finite element methods, respectively.
Similarly, we also compare the results with those obtained with direct AFEM.
\revise{
It is only required to solve the small scale eigenvalue problem on the low dimensional space $V_H+{\rm span}\{u_k\}$
when the numbers of elements of the meshes are $[96,243,1034,3282,10870,37030,128259]$
($k=1, 4, 8, 11, 14, 17, 20$ and $j_k=1,1,1,1,1,1,1$) for linear element and
$[96, 180, 530, 1571, 5026, 15534]$ ($k=1,6,11,16,21,26$ and $j_k=1,1,1,1,1,1$)
for quadratic element, respectively.}

\begin{figure}[ht]
\centering
\includegraphics[width=5.5cm,height=5.5cm]{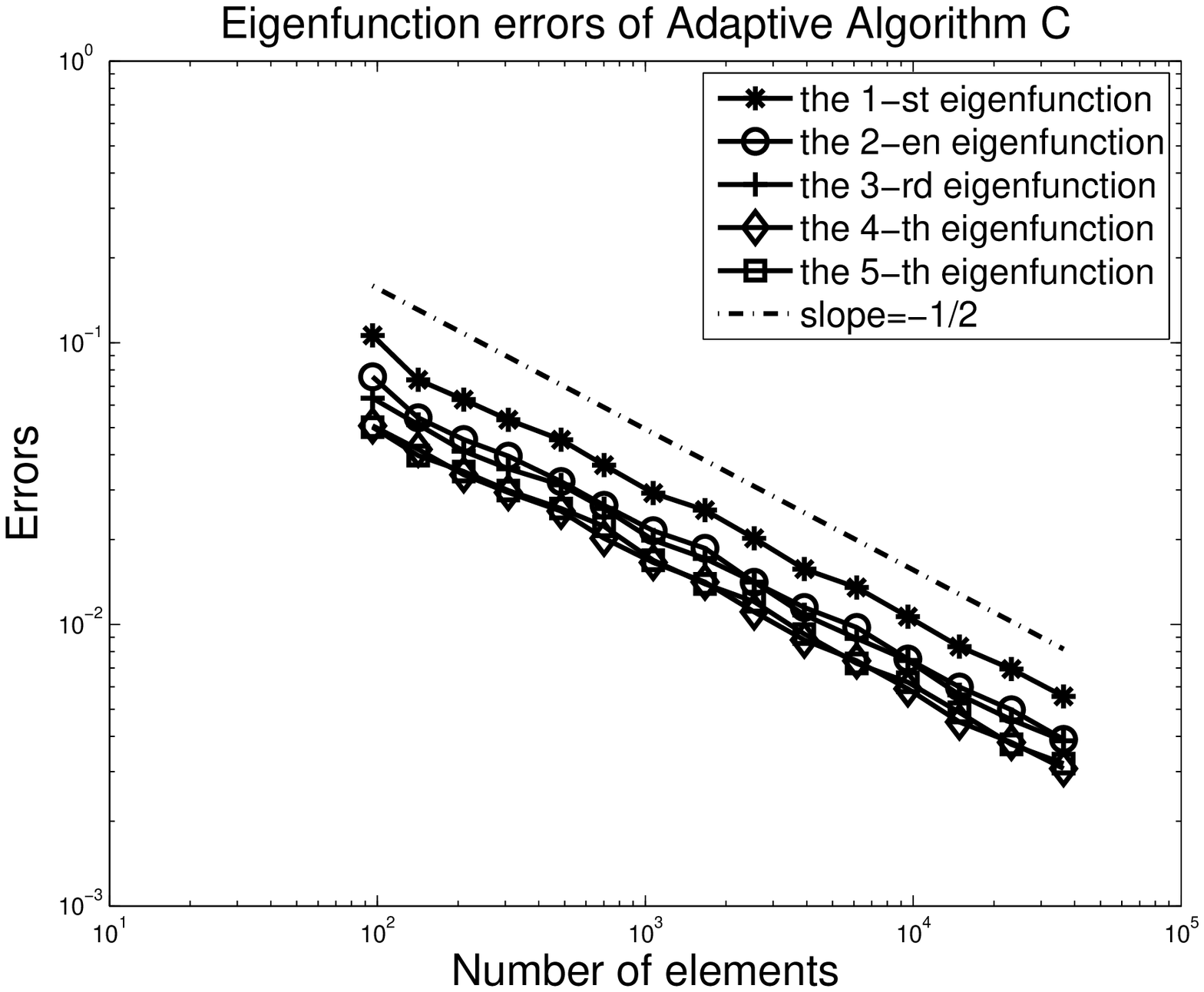}
\includegraphics[width=5.5cm,height=5.5cm]{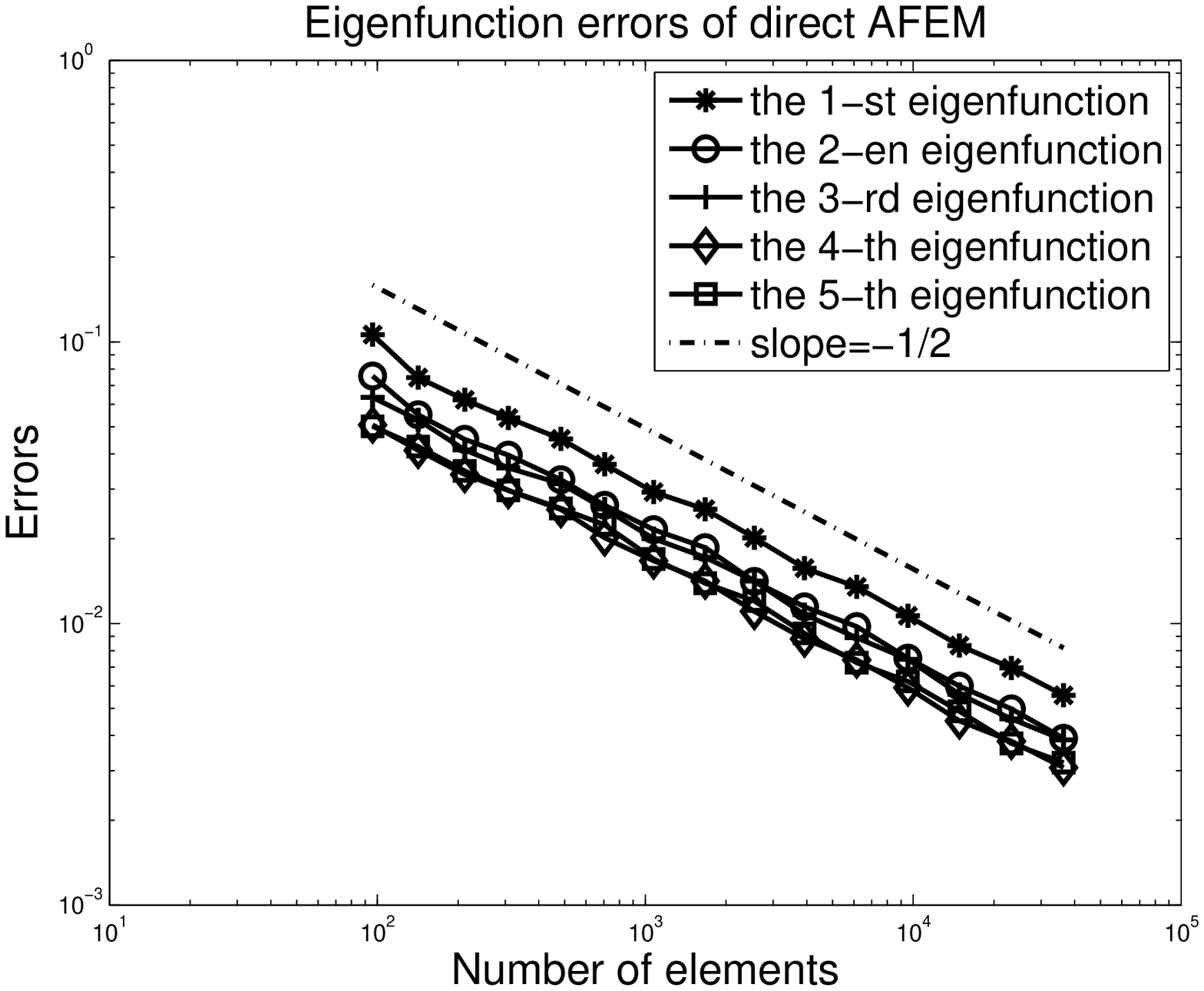}
\caption{The a posteriori error estimates of the eigenfunction approximations
by {\bf Adaptive Algorithm $C$}
and direct AFEM for Example 3 with the linear element} \label{Convergence_AFEM_Exam_3_5_Small}
\end{figure}

\begin{figure}[ht]
\centering
\includegraphics[width=5.5cm,height=5.5cm]{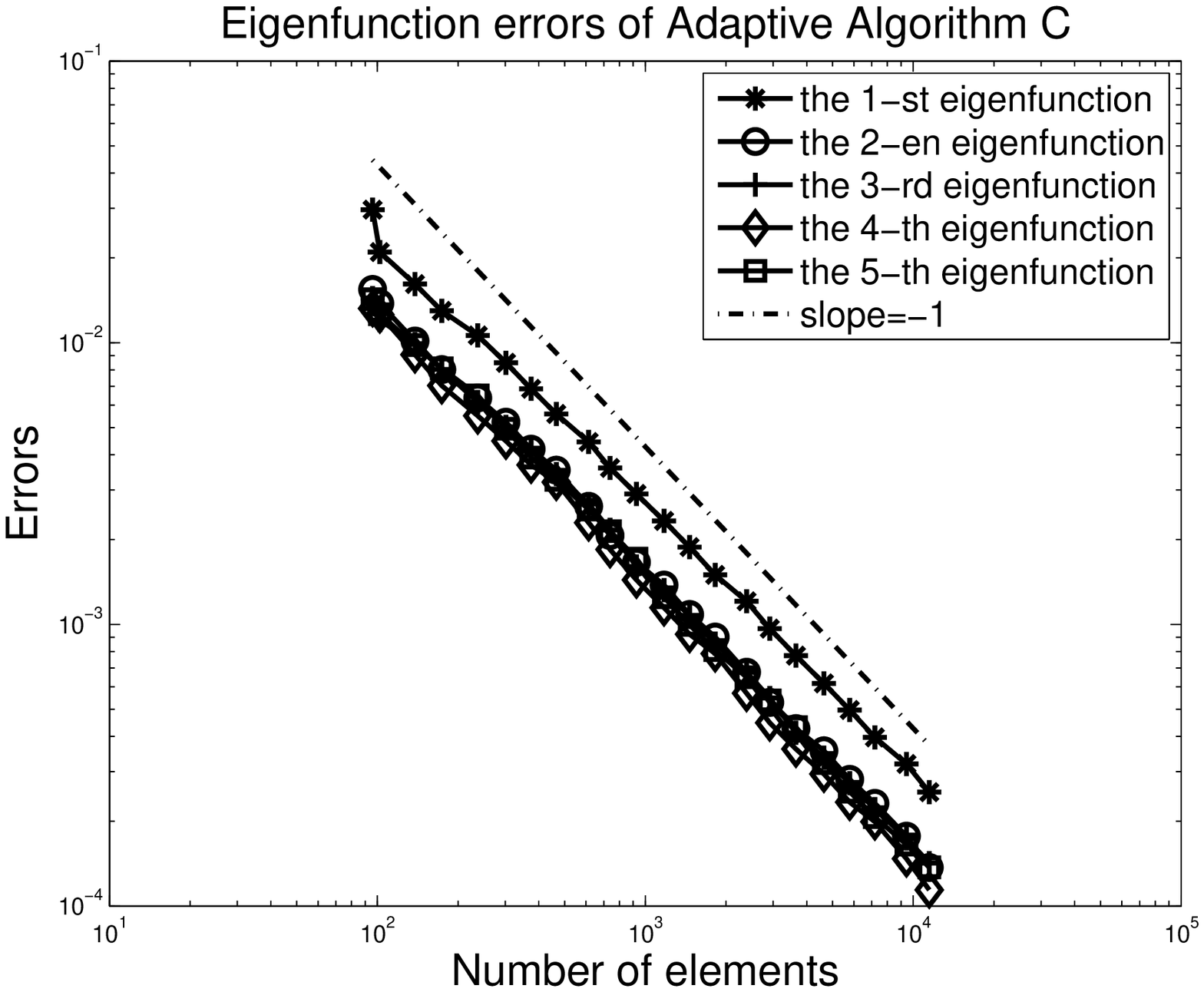}
\includegraphics[width=5.5cm,height=5.5cm]{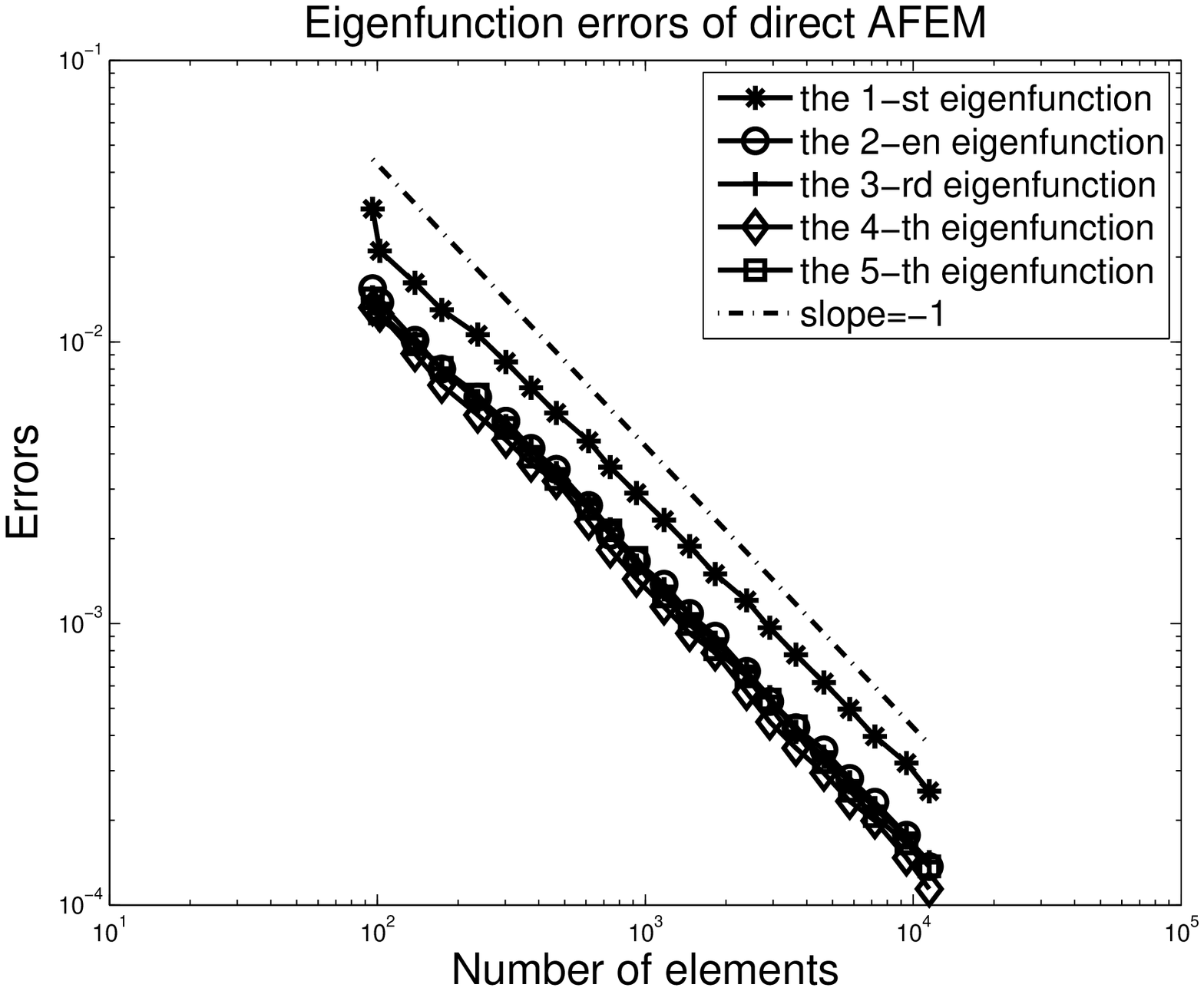}
\caption{The a posteriori error estimates of the eigenfunction approximations
by {\bf Adaptive Algorithm $C$}
and direct AFEM for Example 3 with the quadratic element} \label{Convergence_AFEM_Exam_3_5_Small_Quadratic}
\end{figure}

We also test {\bf Adaptive Algorithm $C$} for $5$ smallest eigenvalue approximations and
their associated eigenfunction approximations.  Figures \ref{Convergence_AFEM_Exam_3_5_Small}
and \ref{Convergence_AFEM_Exam_3_5_Small_Quadratic} show the a posteriori
error estimator produced by {\bf Adaptive Algorithm $C$} and direct AFEM with the linear and quadratic finite element methods, respectively.
In these cases, with {\bf Adaptive Algorithm $C$}, the small scale eigenvalue problem solving only appears
on the meshes with the number of elements: $[96,210,701,2549,9552,36402]$
(when $k= 1,3,6,9 ,12,15$ and $j_k=1,1,1,1,1,1$)
for linear element and $[ 96,237,738 ,2386,7185]$ (when $k=1,5,10,15,20$ and $j_k=1,1,1,1,1$)
for quadratic element, respectively.

From Figures \ref{Convergence_AFEM_Exam_3_First}, \ref{Convergence_AFEM_Exam_3_First_Quadratic},
\ref{Convergence_AFEM_Exam_3_5_Small} and \ref{Convergence_AFEM_Exam_3_5_Small_Quadratic},
we can find the approximations of eigenvalues as well as eigenfunctions by {\bf Adaptive Algorithm $C$}
have the same  convergence behavior as those by the direct AFEM which validates the efficiency of
our proposed scheme.

%=======================================================================================================
{\bf Example 4.}
In the last example, we consider the Laplace eigenvalue problem on three dimensional nonconvex domain
\begin{equation}\label{eigenproblem_Exam_2}
\left\{
\begin{array}{rcl}
-\Delta u &=&\lambda u\ \ \ \ {\rm in}\ \Omega,\\
u&=&0\ \ \ \ \ \ {\rm on}\ \partial\Omega,\\
\|u\|_{a,\Omega}&=&1,
\end{array}
\right.
\end{equation}
where $\Omega=(-1,1)^3\backslash[0, 1)^3$.
Similarly, eigenfunctions with singularities are expected due to the nonconvex property.

\begin{figure}[ht]
\centering
\includegraphics[width=5cm,height=5cm]{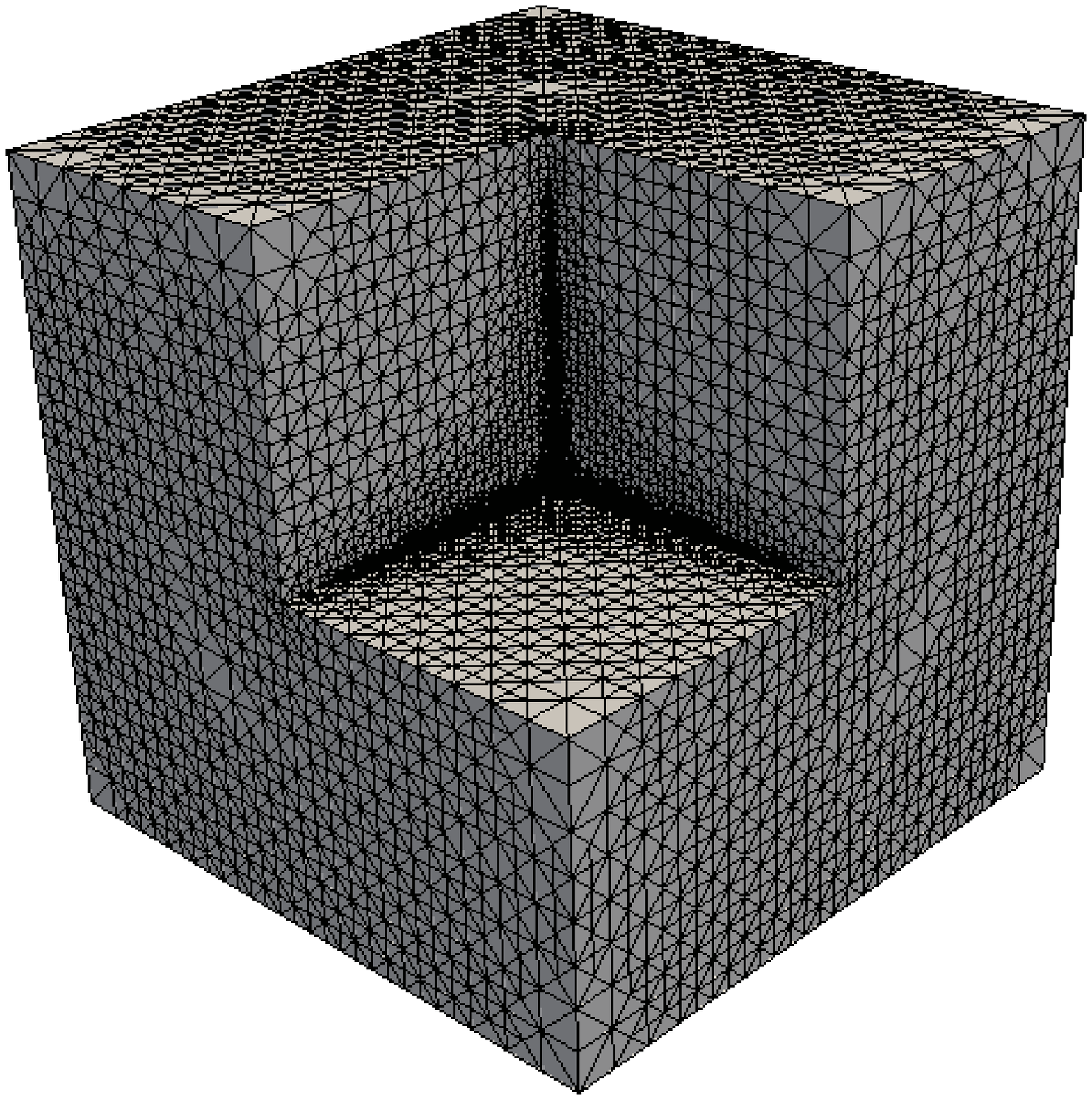} \ \ \ \ \ \ \ \ \
\includegraphics[width=5cm,height=5cm]{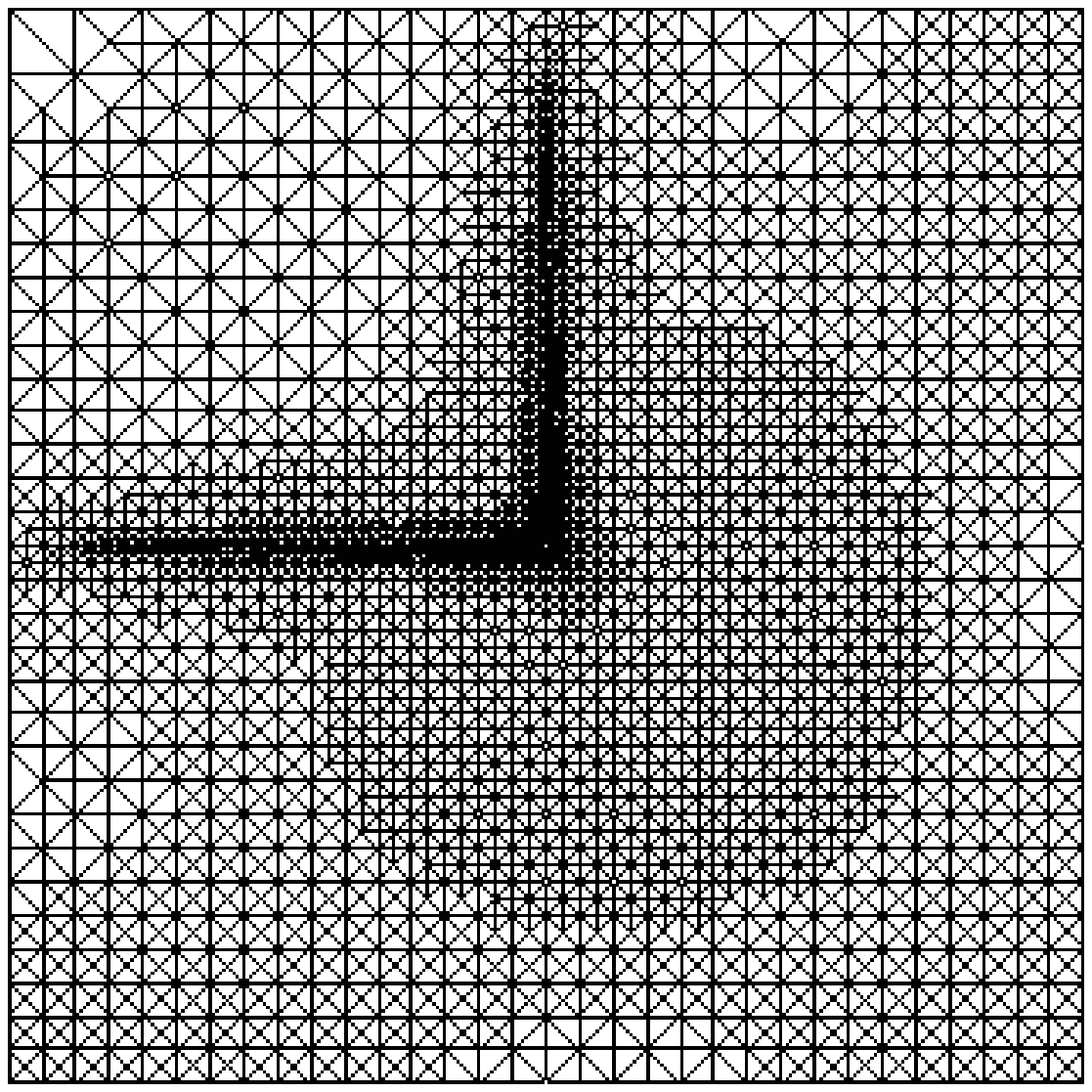}
\caption{The triangulations after adaptive iterations and section for
Example 4 by the linear element}\label{Convergence_AFEM_Exam_4_mesh1}
\end{figure}

\begin{figure}[ht]
\centering
\includegraphics[width=5.5cm,height=5cm]{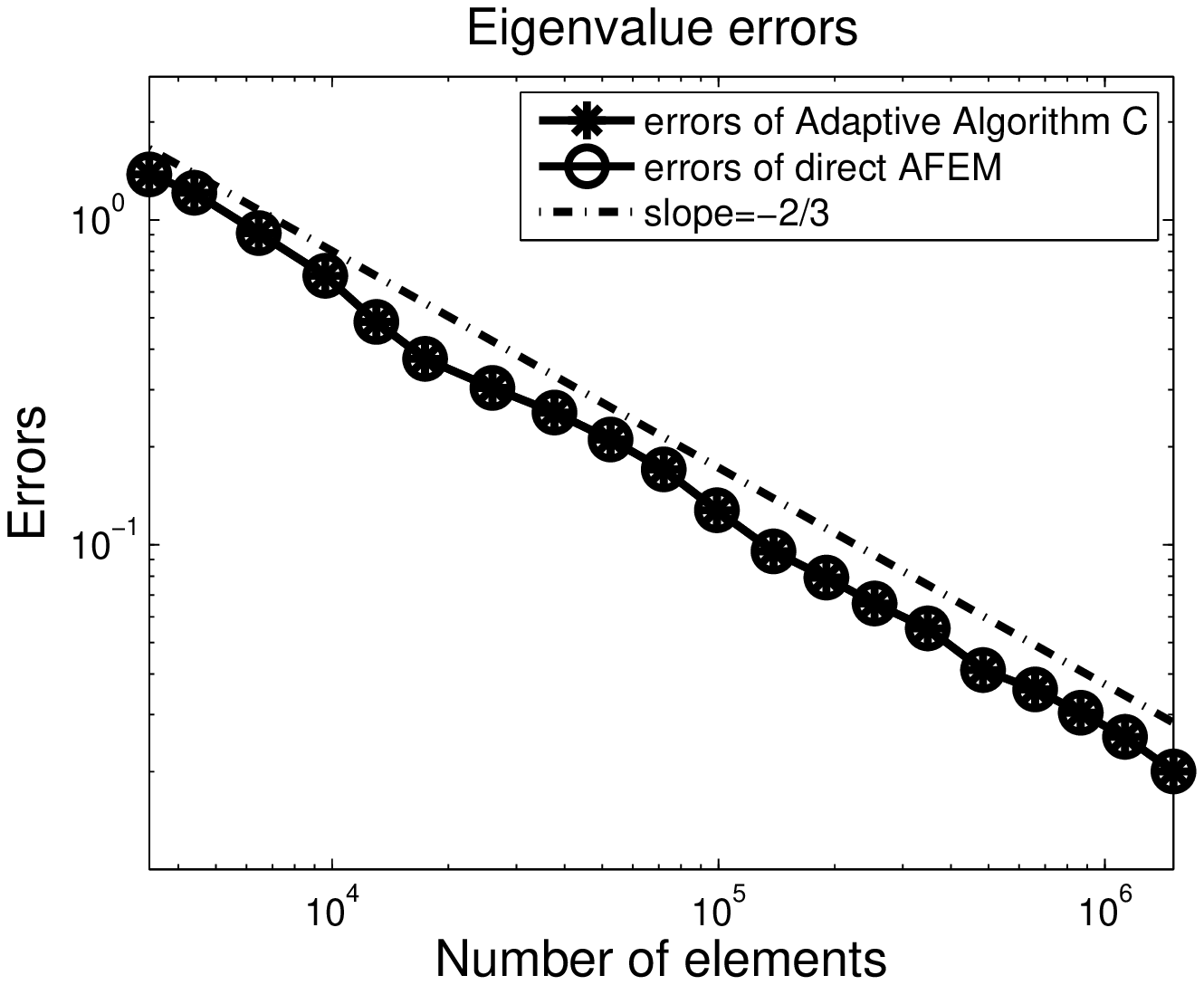}
\includegraphics[width=5.5cm,height=5cm]{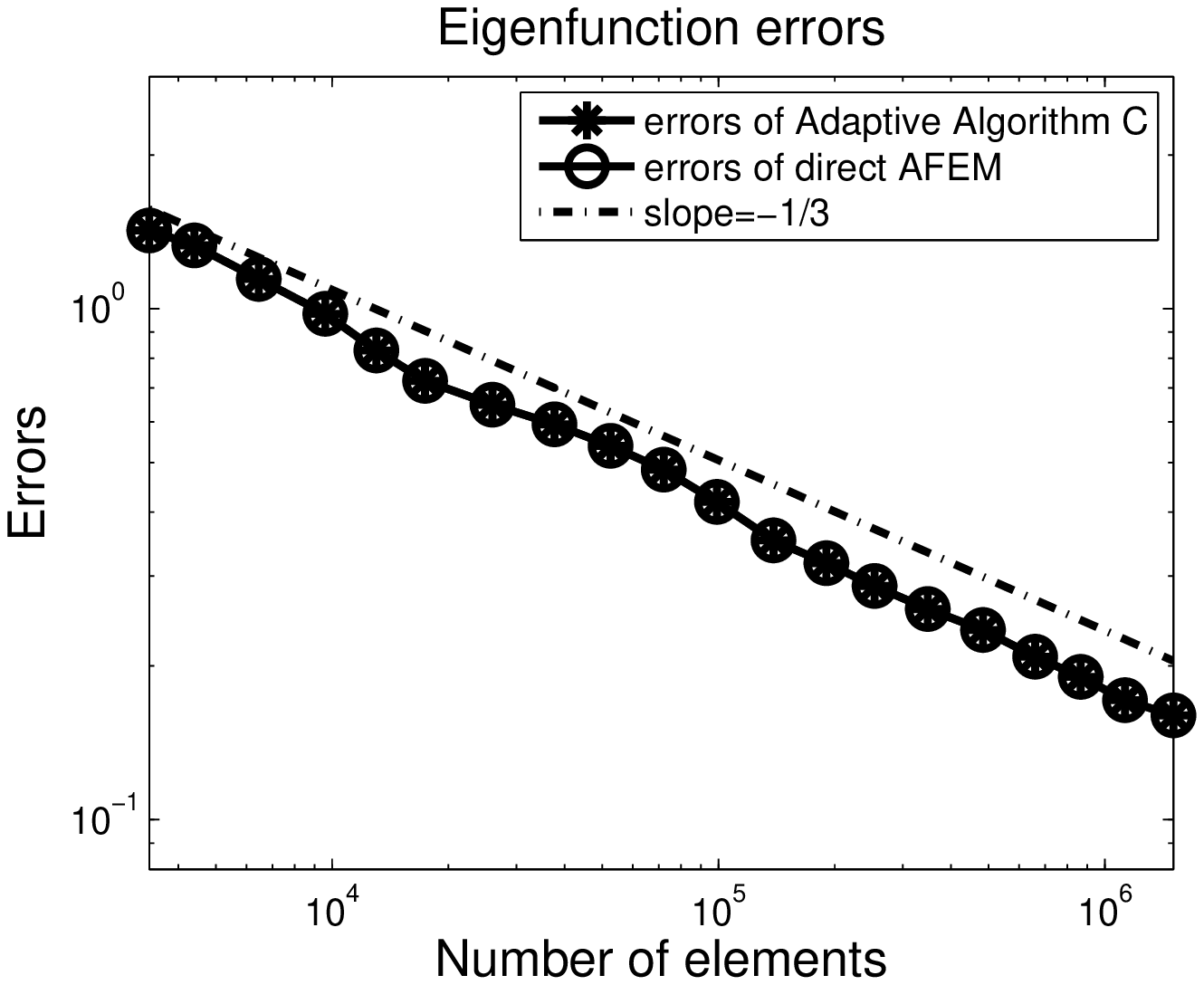}
\caption{The errors of the smallest eigenvalue and the associated
eigenfunction approximations by {\bf Adaptive Algorithm $C$}
and direct AFEM for Example 4 with the linear element}\label{Convergence_AFEM_Exam_4_P1}
\end{figure}

\begin{figure}[ht]
\centering
\includegraphics[width=5cm,height=5cm]{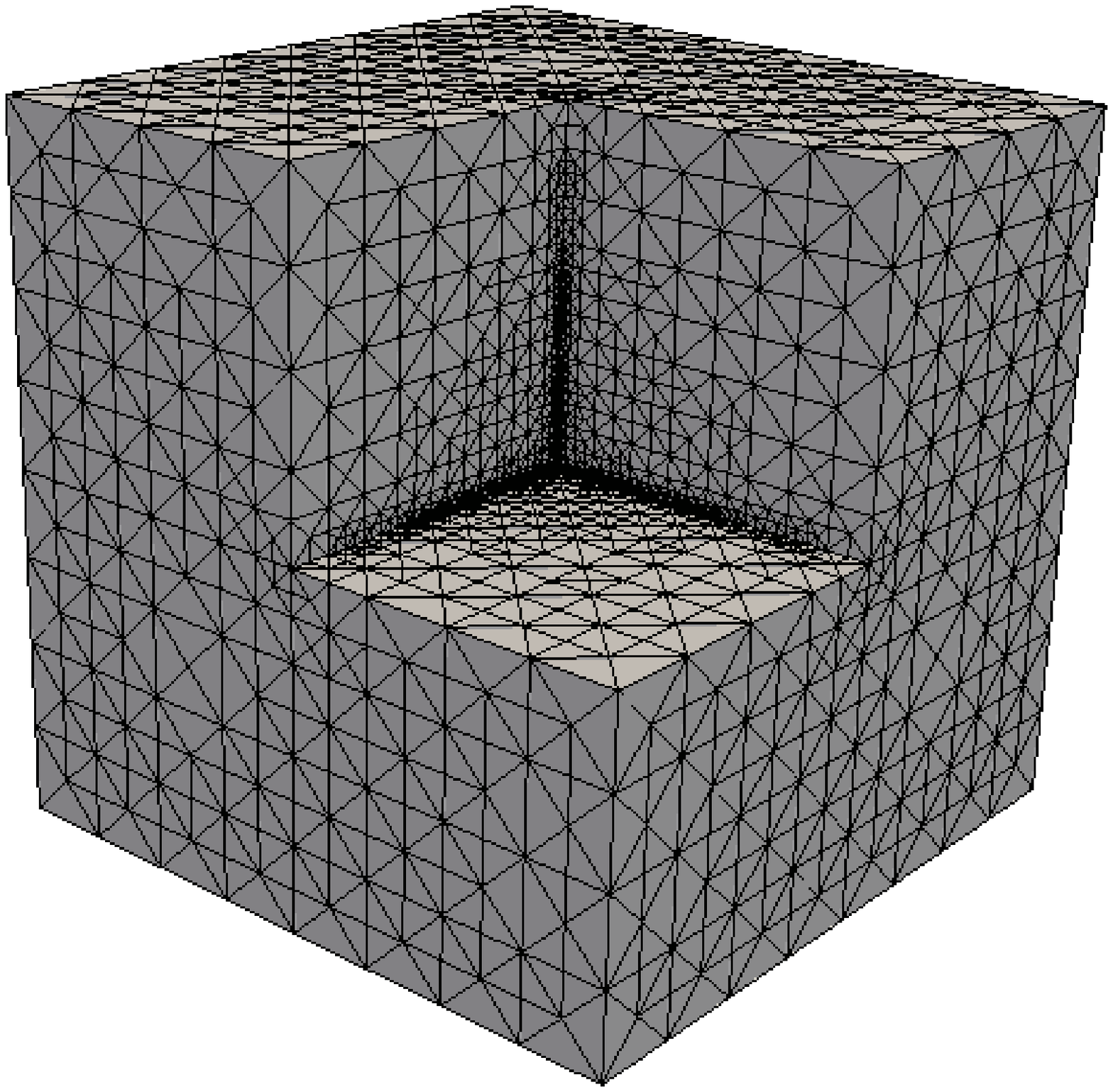}\ \ \ \ \ \ \ \
\includegraphics[width=5cm,height=5cm]{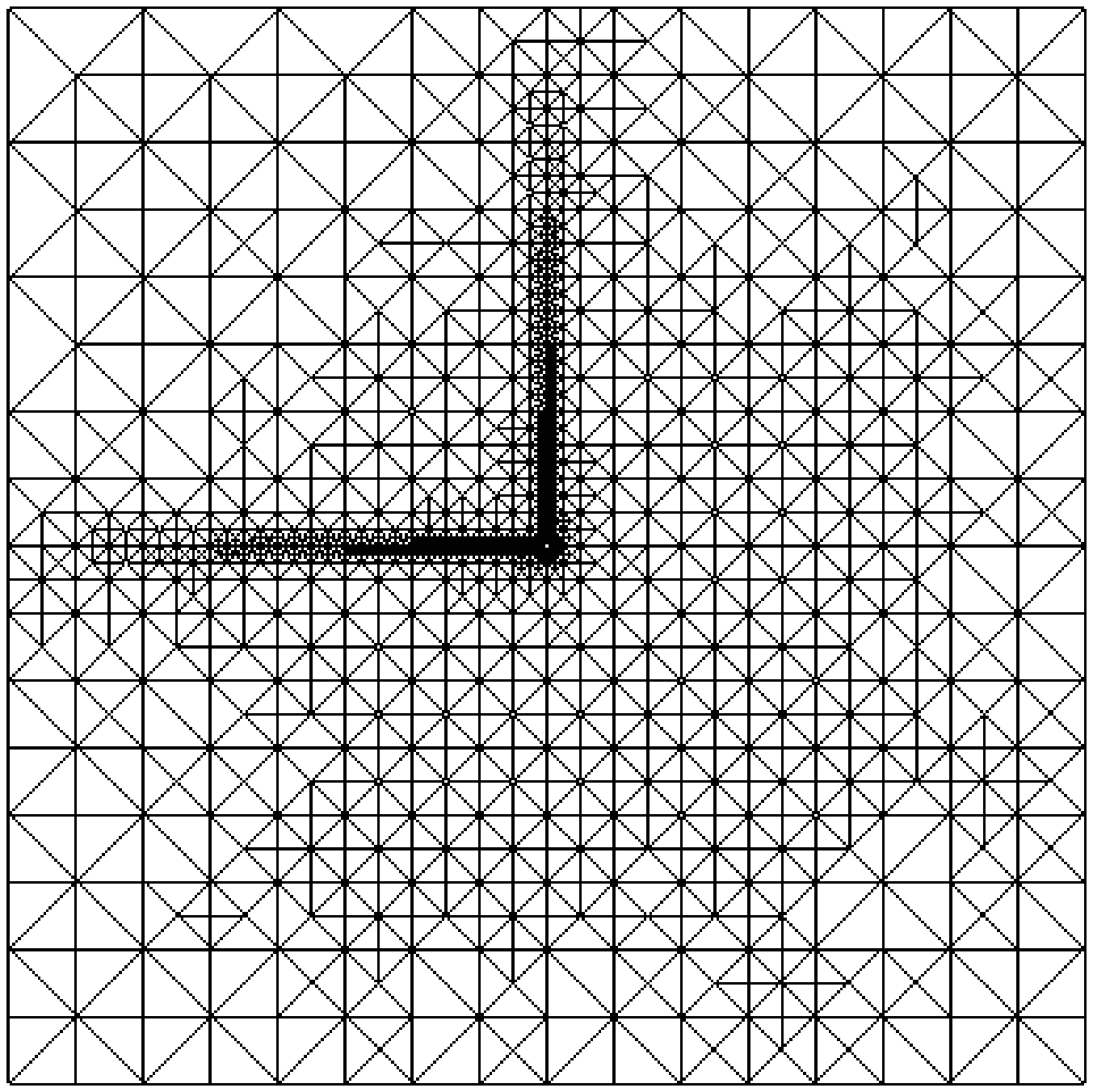}
\caption{The triangulations after adaptive iterations and section for
Example 4 by the quadratic element}\label{Convergence_AFEM_Exam_4_mesh2}
\end{figure}

\begin{figure}[ht]
\centering
\includegraphics[width=5.5cm,height=5cm]{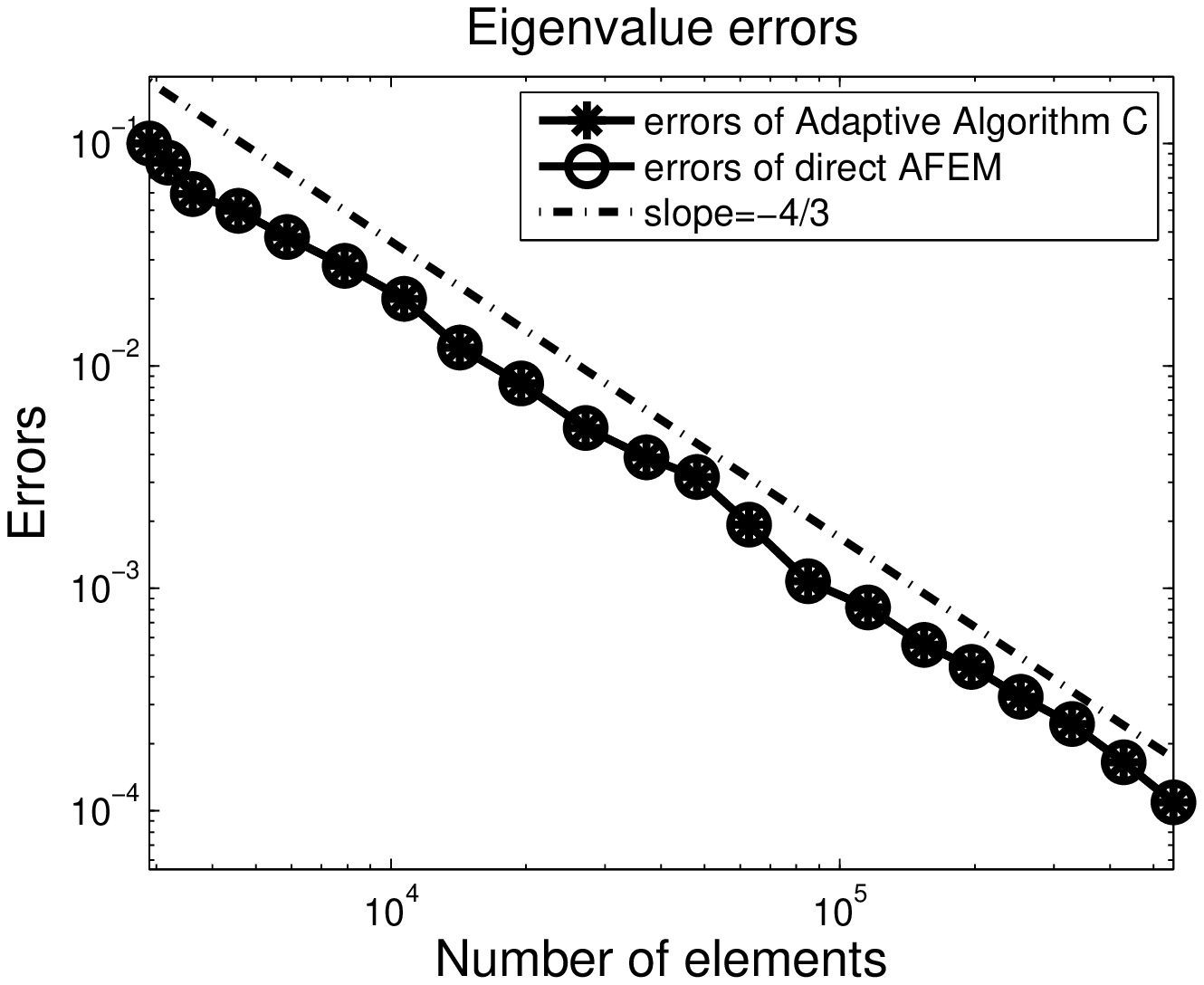}
\includegraphics[width=5.5cm,height=5cm]{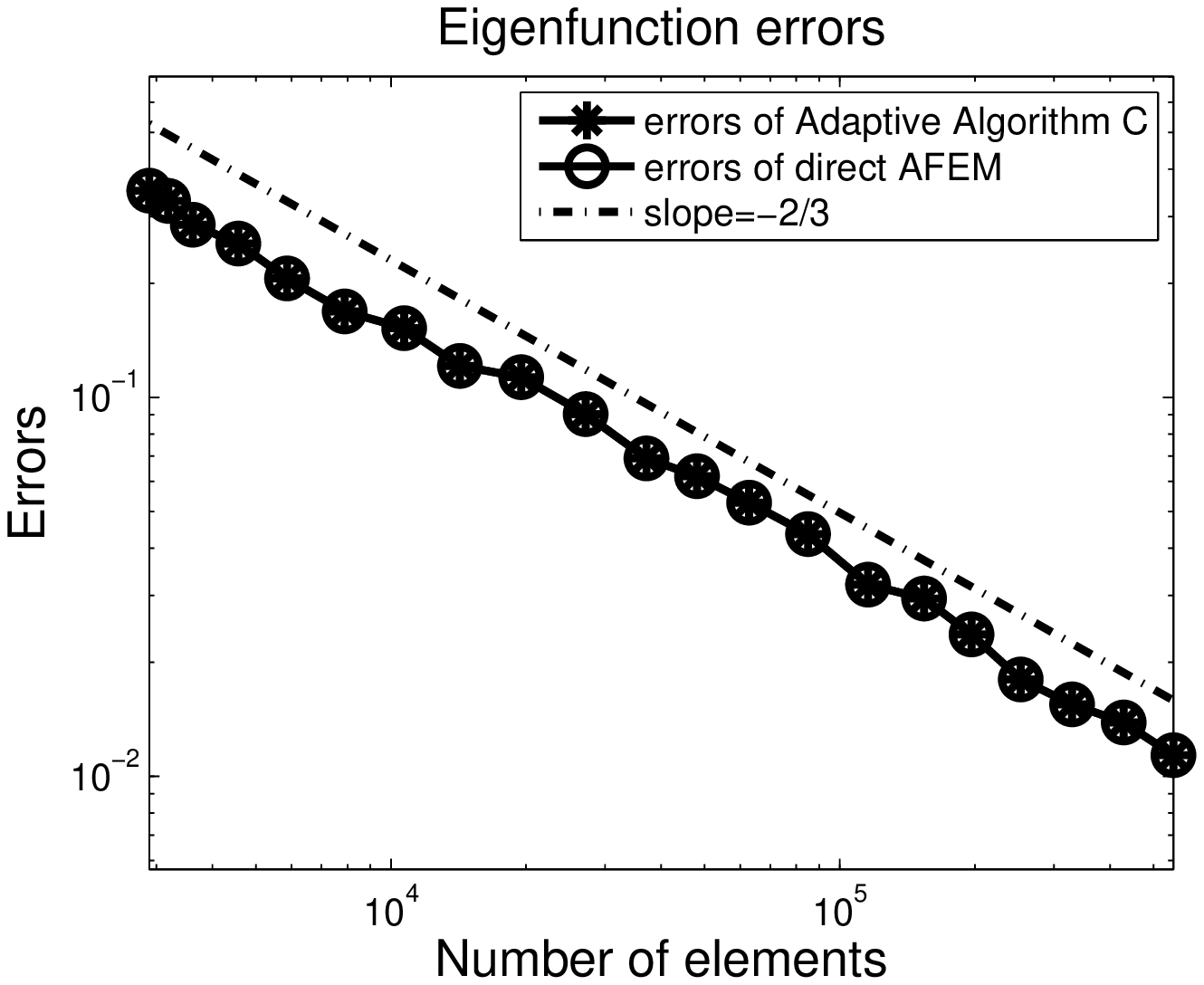}
\caption{The errors of the smallest eigenvalue and the associated
eigenfunction approximations by {\bf Adaptive Algorithm $C$}
and direct AFEM for Example 4 with the quadratic element}\label{Convergence_AFEM_Exam_4_P2}
\end{figure}

In this example, we give the numerical results of {\bf Adaptive Algorithm $C$}
with parameters $\theta_1 = 0.4$ and $\theta_2=0.7$ for linear element
and $\theta_1=0.4$ and $\theta_2=0.4$ for quadratic element, respectively.
First we investigate the numerical results for the first eigenvalue approximation.
Since the exact eigenvalue is not known, an adequately accurate approximation on finer finite element space
is chosen as the exact first eigenvalue for numerical tests.
Figure \ref{Convergence_AFEM_Exam_4_mesh1} and \ref{Convergence_AFEM_Exam_4_mesh2} show the triangulations after
adaptive iterations with the linear and quadratic finite element methods, respectively.
Figures \ref{Convergence_AFEM_Exam_4_P1} and \ref{Convergence_AFEM_Exam_4_P2}
give the corresponding numerical results.
In order to show the efficiency of {\bf Adaptive Algorithm $C$} more clearly, we  also compare
the results with those obtained by direct AFEM.
With {\bf Adaptive Algorithm $C$}, it is only required to solve small scale eigenvalue problems
in the low dimensional space $V_H+{\rm span}\{u_k\}$ when the the numbers of elements of the mesh are $[2688, 7634, 20586, 78102, 287442,1001202]$
($k=1, 5, 8, 12, 16, 20$ and \comm{$j_k=1,1,1,1,1,1$}) for linear element and
$[2688, 4020, 12144, 57188, 264576]$ ($k=1, 5, 9, 14, 19$ and $j_k=1,1,1,1,1$)
for quadratic element, respectively.

\begin{figure}[ht]
\centering
\includegraphics[width=5.5cm,height=5cm]{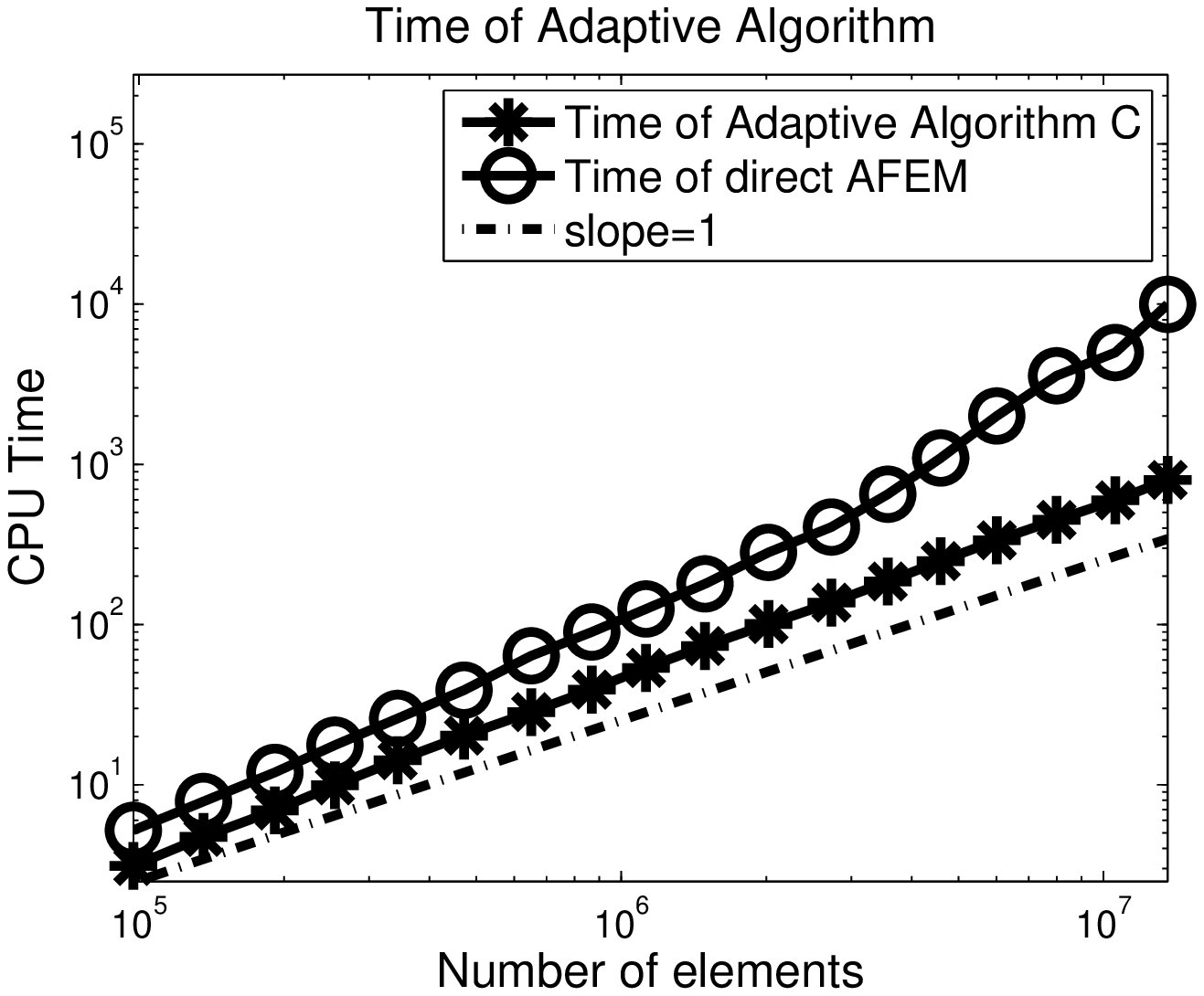}
\caption{The computational time by {\bf Adaptive Algorithm $C$}
and direct AFEM for Example 4 with the linear element}\label{time}
\end{figure}
\revise{In order to show the efficiency of {\bf Adaptive Algorithm $C$} more clearly, we compare the computational time (in second)
of {\bf Adaptive Algorithm $C$} with that of direct AFEM by linear element.
Figure \ref{time} shows the corresponding CPU time results,
which shows  {\bf Adaptive Algorithm $C$} has higher  efficiency than the direct AFEM.}

\begin{figure}[ht]
\centering
\includegraphics[width=5.5cm,height=5cm]{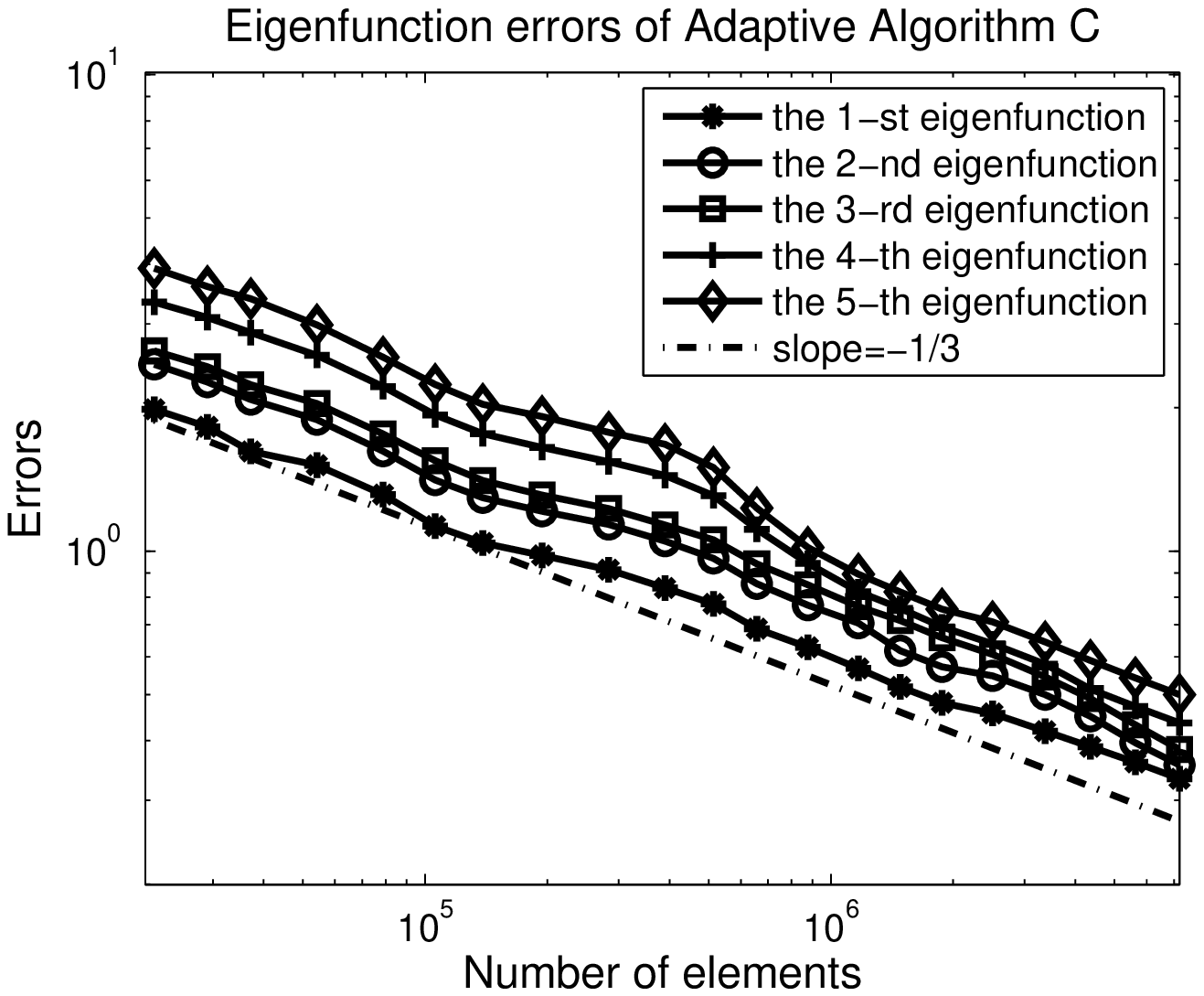}
\includegraphics[width=5.5cm,height=5cm]{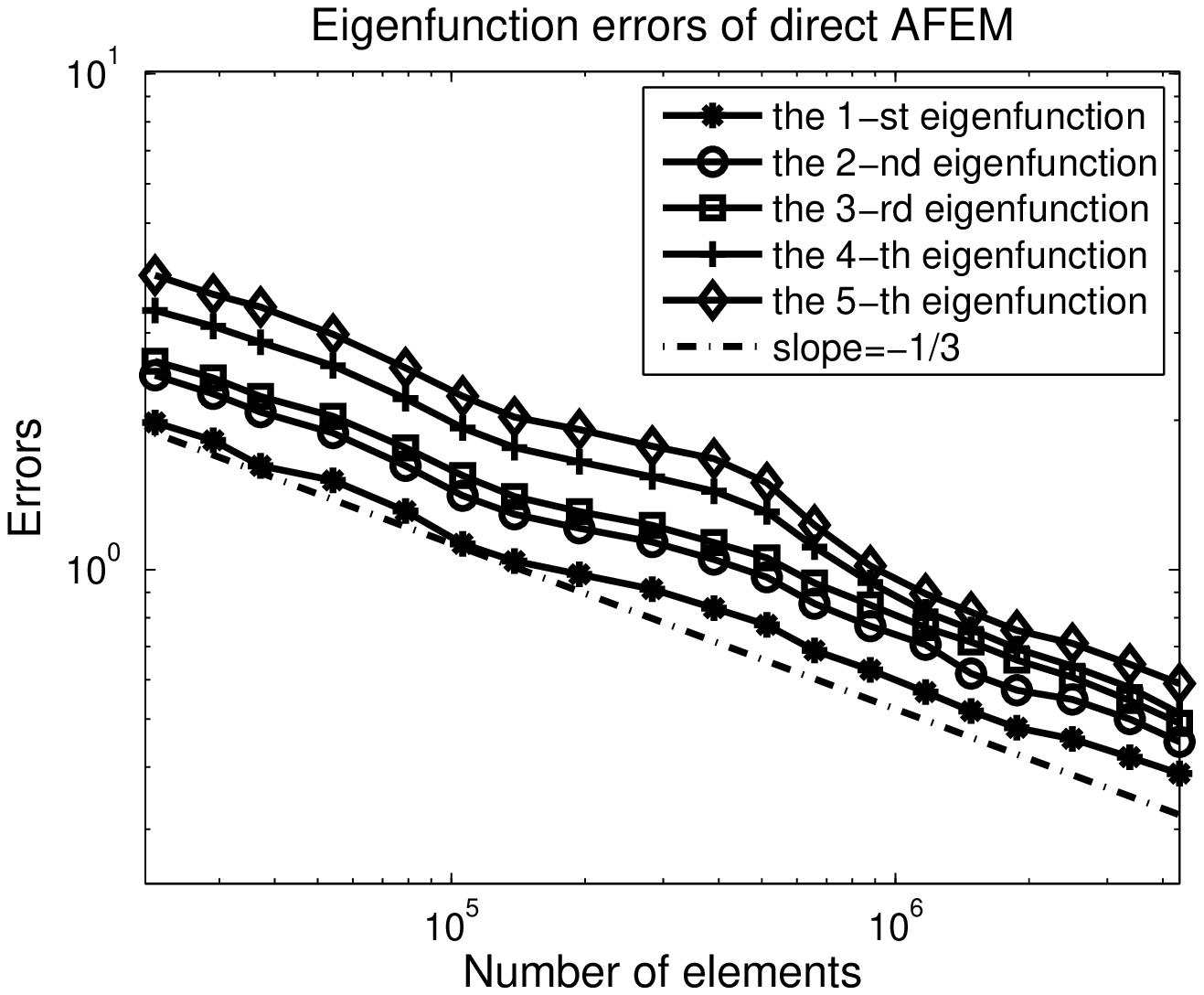}
\caption{The a posteriori error estimates of the eigenpair approximations
by {\bf Adaptive Algorithm $C$}
and direct AFEM for Example 4 with the linear element}\label{Convergence_MulAFEM_Exam_4}
\end{figure}
\revise{
Furthermore, we also test {\bf Adaptive Algorithm $C$} for the smallest $5$ eigenvalue and their associated
eigenfunctions. Figure \ref{Convergence_MulAFEM_Exam_4} shows the a posteriori
error estimators produced by {\bf Adaptive Algorithm $C$} and direct AFEM with the
linear finite element method. In these cases, {\bf Adaptive Algorithm $C$} only solve the small scale eigenvalue problems
on the low dimensional space $V_H+{\rm span}\{u_k\}$ when the the numbers of elements are
$[21504,  67846, 173182,  584308, 2218702]$ ($k=1, 5, 8, 12, 17$ and $j_k=1,1,1,1,1$).}

%-----------------------------------------------------------------------------------------------------
\section{Concluding remarks}
%-----------------------------------------------------------------------------------------------------
\revise{
In this paper, we present an efficient AFEM for eigenvalue problems based on multilevel correction scheme and
the adaptive refinement technique. The most important contribution of this new AFEM is that there is no eigenvalue
solving in the adaptively refined meshes which need much more computation than solving the corresponding linear boundary value problems.
Furthermore, the convergence and quasi-optimal complexity have also been proved based on a relation
between the eigenvalue problem and the associated boundary value problem
(see Theorems \ref{theorem_convergence_eigenpair} and (\ref{theorem_complexity})).
Some numerical experiments for both simple and multiple eigenvalue cases are provided to demonstrate the efficiency of
the proposed AFEM for eigenvalue problems.
It is obvious that this efficient AFEM can be extended to the nonlinear eigenvalue problems and also other type of
nonlinear problems which will be our future work.
}

%The convergence of the type of AFEM here for multiple eigenvalue case
% is different with the analysis for simple eigenvalue case since the exact eigenfunction which is approximated
% may be different between the successive two adaptive steps.  The analysis of the multiple eigenvalue case
% need stronger error estimates which will be one of our topics in the future.
%
%
%
%%-----------------------------------------------------------------------------------------------------

\end{document}